\documentclass{emsart} 
\usepackage{pb-diagram} 
\usepackage{makeidx} 
 \makeindex 
\contact[ikath@mis.mpg.de]{Ines Kath, Max Planck Institute for Mathematics in the Sciences,\\ 
Inselstra{\ss}e 22, D-04103 Leipzig, Germany} 
\contact[martin.olbrich@uni.lu]{Martin Olbrich, University of 
Luxembourg,
Mathematics Laboratory,\\
162 A, avenue de la Fa\"{\i}encerie,
L-1511 Luxembourg,
Grand-Duchy of Luxembourg} 
\newcommand\CC{{\mathbb C}} 
\newcommand\RR{{\mathbb R}} 
\newcommand\NN{{\mathbb N}} 
\newcommand\ZZ{{\mathbb Z}} 
\newcommand\HH{{\mathbb H}} 
 
\newcommand\fb{{\frak b}} 
\newcommand\fg{{\frak{g}}} 
\newcommand\fh{{\frak h}} 
\newcommand\fri{{\frak i}} 
\newcommand\fj{{\frak j}} 
\newcommand\fl{{\frak l}} 
 
\newcommand\fn{{\frak n}} 
\newcommand\fk{{\frak k}} 
\newcommand\fa{{\frak a}} 
\newcommand\fd{{\frak d}} 
\newcommand\fr{{\frak r}} 
\newcommand\fs{{\frak s}} 
 
\newcommand\fz{{\frak z}} 

\newcommand\cZ{{{\cal Z}_Q^2}} 
\newcommand\cH{{{\cal H}_Q^2}} 
\newcommand\cA{{\cal A}} 
\newcommand\cL{{\cal L}} 
\newcommand\HPhi{{\cH(\fl,\Phi_\fl,\fa)}} 
\newcommand\ZPhi{{\cZ(\fl,\Phi_\fl,\fa)}} 

\newcommand\eps{\varepsilon} 
\newcommand{\gl}{\mathop{{\frak g \frak l}}} 
\newcommand{\fsl}{\mathop{{\frak s \frak l}}} 
\newcommand{\fsp}{\frak s \frak p} 
\newcommand\osc{\frak o \frak s \frak c} 
 
\newcommand{\so}{\mathop{{\frak s \frak o}}} 
\newcommand{\su}{\mathop{{\frak s \frak u}}} 
\newcommand{\End}{\mathop{{\rm End}}} 
\newcommand{\Der}{\mathop{{\rm Der}}} 
\newcommand{\Dera}{\mathop{{\rm Der}_a}} 
\newcommand{\GL}{\mathop{{\it GL}}} 
\newcommand{\Hom}{\mathop{{\rm Hom}}} 
\newcommand{\Aut}{\mathop{{\rm Aut}}} 
\newcommand{\Id}{{{\rm id}}} 
\newcommand{\ad}{{{\rm ad}}} 
 
\newcommand{\Ad}{\mathop{{\rm Ad}}} 
\newcommand{\sgn}{\mathop{{\rm sign}}} 
 
\newcommand{\Ann}{\mathop{{\rm Ann}}} 
\newcommand{\Ker}{\mathop{{\rm ker}}} 
\newcommand{\im}{\mathop{{\rm im}}}

\renewcommand{\Im}{\mathop{{\rm Im}}}

\newcommand{\Span}{{{\rm span}}}   

\newcommand\ip{{\langle\cdot \,,\cdot \rangle}} 
\newcommand\lb{{[\cdot\,,\cdot]}} 
\newcommand\ipa{{\langle\cdot \,,\cdot \rangle_\fa}} 
\newcommand\dd{\fd_{\alpha,\gamma}(\fl,\Phi_\fl,\fa)} 
\newtheorem{theo}{Theorem}[section] 
\newtheorem{pr}{Proposition}[section] 
\newtheorem{de}[pr]{Definition}
\newtheorem{ex}[pr]{Example}
\newtheorem{re}[pr]{Remark}
\newtheorem{co}[pr]{Corollary}
\newtheorem{lm}[pr]{Lemma} 
\title{The classification problem for pseudo-Rie\-mannian symmetric 
spaces} 
\author[]{Ines Kath\thanks{Supported by Heisenberg program, DFG.} and Martin Olbrich}

\begin{document}   
\begin{abstract} 
    Riemannian and pseudo-Riemannian symmetric spaces with semisimple 
    trans\-vection group are known and classified for a long time. 
    Contrary to that the description of pseudo-Riemannian symmetric 
    spaces with non-semisimple transvection group is an 
    open problem. In the last years some progress on this problem was 
    achieved. In this article we want to explain these results and 
    some of their applications. 
\end{abstract}     
\maketitle
\tableofcontents 
\section{Introduction}  
There are many basic problems in differential geometry that are 
completely solved for Riemannian manifolds, but that become really 
complicated in the pseudo-Riemannian situation. One of these problems is 
the  
determination of all possible holonomy\index{holonomy} groups
of pseudo-Rie\-mannian 
manifolds. While holonomy groups of Riemannian manifolds are 
classified the problem is open for general pseudo-Riemannian 
manifolds, only the Lorentzian case is solved.  The difficulty is that in general
the holonomy representation of a pseudo-Riemannian manifold 
does not decompose into irreducible summands. Of course, we can 
decompose the representation into indecomposable ones, i.e., into subrepresentations that do 
not have proper non-degenerate invariant subspaces. By de Rham's 
theorem the indecomposable summands are again holonomy 
representations. This reduces the problem to the classification of 
indecomposable holonomy representations. 
Indecomposable holonomy represesentations that are not irreducible
have isotropic invariant subspaces. Such representation are especially difficult to handle 
if these invariant subspaces do not have an invariant 
complement. Manifolds that have an 
indecomposable holonomy representation are called indecomposable. 
Manifolds that are not indecomposable are at least locally a product 
of pseudo-Riemannian manifolds. Hence, we can speak of local factors 
of such a manifold. 

Many open questions in pseudo-Riemannian differential geometry are directly 
related to the unsolved holonomy problem. One of these open questions 
is the classification problem for symmetric spaces. Pseudo-Riemannian 
symmetric spaces\index{symmetric space!pseudo-Riemannian} are in some sense the most simple pseudo-Riemannian 
manifolds. Locally they are characterised by parallelity of the 
curvature tensor. As global manifolds they are defined as follows. A connected pseudo-Riemannian 
manifold $M$ is called a pseudo-Riemannian symmetric space if for any 
$x\in M$ there is an involutive isometry $\theta_{x}$ of $M$ that has 
$x$ as an isolated fixed point. 
In other words, for any $x\in M$ the geodesic reflection at $x$ 
extends to a globally defined isometry of $M$. 

The theory of Riemannian symmetric spaces was developed simultaneously 
with the  
theory of semisimple Lie groups and algebras by E. Cartan during the first 
decades of the twentieth century. It results in a complete 
classification of these spaces, see Helgason's beautiful 
book \cite{H} on the 
subject. The theory remains similar in spirit as long as one is interested in 
pseudo-Riemannian symmetric spaces whose holonomy representation 
is completely reducible as an algebraic representation, i.e., if any invariant subspace of the 
holonomy representation has an invariant complement. 
 These so called reductive 
symmetric spaces were classified by Berger \cite{B} in 1957. This 
classification is essentially 
the classification of 
involutions on real semisimple Lie algebras.

In order to understand non-reductive pseudo-Riemannian symmetric spaces one has to 
consider more general Lie algebras, which are, moreover, equipped with an invariant inner product.
Note that in contrast to the semisimple case this inner product is 
really an additional datum since it is not just a multiple of the Killing 
form. Such a pair $(\fg,\ip)$ consisting of a finite-dimensional real Lie 
algebra and an $\ad$-invariant non-degenerate symmetric bilinear form 
on it is called a {\it metric Lie algebra}\index{metric Lie algebra}. In the literature metric 
Lie algebras 
appear under various different names, e.g.~as quadratic\index{Lie 
algebra!quadratic} or orthogonal Lie algebras\index{Lie 
algebra!orthogonal}.

The transition from symmetric spaces to metric Lie algebras proceeds as follows.
Let $(M,g)$ be a pseudo-Riemannian symmetric space. The group $G$ 
generated by compositions of geodesic reflections $\theta_{x}\circ\theta_{y}$, $x,y\in M$,
is called the transvection group\index{transvection group} of $(M,g)$. It acts transitively on 
$M$. We fix a base point $x_{0}\in M$. The reflection $\theta_{x_{0}}$ 
induces an involutive automorphism $\theta$ of the Lie algebra $\fg$ of $G$, 
and thus a decomposition $\fg=\fg_{+}\oplus\fg_{-}$. The natural 
identification $T_{x_{0}}M\cong\fg_{-}$ induces an 
$\ad(\fg_{+})$-invariant bilinear form $\ip_{-}$ on $\fg_{-}$.
It is an important observation (see \cite{CP}) that $\ip_{-}$ extends 
uniquely to an
$\ad(\fg)$- and $\theta$-invariant inner product $\ip$ on $\fg$.
Thus, starting with a pseudo-Riemannian symmetric space $(M,g)$, we obtain a metric Lie 
algebra $(\fg,\ip)$ together with an isometric 
involution $\theta$ on it. The resulting triple $(\fg,\theta,\ip)$ will be 
called a {\it symmetric triple}\index{symmetric triple} later on.  Up to local isometry, $(M,g)$ can be recovered from 
this structure. We will call 
a symmetric space as well as its associated symmetric triple 
semisimple (reductive, solvable etc.) if its transvection group is semisimple (reductive, solvable etc.).   
The reader will find more details on the correspondence between 
symmetric spaces and symmetric triples in Section \ref{S31}.
For the general theory of symmetric spaces he may consult \cite{Ber, CP, H, KN, koh, L}.
The moral we want to stress at this point is that the understanding 
of metric Lie algebras is crucial for the understanding of symmetric 
spaces.

The present paper focuses on the classification problem for 
pseudo-Riemannian symmetric spaces. The above discussion reduces this 
problem  
to the classification of symmetric triples. It is easy to see that we 
can decompose every symmetric triple into a direct sum  of a semisimple 
one and one whose underlying Lie algebra does not have simple ideals. 
Of course, pseudo-Riemannian symmetric spaces that are associated with semisimple 
symmetric triples are reductive, and thus, as 
explained above, already classified. This is the reason for our 
decision to concentrate here on metric Lie algebras and symmetric triples without simple 
ideals. Thus the investigation of the geometry of semisimple 
symmetric spaces, which is still an active and interesting field, will be 
left almost untouched in this paper. 

The classification of metric Lie algebras appears to be very
difficult. Most likely, one has to accept that one won't get a 
classification in the sense of a list that includes all metric 
Lie algebras for arbitrary index of the inner product. The same is 
true for symmetric triples, the existence of an involution neither 
decreases nor really increases the difficulties. Therefore the aim 
is to develop a structure theory for
metric Lie algebras (and symmetric triples)
that allows a systematic construction and that gives a ``recipe'' how 
to get an explicit classification under suitable 
addiditional conditions, e.g., for small index of the inner product. 
In \cite{KO0,KO1,KO2} we developed a new strategy to reach this aim. The initial idea of this strategy 
is due to B\'erard-Bergery who 
observed that every symmetric triple without simple ideals arises in 
a canonical way by an extension procedure from ``simpler'' Lie algebras with involution. 
We used this idea to give a cohomological description of 
isomorphism classes of metric Lie algebras (and symmetric triples), 
which gives a suitable classification scheme. Here we will present 
this method, called quadratic extension\index{extension!quadratic}, and some of its applications. Moreover, we will survey 
earlier and related 
results due to Cahen, Parker, Medina, Revoy, Bordemann, Alekseevsky, Cort\'es, 
and others concerning metric Lie algebras and symmetric 
triples from this new point of view. 

We do not aim at a complete 
overwiew on the work on metric Lie algebras and symmetric spaces, for 
instance 
the basic results of Astrahancev (see e.g. \cite{A}) will not be 
discussed. However, we try to present a quite complete up-to-date 
account for classification results for metric Lie algebras (Section~\ref{S25}), symmetric spaces 
(Section~\ref{S33}), and symmetric spaces
with certain complex or quaternionic structures (Sections~\ref{S42} 
and \ref{S43}).

Note that metric Lie algebras are of interest in their own 
right, not only in the context of symmetric spaces. They naturally appear in 
various contexts, e.g., in Mathematical physics or in Poisson geometry. As an illustration, we 
shortly discuss the notions of Manin triples and pairs and present a 
new construction method for Manin pairs based on 
the theory of quadratic extensions of metric Lie algebras (Section~\ref{S52}). 
As a further application we study pseudo-Riemannian extrinsic symmetric spaces
by our method (Section~\ref{S51}).

Though being a survey article the paper also contains some new results.
A first group of new results appears in Section \ref{S2} and is due to 
the fact that we develop here a unified theory which works for metric 
Lie algebras, symmetric triples, and symmetric triples with additional 
structures at once. Most of these results are straightforward 
generalisations of the corresponding special results given in the 
original papers \cite{KO1,KO2,KO3}. Proofs that really require new
ideas will be given in the appendix. This generality makes 
Sections~\ref{S22} and \ref{S23} a little bit more technical than 
usual for survey articles. However, having mastered these moderate 
technical difficulties the reader will see in the subsequent sections 
how quite different results
follow easily from one general principle.
The results in Section~\ref{S32} concerning the geometric meaning
of the quadratic extension procedure and the above mentioned 
construction method for Manin pairs appear here for the first time.
In addition, we announce a new result on the structure of 
hyper-K\"ahler symmetric spaces (Theorem~\ref{grob}).

The theory of metric Lie algebras and pseudo-Riemannian symmetric
spaces is far from being complete. In fact, there is a huge amount of 
open problems. The difficulty is to find those questions, which 
really lead to new theoretical 
insight and not just to messy calculations. 
We hope that the questions raised at several places in this paper belong 
to the fist category.

\subsection*{Some conventions} 

We denote by $\NN$ the set of positive integers and we put 
$\NN_{0}:=\NN\cup \{0\}$. 

Let $(\fa,\ip_{\fa})$ be a pseudo-Euclidean vector space. A subspace 
$\fa'\subset \fa$ is called isotropic if $\ip_{\fa}|_{\fa'}=0$.  A 
basis $A_{1},\ldots,A_{p},A_{p+1},\ldots,A_{p+q}$ of $\fa$ is called 
orthonormal if $A_{i}\perp A_{j}$ for $i\not=j$, 
$\langle A_{i},A_{i} \rangle_{\fa}=-1$ for $i=1,\ldots,p$ and $\langle 
A_{j},A_{j} \rangle_{\fa}=1$ for $j=p+1,\ldots,p+q$. In this case 
$(p,q)$ is called signature and $p$ index of $\ip_{\fa}$ (or of $\fa$).
Let $\ip_{p,q}$ be the inner product of signature 
$(p,q)$ on $\RR^{p+q}$  for which the 
standard basis of $\RR^{p+q}$ is an orthonormal basis. Then we call
$\RR^{p,q}:=(\RR^{p+q},\ip_{p,q})$ standard pseudo-Euclidean space. 

We will often  describe a Lie algebra by 
giving a 
basis and some of the Lie bracket relations, e.g. we will write the 
three-dimensional Heisenberg algebra as 
$\fh(1)=\{[X,Y]=Z\}$. In this case we always assume that 
all other brackets of basis vectors vanish. If we do not mention the 
basis explicitly, then we suppose that all basis vectors appear in one 
of the bracket relations (on the left or the right hand side). 

Let $\fl$ be a Lie algebra and let $\fa$ be an $\fl$-module. Then 
$\fa^{\fl}$ denotes the space of invariants in $\fa$, i.e., 
$\fa^{\fl}=\{A\in\fa\mid L(A)=0\} \mbox{ for all } L\in\fl$.
\section{Metric Lie algebras} \label{S2} 
\subsection{Examples of metric Lie algebras}\label{S21} 

The easiest example of a metric Lie algebra is an abelian Lie algebra 
together with an arbitrary (non-degenerate) inner product.  A further 
well-known example is a semisimple Lie algebra equipped with a 
non-zero multiple of its Killing form. 

Let $H$ be a Lie group and $\fh$ its Lie algebra. 
The cotangent bundle $T^{*}H$ can be given a group structure such 
that the associated Lie algebra equals $\fh\,_{\ad^{*}}\ltimes \fh^{*}$. 
Now let $\ip_{\fh}$ be any invariant symmetric bilinear form on 
$\fh$ (which can degenerate). 
We can define on $\fh\,_{\ad^{*}}\ltimes \fh^{*}$ a symmetric 
bilinear form $\ip$ by adding $\ip_\fh$ to the dual pairing of $\fh$ and 
$\fh^{*}$, that is by 
$$\langle H_{1}+Z_{1},H_{2}+Z_{2}\rangle = Z_{1}(H_{2}) + 
Z_{2}(H_{1})+ \langle H_{1},H_{2}\rangle_{\fh}$$ 
for $H_1,H_2\in\fh$, $Z_{1},Z_{2}\in\fh^{*}$. 
It is not hard to prove that $\ip$ is invariant and non-degenerate, 
its signature equals $(\dim \fh,\dim \fh)$. 
Hence $(\fh\,_{\ad^{*}}\ltimes \fh^{*},\ip)$ is a metric Lie algebra. 
In particular, $\ip$ induces a biinvariant metric on $T^{*}H$. 

The following construction is a generalisation of the previous example. 
It is due to Medina and Revoy \cite{MR1}. Starting with 
an $n$-dimensional metric Lie algebra and an arbitrary 
$m$-dimensional  Lie algebra it produces a metric Lie algebra of 
dimension $n+2m$. Let $(\fg,\ip_{\fg})$ be a metric Lie algebra and let 
$(\fh,\ip_{\fh})$ be a Lie algebra with an invariant symmetric 
bilinear form (which can degenerate). Furthermore, let $\pi:\ 
\fh\rightarrow \Dera(\fg,\ip_\fg)$ be a Lie algebra homomorphism 
from $\fh$ into the Lie algebra of all antisymmetric derivations of 
$\fg$. We denote by $\beta\in 
C^{2}(\fg,\fh^{*}):=\Hom(\bigwedge^{2}\fg,\fh^{*})$  the 2-cocycle 
(see \ref{S23} for this notion) 
$$\beta(X,Y)(H):=\langle \pi(H)X,Y \rangle_{\fg},\quad X,Y\in\fg, 
H\in\fh.$$ 
On the vector space $\fd:=\fh^{*}\oplus\fg\oplus\fh$ we define a Lie bracket $\lb$ 
by 
\begin{eqnarray*} 
    \lefteqn{[(Z,X,H),(\tilde Z,\tilde X,\tilde H)]=}\\ 
    && (\beta(X,\tilde X)+\ad_\fh^{*}(H)\tilde Z-\ad_\fh^{*}(\tilde H) 
    Z,\ [X,\tilde X]_{\fg}+\pi(H) \tilde X -\pi(\tilde H)X,\ [H,\tilde 
    H]_{\fh}) 
\end{eqnarray*} 
and an inner product $\ip$  by 
$$\langle (Z,X,H),(\tilde Z,\tilde X,\tilde H) \rangle = 
    \langle X,\tilde X\rangle_{\fg }+\langle 
    H,\tilde H\rangle_{\fh }+Z(\tilde H)+\tilde Z(H)$$ 
for all $Z,\tilde Z\in\fh^{*}$, $X,\tilde X\in\fg$ and $H,\tilde 
H\in\fh$. Then $\fd_{\pi}(\fg,\fh):=(\fd,\ip)$ is a metric Lie 
algebra. It is called double extension\index{extension!double} of $\fg$ by $\fh$ since it can be regarded 
as 
an extension of 
the semi-direct product $\fg\rtimes_{\pi}\fh$ by the abelian Lie 
algebra $\fh^{*}$. If the signature of $\fg$ equals $(p,q)$ and if 
$\dim \fh=m$, then the signature of  $\fd_{\pi}(\fg,\fh)$ equals 
$(p+m,q+m)$. 

The importance of this construction becomes clear from 
the following structure theorem by Medina and Revoy. 
It says that 
we can inductively produce all metric Lie algebras from simple and 
one-dimensional ones by taking direct sums and applying the double 
extension procedure. 
\begin{theo}[Medina/Revoy  \cite{MR1}] \label{TMR} 
If $(\fg,\ip)$ is an indecomposable 
metric Lie algebra, then either $\fg$ is simple or $\fg$ is 
one-dimensional or $\fg$ is a double extension 
$\fd_{\pi}(\tilde\fg,\fh)$ of a metric Lie algebra $\tilde \fg$ by a 
one-dimensional or a simple Lie algebra $\fh$. 
\end{theo}     

We remark that for the special case of solvable metric Lie 
algebras this result can already be found in the form of exercises in
Kac's book \cite{Kac}, Exercises 2.10,11. For a complete proof in this 
case see also \cite{FS}. 

Using Theorem~\ref{TMR} it is not hard to see that any 
indecomposable  non-simple Lorentzian metric Lie algebra\index{metric 
Lie algebra!Lorentzian} 
is the double extension of an abelian Euclidean metric Lie algebra by 
a one-dimensional Lie algebra. This allows the classification of 
isomorphism classes of Lorentzian metric Lie algebras \cite{M}, 
compare Example~\ref{osc} and Theorem~\ref{lorentz}. 

In principle one can try and use this method to classify also metric 
Lie algebras of higher index. This was done in \cite{BK} for index two. 
However, now the following difficulty arises. In general a metric Lie 
algebra of index greater than one can be obtained in many different 
ways by double extension from a lower-dimensional one. Thus in order 
to solve the classification problem we would have to decide under 
which conditions two metric Lie algebras arising in different ways by 
repeated application of the double extension construction (and direct sums) are 
isomorphic. This seems to be very complicated. Therefore we are now looking 
for a way that avoids this difficulty. In the following we will 
develop a structure theory for metric Lie algebras which is more 
adapted to classification problems. The basic idea of this theory goes back to 
B\'erard-Bergery \cite{BB2} who suggested to 
consider indecomposable non-semisimple pseudo-Riemannian symmetric 
spaces as the result of two subsequent extensions. Our starting point is the 
following construction. 

Let $\fl$ be a Lie algebra and let $\rho:\fl\rightarrow \so(\fa)$ be 
an orthogonal representation of $\fl$ on a pseudo-Euclidean vector 
space $(\fa,\ip_{\fa})$. Take $\alpha\in 
C^{2}(\fl,\fa):=\Hom(\bigwedge^{2}\fl,\fa)$ and $\gamma\in 
C^{3}(\fl):=\Hom(\bigwedge^{3}\fl,\RR)$.  We consider the vector space 
$\fd:=\fl^*\oplus\fa\oplus \fl$ 
and define an inner product $\ip$ by 
$$ 
\langle Z_1+A_1+L_1, Z_2+A_2+L_2\rangle :=   
\langle A_1,A_2\rangle _\fa +Z_1(L_2)+Z_2(L_1) 
$$ 
for $Z_1, Z_2\in \fl^*,\, A_1,A_2\in\fa$ and $L_1,L_2\in\fl$. Moreover, 
we define an antisymmetric bilinear map $\lb:\fd\times \fd\rightarrow\fd$ by 
$[\fl^*,\fl^*\oplus\fa]=0$ and 
\begin{eqnarray*} 
\,[L_1,L_2]&=&\gamma(L_1,L_2,\cdot)+\alpha(L_1,L_2)+[L_1,L_2]_\fl\\ 
\,[L,A] &=& -\langle A,\alpha(L,\cdot)\rangle +L(A)\\ 
\,[A_1,A_2] &=& \langle\rho(\cdot) A_1,A_2\rangle \\ 
\,[L,Z] &=& \ad^*(L)(Z) 
\end{eqnarray*} 
for $L,L_1,L_2\in\fl$, $A,A_1,A_2\in\fa$ and $Z\in\fl^*$. 
Then the Jacobi identity for $\lb$ is equivalent to a certain cocycle 
condition for $\alpha$ and $\gamma$. We will denote this condition by 
$(\alpha,\gamma)\in\cZ(\fl,\fa)$ and postpone its exact formulation to 
Section~\ref{S23}. Thus, if $(\alpha,\gamma)\in\cZ(\fl,\fa)$, then 
$\lb$ is a Lie bracket and it is easy to check that $\ip$ is invariant 
with respect to this bracket. This gives the following result. 

\begin{pr}[\rm\cite{KO1}, Prop.~2.4] 
    If $(\alpha,\gamma)\in\cZ(\fl,\fa)$, then 
    $\fd_{\alpha,\gamma}(\fl,\fa):=(\fd,\ip)$ is a metric Lie algebra. 
\end{pr}     
Two special cases of this construction were known 
previously. 
For the case $\alpha=\gamma=0$ our $\fd_{\alpha,\gamma}(\fl,\fa)$ is 
the double extension $\fd_{\rho}(\fa,\fl)$ of the abelian metric Lie 
algebra $\fa$ by $\fl$ (in the sense of Medina and Revoy as explained 
above) and for $\fa=0$ it coincides with the 
$T^{*}$-extension\index{extension!$T^{*}$-} 
introduced by Bordemann \cite{Bor}. 

\begin{ex}\label{osc}{\rm 
    Take $\fl=\RR$. Let $\fa$ be the standard Euclidean vector 
    space $\RR^{2m}$ with (orthonormal) standard basis 
    $e_{1},\ldots,e_{2m}$. Take 
    $\lambda=(\lambda^{1},\ldots,\lambda^{m})\in (\fl^{*})^{m}\cong 
    \RR^{m}$. We define an 
orthogonal representation of $\fl$ 
    on $\fa$ by 
    $$\rho_{\lambda}(L)(e_{2i-1})=\lambda_{i}(L)\cdot e_{2i},\quad 
    \rho_{\lambda}(L)(e_{2i})=-\lambda_{i}(L)\cdot e_{2i-1}$$ 
    for $L\in\fl$ and $i=1,\ldots,m$. We set 
    $\fa_{\lambda}:=(\rho_{\lambda},\fa)$. Then 
    $\osc(\lambda):=\fd_{0,0}(\RR,\fa_{\lambda})$ is a 
    metric Lie algebra of signature $(1,2m+1)$. This Lie algebra 
    is often called oscillator algebra\index{oscillator algebra}. As explained above, $\osc(\lambda)$ 
    can also be considered as double extension of $\fa$ by $\RR$. }
\end{ex}   
In the following we will see that any metric Lie algebra without 
simple ideals is isomorphic to some $\fd_{\alpha,\gamma}(\fl,\fa)$ for 
suitable data $\fl$, $\fa$, $(\alpha,\gamma)\in\cZ(\fl,\fa)$ and how 
this fact can be used to describe isomorphism classes of metric Lie 
algebras. 

\subsection{Metric Lie algebras and quadratic extensions}\label{S22} 

As already mentioned in the introduction we are especially interested in metric Lie 
algebras without simple ideals. In this section we will learn more 
about the structure of such metric Lie algebras. Later on we wish to equip 
metric Lie algebras with additional structures, e.g.~with an 
involution when we want to study symmetric triples or with even more 
structure when we will be studying geometric structures on symmetric spaces. 
For this reason we develop a theory that is equivariant under a Lie 
algebra $\fh$ and a Lie group $K$ acting semisimply on $\fh$ by automorphisms. 
We assume throughout the paper that $K$ has only finitely many 
connected components. 
We suggest to take the trivial case $(\fh,K)=(0,\{e\})$ on first 
reading. In particular, this means that you may omit all maps called $\Phi$ 
in the following.

An 
$(\fh,K)$-module $(V,\Phi_{V})$ consists of a finite-dimensional 
vector space $V$ and a map $\Phi_{V}: \fh\cup K \rightarrow\Hom(V)$ 
such that $\Phi_{V}|_{\fh}: \fh\rightarrow \Hom (V)$ and 
$\Phi_{V}|_{K}: K\rightarrow \GL(V)\subset \Hom(V)$ are representations of $\fh$ 
and 
$K$ satisfying 
$$\Phi_V(k\cdot X) =\Phi_V(k)\Phi_V(X)\Phi_V(k)^{-1}$$ 
for all $k\in K$ and $X\in\fh$. There is a natural notion of a 
$(\fh,K)$-submodule. An $(\fh,K)$-module is called semisimple if for any 
$(\fh,K)$-submodule there is a complementary $(\fh,K)$-submodule. 

\begin{de}\begin{enumerate} 
    \item An $(\fh,K)$-equivariant Lie algebra $(\fl,\Phi_{\fl})$ is a Lie algebra 
$\fl$ 
    that is equipped with the structure of a semisimple $(\fh,K)$-module such that 
$\im(\Phi_{\fl}|_{\fh})\subset \Der (\fl)$ and 
 $\im(\Phi_{\fl}|_{K} )\subset \Aut (\fl)$, where $ \Der (\fl)$ and $\Aut(\fl)$ 
denote the Lie algebra of derivations and the group of automorphisms 
of $\fl$, respectively. 
\item  An $(\fh,K)$-equivariant metric Lie algebra\index{metric Lie 
algebra!$({\frak h},K)$-equivariant} is a triple $(\fg,\Phi,\ip)$ 
such that 
\begin{itemize} 
\item[(i)] $(\fg,\ip)$ is a metric Lie algebra, 
\item[(ii)] $(\fg,\Phi)$ is an $(\fh,K)$-equivariant Lie algebra, 
\item[(iii)] $\Phi(\fh)\subset \Dera(\fg)$ and  $\Phi(K)\subset\Aut(\fg,\ip)$, where $\Dera 
(\fg)$ denotes the Lie algebra of antisymmetric derivations on $(\fg,\ip)$. 
\end{itemize} 
The index (signature) of an $(\fh,K)$-equivariant metric Lie algebra 
$(\fg,\Phi,\ip)$ is the index (signature) of $\ip$.
Sometimes we abbreviate $(\fg,\Phi,\ip)$ as $\fg$. 
\end{enumerate} 
\end{de} 
A homomorphism (resp., isomorphism) of  $(\fh,K)$-equivariant Lie algebras 
$F:(\fl_1,\Phi_{1})\rightarrow 
(\fl_2,\Phi_{2})$ is a 
homomorphism (resp., isomorphism) of Lie algebras  $F:\fl_1\rightarrow 
\fl_{2}$ that satisfies $F\circ \Phi_{1}(h)=\Phi_{2}(h)\circ F$ for 
all $h\in\fh\cup K$. Isomorphisms of $(\fh,K)$-equivariant 
metric Lie algebras are in addition compatible with the given inner 
products. We have a natural notion of direct sums of 
$(\fh,K)$-equivariant (metric) Lie algebras. An $(\fh,K)$-equivariant 
(metric) Lie algebra is called decomposable if it is isomorphic to 
the direct sum of two non-trivial $(\fh,K)$-equivariant (metric) Lie 
algebras. Otherwise it is called indecomposable. 

In the following let $(\fl,\Phi_\fl)$ always be an $(\fh,K)$-equivariant Lie 
algebra. 
\begin{de} 
\begin{enumerate} 
    \item An $(\fl,\Phi_\fl)$-module 
$(\rho,\fa,\Phi_\fa)$  consists of 
\begin{itemize} 
\item[(i)] a semisimple $(\fh,K)$-module $(\fa,\Phi_{\fa})$, 
\item[(ii)] a representation $\rho:\fl\rightarrow \Hom(\fa,\ip)$ that 
satisfies 
$$\rho(\Phi_\fl(k)L)=\Phi_\fa(k)\circ\rho(L)\circ\Phi_\fa(k)^{-1},\quad 
\rho(\Phi_\fl(X)L)=[\Phi_\fa(X),\rho(L)]$$ 
for all $k\in K$, $X\in\fh$ and $L\in\fl$. 
\end{itemize} 
\item An orthogonal 
$(\fl,\Phi_\fl)$-module 
$(\rho,\fa,\ip_\fa,\Phi_\fa)$ consists 
of an  $(\fl,\Phi_\fl)$-module \linebreak $(\rho,\fa,\Phi_\fa)$ and an inner 
product $\ip_{\fa}$ on $\fa$ such that 
\begin{itemize} 
\item[(i)] $(\fa,\ip_{\fa},\Phi_\fa)$ is an abelian $(\fh,K)$-equivariant metric 
Lie algebra, 
\item[(ii)] $\rho$ is an orthogonal representation, i.e. $\rho:\fl\rightarrow 
\so(\fa,\ip)$.   
\end{itemize} 
\end{enumerate} 
We often abbreviate $(\rho,\fa,\Phi_\fa)$ and $(\rho,\fa,\ip_\fa,\Phi_\fa)$ 
as $(\rho,\fa)$ or  $\fa$. 
\end{de} 
Let $(\fl_{i},\Phi_{\fl_{i}})$, $i=1,2$, be two $(\fh,K)$-equivariant Lie 
algebras and let $(\rho_{i},\fa_{i})$, $i=1,2$, be orthogonal 
$(\fl_{i},\Phi_{\fl_{i}})$-modules. Let $S:\fl_1\rightarrow \fl_{2}$ 
be a homomorphism of $(\fh,K)$-equivariant Lie algebras and let 
$U:\fa_{2}\rightarrow \fa_{1}$ be an $(\fh,K)$-equivariant isometric embedding. 
Suppose that 
$$U\circ\rho_{2}(S(L))=\rho_{1}(L)\circ U$$ holds for all $L\in\fl$. 
Then we call $(S,U)$ morphism of pairs. We will write this as 
$(S,U):((\fl_{1},\Phi_{\fl_{1}}),\fa_{1})   \rightarrow 
((\fl_{2},\Phi_{\fl_{2}}),\fa_{2})$, but remember that $S$ and $U$ map 
in different directions. 

We will say that an ideal of an $(\fh,K)$-equivariant Lie algebra 
$(\fg,\Phi)$ is $\Phi$-in\-variant if it is invariant 
under all maps belonging to $\im\Phi$. 
\begin{de}\label{qex} 
Let $(\rho,\fa,\ip_\fa,\Phi_\fa)$ be an orthogonal $(\fl,\Phi_{\fl})$-module. 
A quadratic extension\index{extension!quadratic} of $(\fl,\Phi_{\fl})$ by $\fa$ is given by a 
quadruple $(\fg,\fri,i,p)$, where 
\begin{itemize} 
    \item[(i)] $\fg$ is an $(\fh,K)$-equivariant metric Lie algebra, 
    \item[(ii)] $\fri$ is an isotropic $\Phi$-invariant ideal of $\fg$, 
    \item[(iii)] $i$ and $p$ are homomorphisms of $(\fh,K)$-equivariant 
    Lie algebras constituting an exact sequence 
    $$0\longrightarrow \fa \stackrel i\longrightarrow \fg/\fri 
    \stackrel p\longrightarrow \fl\longrightarrow 0$$ 
    that is consistent with the representation $\rho$ of $\fl$ on $\fa$ 
    and has the property that $i$ is an isometry from $\fa$ to 
    $\fri^{\perp}/\fri$. 
\end{itemize} 
\end{de}     
\begin{ex}[The standard model] \label{sm} {\rm
    First we consider a Lie algebra $\fl$ without further structure, 
    i.e. $\Phi_{\fl}=0$. Let $\fa$ be an orthogonal $\fl$-module and 
    take $(\alpha,\gamma)\in\cZ(\fl,\fa)$. Let 
    $\fd_{\alpha,\gamma}(\fl,\fa)=(\fd,\ip)$ be the metric Lie algebra 
    constructed in Section~\ref{S21}. We identify $\fd/\fl^{*}$ with $\fa\oplus 
\fl$ 
and 
    denote by $i:\ \fa\rightarrow\fa\oplus\fl$ the injection and by 
    $p:\ \fa\oplus\fl \rightarrow\fl$ the projection. Then 
    $(\fd_{\alpha,\gamma}(\fl,\fa),\fl^*,i,p)$ is a quadratic extension of 
    $\fl$ by $\fa$. 
     
    Now suppose that we have in addition a $(\fh,K)$-structure 
    on $\fl$ and $\fa$, i.e.~let $(\fl,\Phi_\fl)$ be an $(\fh,K)$-equivariant Lie 
algebra 
    and let $(\rho,\fa,\ip_\fa,\Phi_\fa)$ be an orthogonal 
    $(\fl,\Phi_{\fl})$-module. Then we can define a 
    map  $\Phi: \fh\cup K \rightarrow\Der(\fd) \cup\Aut(\fd)$ by 
    \begin{eqnarray*} 
        \Phi(X)(Z+A+L) &=& 
-\Phi_{\fl}(X)^{*}(Z)+\Phi_{\fa}(X)(A)+\Phi_{\fl}(X)(L)\\ 
        \Phi(k)(Z+A+L) 
&=&(\Phi_{\fl}(k)^{*})^{-1}(Z)+\Phi_{\fa}(k)(A)+\Phi_{\fl}(k)(L). 
    \end{eqnarray*} 
    Then $\fd_{\alpha,\gamma}(\fl,\Phi_{\fl},\fa):=(\fd,\Phi,\ip)$ is an 
$(\fh,K)$-equivariant metric Lie 
    algebra if $(\alpha,\gamma)$ satisfies a certain 
    natural invariance condition with respect to $\Phi_{\fl}$ and 
    $\Phi_{\fa}$. We write $(\alpha,\gamma)\in\cZ(\fl,\Phi_{\fl},\fa)$ for this 
condition 
    whose exact formulation we will give 
    in Section~\ref{S23}. Hence, if 
    $(\alpha,\gamma)\in\cZ(\fl,\Phi_{\fl},\fa)$, then 
    $(\fd_{\alpha,\gamma},\fl^*,i,p)$ is a quadratic extension of 
    $(\fl,\Phi_{\fl})$ by $\fa$. It is called standard model since, as 
    we will see, any quadratic extension of $(\fl,\Phi_\fl)$ by $\fa$ 
    is in a certain sense equivalent to some 
    $(\fd_{\alpha,\gamma},\fl^*,i,p)$ for a suitable cocycle 
$(\alpha,\gamma)\in\cZ(\fl,\Phi_{\fl},\fa)$. }
\end{ex} 

What makes the theory of quadratic extensions so useful is the fact 
that any $(\fh,K)$-equivariant metric Lie algebra without simple ideals admits such 
a 
structure. Essentially, 
this follows from B\'erard-Bergery's investigations of 
pseudo-Riemannian holonomy representations and symmetric spaces 
in \cite{BB1,BB2,BB3}. He proved that for any metric Lie algebra 
$(\fg,\ip)$ there exists an isotropic ideal $\fri(\fg)\subset\fg$ such that 
$\fri(\fg)^{\perp}/\fri(\fg)$ is abelian. We want to 
describe the construction of this ideal now. However, instead of following 
\cite{BB1,BB2} we will give a description that is more adapted to the 
structure theory that we wish to develop here. In particular, we will 
give an $(\fh,K)$-equivariant formulation. 

Let $(\fg,\Phi,\ip)$ be an $(\fh,K)$-equivariant metric Lie algebra. There 
is a  chain of $\Phi$-invariant ideals 
\begin{align*} 
 \fg = R_{0}(\fg)\supset R_{1}(\fg)\supset R_{2}(\fg)\supset \ 
    \ldots \supset \ R_{l_{-}}(\fg)=0 
\end{align*}   
which is defined by the condition that $R_{k}(\fg)$ is the smallest ideal of 
$\fg$ contained in $R_{k-1}(\fg)$ such that the $\fg$-module 
$R_{k-1}/R_{k}(\fg)$ is semisimple. The ideal $R(\fg):=R_{1}(\fg)$ 
is called nilpotent radical of $\fg$. It has to be distinguished from 
the nilradical (i.e. the maximal nilpotent ideal) $\fn$ and the 
(solvable) radical $\fr$. By Lie's Theorem 
\begin{equation}\label{blub} 
R(\fg)=\fr\cap\fg'=[\fr,\fg]\subset\fn 
\end{equation} 
and $R(\fg)$ acts trivially on any semisimple $\fg$-module \cite{Bou}. 
We define an ideal $\fri(\fg)\subset\fg$ by 
$$ \fri(\fg):= \sum_{k=1}^{l} R_k(\fg) \cap  R_k(\fg)^{\perp} $$ 
and call it the canonical isotropic ideal\index{canonical isotropic 
ideal} of $\fg$. 

\begin{pr}[\cite{BB1,BB2}; \cite{KO1}, Lemma~3.4] 
    If $(\fg,\Phi,\ip)$ is an $(\fh,K)$-equivariant  metric Lie 
    algebra,  then $\fri(\fg)$ is a $\Phi$-invariant isotropic 
    ideal and the $\fg$-module $\fri(\fg)^{\perp}/\fri(\fg)$  is semisimple. 
    If $\fg$ does not contain simple ideals, then the Lie 
    algebra $\fri(\fg)^{\perp}/\fri(\fg)$ is abelian. 
\end{pr}     
In particular, $\fg/\fri(\fg)^{\perp}$ becomes an $(\fh,K)$-equivariant 
Lie algebra  and $\fri(\fg)^{\perp}/\fri(\fg)$ a semisimple orthogonal 
$\fg/\fri(\fg)^{\perp}$-module. Moreover, 
$$ 0\longrightarrow \fri(\fg)^{\perp}/\fri(\fg)\stackrel 
i\longrightarrow  \fg/\fri(\fg)\stackrel 
p\longrightarrow \fg/\fri(\fg)^{\perp}\longrightarrow 0$$ 
is an exact sequence of $(\fh,K)$-equivariant Lie algebras. 
\begin{co}\label{C1} 
   For any $(\fh,K)$-equivariant metric Lie algebra $(\fg,\Phi,\ip)$ 
   without simple ideals the 
   quadruple $(\fg,\fri(\fg),i,p)$ is a quadratic extension of 
   $\fg/\fri(\fg)^{\perp}$ by $\fri(\fg)^{\perp}/\fri(\fg)$. 
\end{co}     
This extension will be called the canonical quadratic 
extension\index{extension!canonical quadratic} 
associated with $(\fg,\Phi,\ip)$. 

\begin{ex} {\rm
     The following example shows that for a given metric Lie 
     algebra $(\fg,\ip)$ there may exist other quadratic extensions 
     $(\fg,\fri,i,p)$ than the canonical one. Let 
     $\fh(1)=\{[X_{1},X_{2}]=X_{3}\}$ be the 
     three-dimensional 
     Heisenberg algebra and let $\sigma^{1},\sigma^2,\sigma^{3}$ be 
     the 
     basis of $\fh(1)^{*}$ that is dual to $X_{1},X_{2},X_{3}$. Let us consider the metric Lie algebra 
     $\fg:=\fd_{0,0}(\fh(1),0)=\fh(1)^{*}\rtimes \fh(1)$. 
     Example~\ref{sm} says that $(\fg,\fh(1)^{*},i,p)$ is a quadratic 
     extension of $\fh(1)$ by $0$, where $(i,p)$ is defined by 
     $$0\stackrel {i=0}\longrightarrow  \fg/\fh(1)^{*}\stackrel   
     p\cong \fh(1) \longrightarrow 0.$$ 
     However, this quadratic extension is not the canonical one. 
     Indeed, we have 
     $R(\fg)=\fg'=\Span\{X_{3},\sigma^{1},\sigma^{2}\}$. In 
     particular, $R(\fg)\subset\fz(\fg)$, hence $R_{2}(\fg)=0$. 
     This implies 
     $$\fri(\fg)=R(\fg)^{\perp}\cap R(\fg)=R(\fg)= 
     \Span\{X_{3},\sigma^{1},\sigma^{2}\}.$$ 
     In particular, the canonical quadratic extension associated 
     with $\fg$ is a quadratic extension 
     of $\fg/\fri(\fg)^{\perp}\cong \RR^{3}\not\cong \fh(1)$ by 
     $\fri(\fg)^{\perp}/\fri(\fg)=0$. }
\end{ex}     
     
Hence, at first glance we have the same difficulty as for double 
extensions, namely, in general an $(\fh,K)$-equivariant metric Lie algebra 
can be obtained in different ways by quadratic extensions. However, 
now we can always distinguish one of these extensions, namely the canonical 
one. As a quadratic extension 
this extension is characterised by the property to be balanced in the following 
sense. 
\begin{de} 
A quadratic extension $(\fg,\fri,i,p)$ of an $(\fh,K)$-equivariant 
Lie algebra $(\fl,\Phi)$ by an 
orthogonal $(\fl,\Phi_{\fl})$-module $\fa$ is 
called balanced\index{extension!balanced quadratic} if $\fri=\fri(\fg)$. 
\end{de}   
Since our aim is to determine isomorphism classes of 
$(\fh,K)$-equivariant metric Lie algebras Corollary~\ref{C1} leads  us to 
the problem to decide for which balanced quadratic extensions 
$(\fg_{1},\fri_{1},i_{1},p_{1})$ and $(\fg_{2},\fri_{2},i_{2},p_{2})$ 
the $(\fh,K)$-equivariant metric Lie algebras $\fg_{1}$ and $\fg_{2}$ 
are isomorphic. We will divide this problem into two steps. First we 
will introduce an equivalence relation for quadratic extensions that 
is stronger than isomorphy of the underlying $(\fh,K)$-equivariant metric 
Lie algebras. We will describe the corresponding equivalence classes. 
In the second step we have to decide which equivalence classes of 
quadratic extensions have isomorphic underlying $(\fh,K)$-equivariant metric 
Lie algebras. 

\begin{de} Two quadratic extensions $(\fg_j,\fri_{j},i_{j},p_{j})$, 
$j=1,2$, of 
$(\fl,\Phi_{\fl})$ by $\fa$ are 
called equivalent if there exists an isomorphism $F: 
\fg_{1}\rightarrow \fg_{2}$ of $(\fh,K)$-equivariant metric Lie algebras 
that maps $\fri_1$ to $\fri_2$ and satisfies $\overline F\circ 
i_1=i_2$ and $p_2\circ\overline F=p_{1}$, where $\overline 
F:\fg_{1}/\fri_{1}\rightarrow \fg_{2}/\fri_{2}$ is the map induced by 
$F$. 
\end{de} 

Similar to the case of ordinary extensions of Lie algebras one can 
describe equivalence classes of quadratic extensions by 
cohomology\index{cohomology} classes. We will introduce a suitable cohomology theory in the next 
section. Actually, for a given $(\fh,K)$-equivariant Lie 
algebra $(\fl,\Phi_{\fl})$ and an orthogonal 
$(\fl,\Phi_{\fl})$-module $\fa$ we will define a cohomology set 
$\HPhi$ and a subset $\HPhi_{b}\subset\HPhi$ such that the following 
holds. 

\begin{theo}\label{T1}  There is a bijective map $\Psi$ from the set of 
equivalence classes of quadratic extensions of $(\fl,\Phi_{\fl})$ by 
$\fa$ to $\HPhi$. The image under $\Psi$ of the subset of all equivalence classes 
of balanced extensions  equals $\HPhi_{b}\subset\HPhi$. 
\end{theo} 

The set $\HPhi$ consists of equivalence classes 
$[\alpha,\gamma]$ of cocycles 
$(\alpha,\gamma)\in\cZ(\fl,\Phi_\fl,\fa)$ with respect to a certain 
equivalence relation. The inverse of $\Psi$ then maps $[\alpha,\gamma]$ 
to the equivalence class of the standard model $\dd$ (Example~\ref{sm}). For an explicit description of the map 
$\Psi$ see Section \ref{S24}. 
For a proof of this theorem in the non-equivariant case see 
\cite{KO1}, Thm.~2.7 and Thm.~3.12.

\subsection{Quadratic cohomology}\label{S23} 
The aim of this section is the exact definition 
of the cohomolgy sets that appear in Theorem \ref{T1}. Since quadratic 
extensions are not ordinary Lie algebra extensions we cannot 
expect to describe them by usual 
Lie algebra cohomology.   We need a kind of non-linear cohomology. 
Such cohomology sets were first introduced by Grishkov \cite{Gr} in a 
rather general setting. For the 
special case of cohomology needed for quadratic extensions 
we gave a self-contained presentation in \cite{KO1}. Neither 
\cite{Gr} nor \cite{KO1} deals with the equivariant situation.
As we will see here, the $(\fh,K)$-action can be easily incorporated. 

Let us first recall the construction of the usual Lie algebra 
cohomology\index{cohomology!Lie algebra}. Let $\rho:\ \fl\rightarrow \gl(\fa)$ be a 
representation of a Lie algebra $\fl$ on a vector space $\fa$. Then we have the 
standard Lie algebra 
cochain complex $(C^*(\fl,\fa),d)$, where 
$C^p(\fl,\fa)= \Hom(\textstyle{\bigwedge ^p}\fl,\fa)$ and   
$d:C^p(\fl,\fa)\rightarrow 
C^{p+1}(\fl,\fa)$ is defined 
by 
\begin{eqnarray*} 
    d\tau(L_1,\dots,L_{p+1})&=&\sum_{i=1}^{p+1}(-1)^{i-1} 
    \rho(L_{i})\tau(L_{1},\ldots, \hat L_i,\ldots, L_{p+1})\\   
    &&+\sum_{i<j}(-1)^{i+j}\tau([L_{i},L_{j}],L_{1},\ldots,\hat 
    L_i,\ldots,\hat L_{j},\ldots,L_{p+1}). 
\end{eqnarray*} 
The corresponding cohomology groups are 
denoted by $H^p(\fl,\fa)$. 
In the special case where $\fa$ is the one-dimensional trivial representation 
of $\fl$  we denote the standard cochain complex also by 
$(C^*(\fl),d)$ and the cohomology groups by $H^{p}(\fl)$. 

Suppose that we have two Lie algebras $\fl_{i}$, $i=1,2$ and orthogonal 
$\fl_{i}$-modules $\fa_{i}$ and that 
$(S,U):(\fl_{1},\fa_{1})\rightarrow (\fl_2,\fa_2)$ is a morphism of 
pairs. Then we have pull back maps 
\begin{eqnarray*} 
(S,U)^{*}:&&C^{p}(\fl_{2},\fa_{2})\longrightarrow 
C^{p}(\fl_{1},\fa_{1}) \\&& 
(S,U)^{*}\alpha(L_{1},\ldots,L_{p}):=U\circ\alpha(S(L_{1}),\ldots,S(L_{p}))\\[2ex] 
(S,U)^{*}: && C^{p}(\fl_{2})\longrightarrow C^{p}(\fl_{1}),  \\ 
&&(S,U)^{*}\gamma(L_{1},\ldots,L_{p}):=\gamma(S(L_{1}),\ldots,S(L_{p})). 
\end{eqnarray*} 

Now let $(\fl,\Phi_\fl)$ be an $(\fh,K)$-equivariant Lie algebra and let $\fa$ 
be an orthogonal $(\fh,K)$-module. Then we can consider 
\begin{equation}\label{EX} 
    (e^{\Phi_{\fl}(X)}, e^{-\Phi_{\fa}(X)}):\ (\fl,\fa)\longrightarrow   
    (\fl,\fa),\ X\in\fh, 
\end{equation} 
and 
\begin{equation}\label{Ek} 
    (\Phi_{\fl}(k),(\Phi_{\fa}(k))^{-1}):\ (\fl,\fa)\longrightarrow   
    (\fl,\fa),\ k\in K, 
\end{equation} 
as morphism of pairs (without $(\fh,K)$-structure). Let 
$C^p(\fl,\fa)^{(\fh,K)} \subset C^p(\fl,\fa)$ denote the 
subspace of cochains that are invariant under these morphisms of 
pairs for all $X\in\fh$ and $k\in K$. 

We define a product $C^p(\fl,\fa)^{(\fh,K)}\times 
C^q(\fl,\fa)^{(\fh,K)} \rightarrow C^{p+q}(\fl)^{(\fh,K)}$ by 
$$\langle\cdot\wedge\cdot\rangle: C^p(\fl,\fa)^{(\fh,K)}\times 
C^q(\fl,\fa)^{(\fh,K)}\stackrel{\wedge}{\longrightarrow} 
C^{p+q}(\fl,\fa\otimes \fa)^{(\fh,K)} \ 
\stackrel{\ipa}{\longrightarrow} \ C^{p+q}(\fl)^{(\fh,K)}.$$ 
Now we define the set of quadratic $1$-cochains to be 
$${\cal C} ^1_Q(\fl,\Phi_{\fl},\fa)= C^1(\fl,\fa)^{(\fh,K)}\oplus 
C^2(\fl)^{(\fh,K)}. 
$$ 
This set is a group with group operation defined by 
$$(\tau_1,\sigma_1)*(\tau_2,\sigma_2)=(\tau_1+\tau_2, \sigma_1 
+\sigma_2 +\textstyle{\frac12} \langle \tau_1\wedge 
\tau_2\rangle)\,.$$ 
Let us consider the set
$${\cal Z} ^{2}_Q(\fl,\Phi_{\fl},\fa)=\{(\alpha,\gamma) \in 
C^{2}(\fl,\fa)^{(\fh,K)}\oplus 
C^{3}(\fl)^{(\fh,K)} \mid d\alpha=0,\ 
d\gamma=\textstyle{\frac12}\langle\alpha \wedge\alpha\rangle\} $$ 
whose elements are called quadratic $2$-cocycles. Then the group
${\cal C} ^{1}_Q(\fl,\Phi_{\fl},\fa)$ acts from the right on ${\cal Z} 
^{2}_Q(\fl,\Phi_{\fl},\fa)$ by 
$$(\alpha,\gamma)(\tau,\sigma)=\left(\,\alpha +d\tau,\gamma +d\sigma 
+\langle(\alpha +\textstyle{\frac12} d\tau)\wedge\tau\rangle\,\right).$$ 
We define the quadratic cohomology\index{cohomology!quadratic} set as the orbit space 
$${\cal H}^{2}_Q(\fl,\Phi_{\fl},\fa):={\cal Z}^{2}_Q(\fl,\Phi_{\fl},\fa)/ {\cal 
C}^{1}_Q(\fl,\Phi_{\fl},\fa).$$ 
The equivalence class of 
$(\alpha,\gamma)\in {\cal Z} ^{2}_Q(\fl,\Phi_{\fl},\fa)$ in ${\cal 
H}^{2}_Q(\fl,\Phi_{\fl},\fa)$ is denoted by $[\alpha,\gamma]$. 
As usual, if $(\fh,K)$ is trivial, then we omit $\Phi_{\fl}$ in all the 
notation above. 

Now let us consider a morphism of pairs 
$(S,U):((\fl_1,\Phi_1),\fa_1)\rightarrow ((\fl_2,\Phi_2),\fa_2)$. 
As discussed above $(S,U)$ acts on $C^{p}(\fl,\fa)$. By 
$(\fh,K)$-equivariance $(S,U)$ maps the subspace 
$C^{p}(\fl_{2},\fa_{2})^{(\fh,K)}\subset C^{p}(\fl_{2},\fa_{2})$ to 
$C^{p}(\fl_{1},\fa_{1})^{(\fh,K)}\subset C^{p}(\fl_{1},\fa_{1})$. 
It is not hard to prove that this map induces a map 
$$(S,U)^{*}:\ \cH(\fl_{2},\Phi_{\fl_{2}},\fa_{2})\longrightarrow 
\cH(\fl_{1},\Phi_{\fl_{1}},\fa_{1})$$ 
(cf.~\cite{KO1} for a proof in the non-equivariant case). 

In particular, for a given $(\fh,K)$-equivariant Lie algebra 
$(\fl,\Phi_{\fl})$ and an orthogonal $(\fl,\Phi_{\fl})$-module 
$\fa$ the morphisms of pairs in (\ref{EX}) and (\ref{Ek}) 
induce maps on $\cH(\fl,\fa)$. Let $\cH(\fl,\fa)^{(\fh,K)}\subset\cH(\fl,\fa)$ 
denote 
the subset of all cohomology classes that are invariant under these 
maps for all $X\in\fh$ and $k\in K$. 

\begin{pr}\label{inva} 
    Under our assumptions on $(\fh,K)$ the inclusion 
    $$\cZ(\fl,\Phi_{\fl},\fa)\hookrightarrow \cZ(\fl,\fa)$$ induces a 
    bijection 
    $$\HPhi     \longrightarrow \cH(\fl,\fa)^{(\fh,K)}.$$
\end{pr} 
You can find a proof of this proposition in the appendix. 

Next we will define the subset $\HPhi_{b}\subset\HPhi$, which plays 
an important role in Theorem~\ref{T1}. It was first introduced in 
\cite{KO1}, where you can also find a proof of the fact that its 
elements correspond exactly to those extensions that are balanced. 

As usual, an $(\fl,\Phi_{\fl})$-module 
$(\rho,\fa)$ will be called semisimple if every   
$(\fl,\Phi_{\fl})$-submodule has an  $(\fl,\Phi_{\fl})$-invariant 
complement. This is the case if and only if $\rho$ is semisimple. 
In the following definition we need the notion of the socle $S(\fl)$ of a 
Lie algebra $\fl$, which is the maximal ideal of $\fl$ on which $\fl$ acts 
semisimply. 

\begin{de} Let $(\fl,\Phi_{\fl})$ be an $(\fh,K)$-equivariant Lie 
algebra, let $(\rho,\fa)$ be a semisimple orthogonal 
$(\fl,\Phi_{\fl})$-module and take $(\alpha,\gamma)\in \ZPhi$. 
Since $\rho$ is semisimple we have a decomposition $\fa=\fa^\fl\oplus 
\rho(\fl)\fa$ and a corresponding decomposition 
$\alpha=\alpha_0+\alpha_1$.  Let $m$ be such that $R_{m+1}(\fl)=0$. 
Then $(\alpha,\gamma)\in \ZPhi$ is called balanced if it satisfies the 
following conditions $(A_{k})$ and $(B_{k})$ for $0\le k\le m$. 
\begin{enumerate} 
\item[$(A_0)$] 
Let $L_0\in \fz(\fl)\cap \ker \rho$ be such that there exist 
elements 
$A_0\in \fa$ and $Z_0\in \fl^*$ satisfying 
for all $L\in\fl$ 
\begin{enumerate} 
\item[(i)] $\alpha(L,L_0)=\rho(L) A_0 $, 
\item[(ii)] $\gamma(L,L_0,\cdot)=-\langle A_0,\alpha(L,\cdot)\rangle_\fa +\langle 
Z_0, [L,\cdot]_\fl\rangle$ as an element of $\fl^*$, 
\end{enumerate} 
then $L_0=0$. 
\item[$(B_0)$] The subspace $\alpha_0(\ker \lb_\fl)\subset 
\fa^\fl$ is non-degenerate. 
\item[$(A_k)$] $(k\ge 1)$\\ 
Let $\fk\subset S(\fl)\cap R_k(\fl)$ be an $\fl$-ideal such that there exist 
elements 
$\Phi_1\in \Hom(\fk,\fa)$ and $\Phi_2\in \Hom(\fk,R_k(\fl)^*)$ satisfying 
for all $L\in\fl$ and $K\in\fk$ 
\begin{enumerate} 
\item[(i)] $\alpha(L,K)=\rho(L)\Phi_1(K)-\Phi_1([L,K]_\fl)$, 
\item[(ii)] $\gamma(L,K,\cdot)=-\langle \Phi_1(K),\alpha(L,\cdot)\rangle_\fa 
+\langle 
\Phi_2(K), [L,\cdot]_\fl\rangle +\langle \Phi_2([L,K]_\fl), \cdot \rangle$ as an 
element of $R_k(\fl)^*$, 
\end{enumerate} 
then $\fk=0$. 
\item[$(B_k)$] $(k\ge 1)$\\ 
Let $\fb_k\subset\fa$ be the maximal submodule such that the system of equations 
$$ \langle\alpha(L,K), 
B\rangle_\fa=\langle\rho(L)\Phi(K)-\Phi([L,K]_\fl),B\rangle_\fa\ , 
\quad L\in\fl, K\in R_k(\fl), B\in\fb_k, $$ 
has a solution $\Phi \in \Hom(R_k(\fl),\fa)$. Then $\fb_k$ is non-degenerate. 
\end{enumerate} 
\end{de} 
One can prove that for a cocycle the property to be balanced depends only 
on its cohomology class. Hence we may call a cohomology class 
$[\alpha,\gamma]\in\HPhi$ balanced\index{cohomology class!balanced} if $(\alpha,\gamma)\in \ZPhi$ is 
balanced. 

For an $(\fh,K)$-equivariant Lie algebra $(\fl,\Phi_{\fl})$ and 
an orthogonal $(\fl,\Phi_{\fl})$-module $(\rho,\fa)$ let 
$\HPhi_{b}\subset\HPhi$ 
be the set of all balanced cohomology classes if $\rho$ is semisimple 
and put $\HPhi_{b}:=\emptyset$ if $\rho$ ist not semisimple. 
This finishes the definition of the cohomology sets used in Theorem 
\ref{T1}. 

\begin{ex}\label{Eh1} 
    Take $\fl=\fh(1)=\{[X,Y]=Z\}$ and let $(\rho,\fa)$ be a semisimple orthogonal 
    $\fl$-module. Then the following two maps are bijective: 
    $$Z_\fl:=\{\alpha\in C^2(\fl,\fa)\mid \alpha(X,Y)=0, 
    \alpha(Z,\fl)\subset \fa^\fl\} \longrightarrow H^{2}(\fl,\fa),\   
    \alpha \mapsto [\alpha],$$ 
    $$ Z_{\fl,b}:=\{\alpha\in Z_{\fl}\mid \alpha\not=0, \alpha(Z,\fl) 
    \subset \fa^{\fl} 
    \mbox{\rm \,  
    is non-degenerate}\}\longrightarrow \cH(\fl,\fa)_{b},\   
    \alpha \mapsto [\alpha,0].$$ 
\end{ex} 
\proof Let us consider the first map. It is well-defined and 
injective. We will prove that it is surjective. Take $a\in 
H^{2}(\fl,\fa)$. Since $\fl$ is nilpotent and $\fa$ is semisimple we have 
$H^{*}(\fl,\fa)=H^{*}(\fl,\fa^{\fl})$ \cite{D}. Hence we can represent $a$ by a 
cocycle 
$\alpha$ that satisfies $\alpha(\fl,\fl)\subset \fa^{\fl}$. We 
define $\tau\in C^{1}(\fl,\fa)$ by $\tau(X)=\tau(Y)=0$, 
$\tau(Z)=\alpha(X,Y)$. Then $\tilde\alpha:=\alpha+d\tau$ 
satisfies $\tilde \alpha(X,Y)=\alpha(X,Y) - \tau([X,Y]) 
=0$ and $\tilde \alpha(Z,\fl)\subset \fa^\fl$. Since   
$[\alpha]=[\tilde \alpha]$ the assertion follows. 

Now let us turn to the second map. First we have to check that it is also 
well-defined. If $\alpha\in Z_\fl$, then $\langle \alpha\wedge \alpha 
\rangle =0$. Hence $(\alpha,0)\in\cZ(\fl,\fa)$. We have to show 
that the cocycle $(\alpha,0)$ is balanced if 
$\alpha(\fl,\fl)$ is non-degenerate and $\alpha\not=0$. We have to check the 
admissibility conditions 
$(A_k)$, $(B_k)$ for $k=0,1$ (note that $R_{2}(\fl)=0$). Since $\rho$ is semisimple 
and $Z\in R(\fl)$ we 
have 
$\rho(Z)=0$, hence $\fz(\fl)\cap\Ker \rho=\RR\cdot Z$. 
Because of $\fa^{\fl}\supset\alpha(Z,\fl)\not=0$ Conditions  $(A_0)$ 
and $(A_{1})$ are satisfied. Conditions $(B_{0})$ and 
$(B_{1})$ hold since $\alpha(Z,\fl)\subset\fa^{\fl}$ is non-degenerate by 
assumption. Hence the second map is well-defined. Obviously it is 
injective. Let us prove surjectivity. Suppose $a\in \cH(\fl,\fa)_{b}$. 
By surjectivity of our first map we can represent $a$ by a balaced 
cocycle $(\alpha,\gamma)$ with $\alpha\in Z_{\fl}$. Clearly, 
$\alpha=0$ would contradict Condition $(A_{0})$ (choose 
$L_{0}=Z$, $A_{0}=0$, $Z_{0}=0$). Hence $0\not=\alpha\in Z_{\fl}$. Now it is 
easy to see that there exists a cochain $\tau\in C^{1}(\fl,\fa)$ such 
that $d\tau=0$ and 
$\langle\alpha\wedge\tau\rangle=\gamma$. Consequently, 
$a=[\alpha,\gamma]=[\alpha,0]$, which proves the assertion. 
\qed 

Since we want to construct indecomposable metric Lie algebras by 
quadratic extensions we also need a notion of indecomposability for 
quadratic cohomology classes.

\begin{de}\label{kosel} 
A non-trivial decomposition of a pair $((\fl,\Phi_{\fl}),\fa)$ 
consists of two non-zero morphisms of pairs 
$ (q_i,j_i):((\fl,\Phi_{\fl}),\fa)\longrightarrow 
((\fl_i,\Phi_{\fl_i}),\fa_i)$, $i=1,2,$ 
such that $(q_1,j_1)\oplus (q_2,j_2):((\fl,\Phi_{\fl}),\fa)\rightarrow 
((\fl_1,\Phi_{\fl_1}),\fa_1)\oplus((\fl_2,\Phi_{\fl_2}),\fa_2)$ is 
an isomorphism. 

A cohomology class $\varphi\in 
\cH(\fl,\Phi_{\fl},\fa)$ is called decomposable if it can 
be written as a sum 
$$ \varphi=(q_1,j_1)^*\varphi_1 +(q_2,j_2)^*\varphi_2 $$ 
for a non-trivial decomposition $(q_i,j_i)$ of $((\fl,\Phi_\fl),\fa)$ and 
certain $\varphi_i\in \cH(\fl_i,\Phi_{\fl_i},\fa_i)$, $i=1,2$. 
Here addition is induced by addition in the vector space $C^{2}(\fl,\fa)\oplus 
C^{3}(\fl)$. A cohomology class which is not decomposable is called 
indecomposable. 
\end{de} 
Then we have the following relation between indecomposbility of 
cohomology classes and indecomposability of metric Lie algebras.
\begin{pr}[\cite{KO1}, Prop.~4.5]\label{hermlin} 
An $(\fh,K)$-equivariant metric Lie algebra $(\fg,\Phi,\ip)$ is 
indecomposable if and only if the image under $\Psi$ of the canonial 
quadratic extension associated with $(\fg,\Phi,\ip)$ is an 
indecomposable cohomology class. 
\end{pr} 

\subsection{A classification scheme}\label{S24} 

According to Corollary~\ref{C1}, each $(\fh,K)$-equivariant metric Lie algebra 
without simple ideals 
comes with a distinguished structure of a quadratic extension which is balanced. By 
Thm.~\ref{T1} balanced quadratic extensions are characterised by balanced quadratic 
cohomology classes. 
We obtain a (functorial) assignment 
\begin{eqnarray} 
\lefteqn{\hspace{-1.5cm}\{ (\fh,K)\mbox{-equivariant metric Lie algebras } 
(\fg,\Phi,\ip) \mbox{ without simple ideals}\}}\nonumber\\ 
&\Longrightarrow& \left\{\mbox{quadruples } 
\left(\fl,\Phi_\fl,\fa,[\alpha,\gamma]\in 
\HPhi_b\right)\right\}\ \label{rueffel} 
\end{eqnarray} 
In order to make (\ref{rueffel}) concrete let us write down the map $\Psi$ 
appearing 
in Thm.~\ref{T1} explicitly. 
Let $(\fg,\fri,i,p)$ be a quadratic extension of $(\fl,\Phi_{\fl})$ by 
$\fa$. Let $\tilde p:\fg\rightarrow \fl$ be the map induced by $p$. We 
choose an $(\fh,K)$-equivariant section $s: \fl\rightarrow \fg$ of 
$\tilde p$ with isotropic 
image (which exists by semisimplicity 
of $\Phi$) and define $\alpha\in C^{2}(\fl,\fa)^{(\fh,K)}$ 
and $\gamma\in C^{3}(\fl)^{(\fh,K)}$ by 
\begin{eqnarray*} 
i(\alpha(L_1,L_2))&:=&[s(L_1),s(L_2)]-s([L_1,L_2]_\fl)\, +\, \fri\quad 
\in\, 
\fg/\fri 
\\ 
\gamma(L_1,L_2,L_3)&:=& \langle\, [s(L_1),s(L_2)], s(L_3)\rangle \ . 
\end{eqnarray*} 
Then $(\alpha,\gamma)\in{\cal Z} ^{2}_Q(\fl,\Phi_{\fl},\fa)$, 
and $[\alpha,\gamma]\in\HPhi$ is the desired cohomology class. 
Note that $[\alpha,\gamma]$ does not depend on the choice of $s$ while 
$(\alpha,\gamma)$ does. 

It will turn out that the data on the right hand side 
of (\ref{rueffel}) give a very useful description of the set of isomorphism classes 
of 
$(\fh,K)$-equivariant metric Lie algebras. In fact, the resulting classification 
scheme, 
Theorem \ref{class} below, is the basis of most of the classification results for 
metric Lie algebras 
and symmetric spaces that will be presented in this article. 

Let $(\fl,\Phi_{\fl})$ be an $(\fh,K)$-equivariant Lie algebra. We consider the 
category $M^{ss}_{\fl,\Phi_{\fl}}$ 
of semisimple orthogonal $(\fl,\Phi_{\fl})$-modules, where the morphisms between 
two modules 
$\fa_1,\fa_2$ are given 
by morphism of pairs $(S,U): ((\fl,\Phi_{\fl}),\fa_1)\rightarrow 
((\fl,\Phi_{\fl}),\fa_2)$. We denote the automorphism group of an object $\fa$ of 
$M^{ss}_{\fl,\Phi_{\fl}}$ 
by $G_{\fl,\Phi_{\fl},\fa}$. 
The natural right action of $G_{\fl,\Phi_{\fl},\fa}$ on $\HPhi$ leaves $\HPhi_{b}$ 
invariant. 

\begin{theo}[compare \cite{KO1}, Thm.~4.6]\label{class} 
Let $\cL$ be a complete set of representatives of iso\-morphism classes of 
$(\fh,K)$-equivariant Lie algebras. 
For each $(\fl,\Phi_{\fl})\in\cL$ we choose a complete set of representatives 
$\cA_{\fl,\Phi_{\fl}}$ of isomorphism classes of objects in~$M^{ss}_{\fl,\Phi_{\fl}}$. 

Then (\ref{rueffel}) descends to a bijective map from the set of isomorphism 
classes 
of 
 $(\fh,K)$-equivariant metric Lie algebras without simple ideals 
 to the union of orbit spaces 
 $$\coprod_{(\fl,\Phi_{\fl})\in\cL}\  \coprod_{\fa\in\cA_{\fl,\Phi_{\fl}}}\   
 \HPhi_{b}/ G_{\fl,\Phi_{\fl},\fa}\ .$$ 
The inverse of this map sends the orbit of $[\alpha,\gamma]\in\HPhi_{b}$ 
to the isomorphism class of $\dd$. 

An $(\fh,K)$-equivariant metric Lie algebra is indecomposable if and only if 
the corresponding $G_{\fl,\Phi_{\fl},\fa}$-orbit in $\HPhi_{b}$ consists of 
indecomposable 
cohomology classes. 
\end{theo} 
In particular, the theorem says that the $(\fh,K)$-equivariant metric Lie algebras 
$\dd$ 
exhaust all isomorphism classes of $(\fh,K)$-equivariant metric Lie algebras 
without 
simple ideals. 

Its proof is word by word the same as the one for the non-equivariant case that 
is given in complete detail in \cite{KO1}. The main ingredient is Theorem~\ref{T1}, 
which actually involves the 
construction of the standard model $\dd$. In addition, one 
has to show that two metric Lie algebras are isomorphic if and only if their 
associated balanced 
cohomology classes are mapped to each other by an isomorphism of pairs. 
This is a rather straightforward 
consequence of the functoriality of the canonical isotropic ideal 
$\fri(\fg)$ (\cite{KO1}, 
Prop.~4.2). The statement about indecomposability follows from 
Prop.~\ref{hermlin}.

Theorem \ref{class} provides a complete set of invariants for equivariant metric 
Lie 
algebras 
and therefore structures the set of all isomorphism classes of them in a certain 
way. 
At first glance, however, Theorem \ref{class} might look rather useless for 
classification problems. Indeed, the classification 
scheme seems to involve the classification of all (equivariant) Lie algebras. 
While this reflects the real difficulty of the problem there are at least two 
circumstances that, 
nevertheless, 
allow for interesting applications of the theorem. First of all, the index of $\dd$ is at least $\dim\fl$. Thus for classification of equivariant metric Lie 
algebras of small 
index only Lie algebras $\fl$ of small dimensions are needed. In general, it is a 
good strategy 
to look for certain subclasses of all equivariant metric Lie algebras such that the 
corresponding ingredients 
in Theorem \ref{class} become manageable. 

Secondly, a great many 
Lie algebras $\fl$ 
satisfy $\cH(\fl,\fa)_{b}=\emptyset$ for all semi-\linebreak simple orthogonal 
$\fl$-modules $\fa$ 
(in view of Proposition \ref{inva} this also implies \linebreak 
$\HPhi_{b}=\emptyset$ 
for 
any equivariant structure $\Phi_\fl$ on $\fl$). We call 
these Lie algebras non-admissible and the remaining ones, which really occur in 
Theorem \ref{class}, admissible\index{Lie algebra!admissible}. E.g., the two-dimensional non-abelian Lie algebra 
and Heisenberg Lie algebras of dimension $\ge 5$ are non-admissible (\cite{KO1}, 
Prop.~5.2). It is easy to see that all reductive Lie algebras 
are admissible and that the class of admissible Lie algebras is closed under 
forming 
direct sums. 
Solvable admissible Lie algebras $\fl$ with 
$\dim\fl'\le 2$ are classified in \cite{Knil} and \cite{KO1}, Section 5. Up to 
abelian summands, 
there are only finitely many of them. This result (together with its proof) shows 
that an a priori 
classification of Lie algebras (within a certain class) is not really needed for 
applications of 
Theorem~\ref{class}. Up to now it is unknown how large the class of admissible Lie 
algebras really 
is. Any new structure result for admissible Lie algebras would give highly 
desirable 
new insight 
into the world of (equivariant) metric Lie algebras. 

Contrary to $\cL$, the determination of $\cA_{\fl,\Phi_{\fl}}$ for given 
$(\fl,\Phi_{\fl})$ is usually 
no problem (since only semisimple modules are needed). Also the computation of the 
cohomology 
sets $\HPhi_{b}$ succeeds in many interesting cases (techniques of homological 
algebra are helpful here). 
To achieve a true classification, 
one would have to determine, as a last step, the orbit space of the action of the 
group  $G_{\fl,\Phi_{\fl},\fa}$ on the set $\HPhi_{b}$. This may lead to unsolved 
classification problems 
again, e.g., to the description of the orbit space 
$\bigwedge^{3}\RR^{n}/\GL(n,\RR)$ for large $n$ (compare \cite{KO0}, Sect.~5). But often $\HPhi_{b}$ is so small such that things can be carried out 
explicitly.

\subsection{Classification results for metric Lie algebras}\label{S25} 
For trivial $(\fh,K)$, i.e.~for $(\fh,K)=(0,\{e\})$ Theorem~\ref{class} 
specialises to a classification scheme for metric Lie algebras. 
As mentioned above for certain classes of metric Lie algebras this can be used to 
obtain a 
full classification 
(in the sense of a list).  In particular we can use it to classify 
metric Lie algebras of small index.

Let us consider the Lorentzian case, i.e. the case where the 
metric Lie algebra has index~1. We already know examples of 
Lorentzian metric Lie algebras\index{metric Lie algebra!Lorentzian}, namely 
the oscillator algebras\index{oscillator algebra} $\osc(\lambda)$ 
defined in 
Section~\ref{S21}. For a classification of indecomposable Lorentzian 
metric Lie algebras we have to restrict to $\fl$ 
with $\dim \fl \le 1$ and to Euclidean $\fl$-modules $\fa$  in our 
classification scheme. 
In this way we reproduce the following well-known result, which was 
originally proved using double extensions. 

\begin{theo}[Medina \cite{M}] \label{lorentz} Each indecomposable non-simple metric 
Lie algebra of signature $(1,q)$, $q>0$, is isomorphic to an oscillator algebra 
$\osc(\lambda)$ for exactly one 
$\lambda=(\lambda^{1},\ldots,\lambda^{m})\in (\RR^{*})^{m}\cong 
\RR^{m}$, $ q=2m+1$, with  $\lambda^{1}=1\le \lambda^{2}\le \ldots \le 
\lambda^{m}$. 
\end{theo}     
For the study of Lorentzian metric Lie algebras the method of double 
extensions and the method of quadratic extensions are in some sense 
equivalent, since every indecomposable non-simple Lorentzian metric 
$\fg$ Lie algebra has exactly one isotropic ideal, namely its centre 
$\fz(\fg)$. However, as already mentioned, for classification of metric Lie 
algebras of higher index quadratic extensions are more useful than 
double extensions, see for example the classification of 
metric Lie algebras of index 2 in \cite{KO0} and the classification of metric 
Lie algebras of index 3 in \cite{KO1}. Both of them are based 
on Theorem~\ref{class}. 

As a further example that shows that Theorem~\ref{class} is a useful mean for 
concrete classification 
problems let us consider nilpotent metric Lie algebras\index{metric 
Lie algebra!nilpotent} of small 
dimension. Favre and Santharoubane classified such Lie algebras up to 
dimension~7 in \cite{FS}. Their proof is based on the double extension 
method. In \cite{Knil}, Theorem~\ref{class} has been used to give a 
classification of nilpotent metric Lie algebras of dimension at most~10. 
Most of the isomorphism classes are isolated ones, however, in dimension 10 
also 1-parameter families occur. 
Here we will restrict ourselves to nilpotent metric Lie algebras of dimension at 
most 9.

In the following theorem $\fa$ stands for a pseudo-Euclidean vector 
space that we consider always as a trivial $\fl$-module. As usual, let $e_1,\ldots, 
e_{p+q}$ be the standard basis of $\RR^{p,q}$. Furthermore, 
$\sigma^1,\ldots,\sigma^{l}\in\fl^{*}$ will denote the dual basis of a given 
basis $X_{1},\ldots,X_{l}$ of $\fl$ and $\sigma^{{i_1}\ldots 
i_{j}}:=\sigma^{i_1}\wedge \ldots \wedge \sigma^{i_{j}}$. 
\begin{theo} [\cite{Knil}]
If $(\fg,\ip)$ is an indecomposable non-abelian nilpotent metric Lie 
algebra of dimension at most 9, then it is isomorphic to 
$\fd_{\alpha,\gamma}(\fl,\fa)$ for exactly one of the data in the 
following list 
\begin{enumerate} 
    \item 
    $\fl=\fg_{4,1}=\{[X_{1},X_{2}]=X_{3},[X_{1},X_{3}]=X_{4}\}$, 
    $\ \fa\in\{\RR^{0,1},\RR^{1,0}\}$,\\[0.7ex] 
    $\alpha=\sigma^{14}\otimes e_{1}$, 
    $\gamma\in\{0,\sigma^{234},\sigma^{134},-\sigma^{134}\}$; 
    \item 
    $\fl=\fh(1)\oplus\RR=\{[X_{1},X_{2}]=X_{3}\}\oplus\RR\cdot X_{4}$, 
    $\ \fa\in\{\RR^{0,1},\RR^{1,0}\}$,\\[0.7ex] 
    $(\alpha,\gamma)=(\sigma^{13}\otimes  e_{1},\sigma^{234})$; 
    \item 
    $\fl=\RR^{4}$, $\ \fa\in\{\RR^{0,1},\RR^{1,0}\}$, 
    $\ (\alpha,\gamma)=(\sigma^{13}\otimes  e_{1},\sigma^{234})$; 
    \item 
    $\fl=\fh(1)=\{[X_{1},X_{2}]=X_{3}\}$, 
    \begin{enumerate} 
        \item $ \fa\in\{\RR^{0,1},\RR^{1,0}\}$, $\ 
(\alpha,\gamma)=(\sigma^{13}\otimes   
        e_{1},0)$, 
        \item $ \fa\in\{\RR^{0,2},\RR^{2,0},\RR^{1,1}\}$, 
        $\ (\alpha,\gamma)=(\sigma^{13}\otimes  e_{1}+\sigma^{23}\otimes   
e_{2},0)$; 
    \end{enumerate}         
    \item 
    $\fl=\RR^{3}$, 
    \begin{enumerate} 
        \item $\fa=0$, $\ (\alpha,\gamma)=(0,\sigma^{123})$, 
        \item $ \fa\in\{\RR^{0,2},\RR^{2,0},\RR^{1,1}\}$, $\ 
        (\alpha,\gamma)=(\sigma^{12}\otimes e_{1}+\sigma^{13}\otimes 
        e_{2},0)$, 
        \item $ \fa\in\{\RR^{0,3},\RR^{2,1},\RR^{1,2},\RR^{3,0}\}$, 
        $\ (\alpha,\gamma)=(\sigma^{12}\otimes e_{1}+\sigma^{13}\otimes 
        e_{2}+\sigma^{23}\otimes e_{3},0)$; 
    \end{enumerate}         
    \item 
    $\fl=\RR^{2}$, $\ \fa\in\{\RR^{0,1},\RR^{1,0}\}$, 
    $\ (\alpha,\gamma)=(\sigma^{12}\otimes  e_{1},0)$. 
\end{enumerate}     
\end{theo}

\section{Symmetric spaces}\label{S3} \index{symmetric space!pseudo-Riemannian|(}
\subsection{Symmetric triples and quadratic extensions} 
\label{S31}

In this section we are concerned with $\ZZ_2$-equi\-variant objects, i.e., 
$(\fh,K)$-equivariant modules, Lie algebras, metric Lie algebras etc., where 
$K=\ZZ_2$ is the group consisting of two elements and $\fh=0$. 
Any equivariant structure $\Phi$, $\Phi_\fl\dots$ is determined by its 
value on the nontrivial 
element of $\ZZ_2$, which is an involution. We will denote this involution by 
$\theta$, $\theta_\fl\dots$ . 
We will keep the notation of Section \ref{S2} but with all $\Phi$'s replaced 
by the corresponding $\theta$'s. Any $\ZZ_2$-module $V$ has a decomposition 
$V=V_+\oplus V_-$ into the 
$(\pm 1)$-eigenspaces of $\theta_V$. 

As explained in the introduction, one associates with a pseudo-Riemannian symmetric space 
the Lie algebra of its transvection group\index{transvection group} together with a natural 
non-degenerate symmetric bilinear form and an isometric involution on it. This 
leads to the 
notion of a symmetric triple\index{symmetric triple}. 
\begin{de}\label{sytri} 
    \mbox{} 
\begin{enumerate} 
\item[(a)] A $\ZZ_2$-equivariant Lie algebra $(\fg,\theta)$ is called proper if 
$[\fg_-,\fg_-]=\fg_+$. 
\item[(b)] A symmetric triple is a proper $\ZZ_2$-equivariant metric Lie algebra 
$(\fg,\theta,\ip)$. 
\item[(c)] The index (signature) of a symmetric triple $(\fg,\theta,\ip)$ 
is the index (signature) of the 
symmetric bilinear form 
$\ip|_{\fg_-}$. 
\end{enumerate} 
\end{de} 
The Lie algebra of the transvection group of a pseudo-Riemannian symmetric space 
$M$ of 
index $p$ carries the structure of a 
symmetric triple of the same index in a natural way. We call it the symmetric 
triple of $M$. 
The notions of isomorphy and decomposability carry over from general 
$\ZZ_2$-equivariant metric Lie algebras 
to symmetric triples. Then we have 
\begin{pr}[see e.g.~\cite{CP}]\label{otto} 
The assignment which sends each pseudo-Rie\-mannian symmetric space to its 
symmetric 
triple induces a bijective map 
between isometry classes 
of simply connected symmetric spaces and isomorphism classes of symmetric triples. 
A 
symmetric space is indecomposable if and only if its symmetric triple is so. 
\end{pr} 
All (not necessarily simply connected) symmetric spaces corresponding to a given 
symmetric triple $(\fg,\theta,\ip)$ 
can be easily classified. Let us shortly discuss this classification. For the facts 
used in this discussion we refer to 
\cite{koh} and \cite{H}, Ch.VII, \S\S 8,9.   

Let $\tilde G$ be the connected and simply connected 
Lie group with Lie algebra $\fg$. Then $\theta$ integrates 
to an involutive automorphism $\tilde \theta:\tilde G\rightarrow \tilde G$. 
The group  $\tilde G^{\tilde \theta}$ is 
connected. 
Let $Z(\tilde G)$ be the center of $\tilde G$, and set $Z_0:=Z(\tilde G)\cap\tilde 
G^{\tilde \theta}$. Then $Z_0$ is discrete. Choose a discrete $\tilde\theta$-stable 
subgroup $Z\subset Z(\tilde G)$ containing $Z_0$ and set $G:=\tilde G/Z$. Then 
$\tilde\theta$ induces an automorphism $\theta: G\rightarrow G$. Its fixed point 
group 
$G^\theta$ has at most finitely connected components. The connected component 
$G^\theta_0$ satisfies 
$ 
G^\theta_0\cap Z(G)=\{e\}$. 
We choose a group $G_+$ such that $G^\theta_0\subset G_+\subset G^\theta$ and 
$G_+\cap Z(G)=\{e\}$ and define 
$M:=G/G_+$. 
Let $x_0=eG_+$ be the base point of the homogeneous space $M$. Then $\fg_-\cong 
T_{x_0}M$, 
and $\ip|_{\fg_-}$ defines a $G_+$-invariant scalar product on $T_{x_0}M$, which 
extends uniquely to a $G$-invariant 
pseudo-Riemannian metric on $M$. Moreover, $\tilde \theta$ induces an involutive 
isometry $\theta_{x_0}$ of $M$ 
having $x_0$ as an isolated fixed point. Conjugating $\theta_{x_0}$ with elements 
of 
$G$ we get involutive 
isometries $\theta_x$ for all $x\in M$. Thus $M$ has the structure of a 
pseudo-Riemannian symmetric space. Moreover, 
$G$ is the transvection group of~$M$. 

If $Z$ and $G_+$ run through all possible choices as above, then the resulting 
spaces $M$ exhaust the isometry classes of pseudo-Riemannian spaces having a 
symmetric triple 
isomorphic to $(\fg,\theta,\ip)$. Two such spaces are 
isometric if and only if the defining data $(Z,G_+)$ are conjugated 
by the automorphism group $\Aut(\tilde G,\tilde\theta,\ip)$ consisting of all 
automorphisms of $\tilde G$ that respect the involution as well as the 
pseudo-Riemannian metric on $\tilde G$ induced by $\ip$. 
The simply connected symmetric space associated with $(\fg,\theta,\ip)$ arises if 
we 
choose 
$Z=Z_0$, $G_+=G^\theta_0$. 

Thus the classification of symmetric spaces is reduced to the classification of 
symmetric triples, i.e., 
proper $\ZZ_2$-equivariant metric Lie algebras. The theory developed in Section 
\ref{S2} associates with every $\ZZ_2$-equivariant metric Lie algebra 
$(\fg,\theta,\ip)$ (without simple ideals) via its canonical quadratic 
extension a quadruple $(\fl,\theta_\fl,\fa,[\alpha,\gamma])$, where 
\begin{itemize} 
\item $(\fl,\theta_\fl)$ is a $\ZZ_2$-equivariant Lie algebra, 
\item $\fa$ is a semisimple orthogonal $(\fl,\theta_\fl)$-module, and 
\item $[\alpha,\gamma]\in \cH(\fl,\theta_\fl,\fa)_b$ is a balanced quadratic 
cohomology class. 
\end{itemize} 
Then $(\fg,\theta,\ip)$ is isomorphic to the $\ZZ_2$-equivariant metric Lie algebra 
$\fd_{\alpha,\gamma}(\fl,\theta_\fl,\fa)$. 

Because we want to apply the classification scheme Theorem \ref{class} to symmetric 
triples we have to express the properness condition for $(\fg,\theta,\ip)\cong 
\fd_{\alpha,\gamma}(\fl,\theta_\fl,\fa)$ 
in Definition \ref{sytri} in terms of $(\fl,\theta_\fl,\fa,[\alpha,\gamma])$. 

\begin{pr}[\cite{KO2}, Sect.~5]\label{miau} 
Let 
$(\fl,\theta_\fl,\fa,[\alpha,\gamma])$ be a quadruple as above. 
By semisimplicity $\fa=\fa^\fl\oplus \rho(\fl)\fa$, and we have a 
corresponding decomposition $\alpha=\alpha_0+\alpha_1$. 
Then $\fd_{\alpha,\gamma}(\fl,\theta_\fl,\fa)$ is a symmetric triple if and only if 
\begin{itemize} 
\item[$(T_1)$] the $\ZZ_2$-equivariant Lie algebra $(\fl,\theta_\fl)$ is proper, 
and 
\item[$(T_{2})$] $\quad \fa^\fl_+=\alpha_0(\Ker \lb_{\fl_-})\ .$ 
\end{itemize} 
\end{pr} 

\begin{de} 
Let $(\fl,\theta_\fl)$ be a proper $\ZZ_2$-equivariant Lie algebra, and let 
$\fa$ be a semisimple orthogonal $(\fl,\theta_\fl)$-module. A quadratic extension 
$(\fg,\fri,i,p)$ 
of $(\fl,\theta_\fl)$ by $\fa$ is called 
admissible\index{extension!admissible quadratic} if it is balanced and $\fg$ is 
proper, i.e., a symmetric triple. 
A cohomology class $[\alpha,\gamma]\in \cH(\fl,\theta_\fl,\fa)$ is called 
admissible\index{cohomology class!admissible}, 
if it is balanced 
and satisfies $(T_{2})$. We denote the set of all admissible quadratic cohomology 
classes by   
$\cH(\fl,\theta_\fl,\fa)_\sharp$ and its subset of all indecomposable admissible 
classes by $\cH(\fl,\theta_\fl,\fa)_0$. 
\end{de} 
By Proposition \ref{miau} admissible cohomology classes correspond to admissible 
quadratic extensions. 

Let $(\fl,\theta_{\fl})$ be a proper $\ZZ_2$-equivariant Lie algebra. As in Section 
\ref{S24} we consider the category $M^{ss}_{\fl,\theta_{\fl}}$ 
of semisimple orthogonal $(\fl,\theta_{\fl})$-modules and morphisms 
of pairs $((\fl,\theta_{\fl}),\fa_1)\rightarrow 
((\fl,\theta_{\fl}),\fa_2)$. The automorphism group $G_{\fl,\theta_{\fl},\fa}$ of 
an 
object $\fa$ 
of $M^{ss}_{\fl,\theta_{\fl}}$ 
acts on $\cH(\fl,\theta_\fl,\fa)_0$ and on $\cH(\fl,\theta_\fl,\fa)_\sharp$ from 
the 
right. 

Combining Proposition \ref{miau} with Theorem \ref{class} we arrive at the 
following 
classification scheme 
for symmetric triples. We prefer to formulate it for indecomposable symmetric 
triples. 
One gets the corresponding statement for general symmetric triples if one replaces 
$\cH(\fl,\theta_\fl,\fa)_0$ 
by $\cH(\fl,\theta_\fl,\fa)_\sharp$. 

 \begin{theo}[\cite{KO2}, Section 6]\label{symmcl} 
Let $\cL_p$ be a complete set of representatives of isomorphism classes of proper 
$\ZZ_2$-equivariant Lie algebras. 
For each $(\fl,\theta_{\fl})\in\cL_p$ we choose a complete set of representatives 
$\cA_{\fl,\theta_{\fl}}$ of isomorphism classes of objects in 
$M^{ss}_{\fl,\theta_{\fl}}$. 

Then there is a bijective map from the set of isomorphism classes of 
non-semisimple indecomposable symmetric triples   
 to the union of orbit spaces 
 \begin{equation}\label{Esymmcl} 
     \coprod_{(\fl,\theta_{\fl})\in\cL_p}\  \coprod_{\fa\in\cA_{\fl,\theta_{\fl}}}\   
 \cH(\fl,\theta_\fl,\fa)_0/ G_{\fl,\theta_{\fl},\fa}\ . 
 \end{equation} 
The inverse of this map sends the orbit of 
$[\alpha,\gamma]\in\cH(\fl,\theta_\fl,\fa)_0$ 
to the isomorphism class of $\fd_{\alpha,\gamma}(\fl,\theta_\fl,\fa)$. 
\end{theo} 
The theorem says in particular that the set 
$$\{ 
\fd_{\alpha,\gamma}(\fl,\theta_\fl,\fa)\:|  
\:(\fl,\theta_{\fl})\in\cL_p,\fa\in\cA_{\fl 
,\theta_{\fl}}, 
[\alpha,\gamma]\in\cH(\fl,\theta_\fl,\fa)_0\}$$ 
exhausts all isomorphism classes of non-semisimple indecomposable symmetric 
triples. 
As discussed at the end of Section \ref{S24} one does not need the whole set 
$\cL_p$. 
Only those $(\fl,\theta_\fl)$ 
such that the Lie algebra $\fl$ is admissible really occur. 

We will explain in Section \ref{S33} how Theorem \ref{symmcl} can be used in order 
to 
give a full classification 
of all symmetric triples of index at most $2$. Before doing this we want to discuss 
the implications of the theory of quadratic extensions for the geometry of 
symmetric 
spaces.

\subsection{The geometry of quadratic extensions} \label{S32}

Let $M$ be a pseudo-Riemannian symmetric space without semisimple 
local factors. As discussed in the last section, its (local) geometry
is completely determined by its symmetric triple and, moreover, this symmetric triple carries
the additional structure of an admissible quadratic extension. It is this structure that leads to
the classification scheme Theorem \ref{symmcl}. But what does this structure mean for the geometry
of $M$? This is the question we are going to discuss now. In particular, 
it will turn out that any  pseudo-Riemannian symmetric space without semisimple 
local factors
$M$ comes with a distinguished fibration $q: M\rightarrow N$ over an {\em 
affine} symmetric space $N$  such that
all fibres are flat.

Let us first recall the notion of an affine symmetric space\index{symmetric 
space!affine}. 
A connected manifold with connection $(M,\nabla)$ is
called an affine symmetric space if for each $x\in M$ there is an
involutive affine transformation $\theta_x$ of
$(M,\nabla)$ such that $x$ is an isolated fixed point of $\theta_x$. Note that $\theta_x$, if it exists,
is uniquely determined by $\nabla$.
Forgetting about the metric and only remembering the Levi-Civita connection we can consider
any pseudo-Riemannian symmetric space as an affine symmetric space. There are, however, many affine
symmetric spaces that do not admit any symmetric pseudo-Riemannian metric.
In exactly the same way as in the pseudo-Riemannian case one constructs 
the group of transvections\index{transvection group} of $(M,\nabla)$, which
acts transitively on $M$. Its Lie algebra comes with an involution but 
without scalar product. We will call this Lie algebra
with involution the symmetric pair\index{symmetric pair} of $M$. 
\begin{de}\label{sypa}
A symmetric pair is a proper $\ZZ_2$-equivariant Lie algebra $(\fg,\theta)$ satisfying
$\fz(\fg)\subset\fg_-$.
\end{de}
It follows from the effectivity of the action of the transvection group that the symmetric pair of $M$ is
indeed a symmetric pair in the sense of this definition.
Moreover, if $(\fg,\theta,\ip)$ is a symmetric triple, then $(\fg,\theta)$ is a symmetric pair. Indeed,
if $X\in \fz(\fg)$, then $ \langle X, [Y,Z] \rangle =\langle [X,Y], Z\rangle =0$ holds for all $Y,Z\in \fg_-$,
i.e., $X\in [\fg_-,\fg_-]^\perp=\fg_-$.
We have the following analogue of Proposition \ref{otto}.
\begin{pr}\label{mops} 
The assignment which sends each affine symmetric space to its symmetric pair induces a bijective map between affine
diffeomorphism classes
of simply connected affine symmetric spaces and isomorphism classes of symmetric pairs.
\end{pr}
Also the description of all affine symmetric spaces corresponding to a given symmetric pair proceeds in the
same way as in the pseudo-Riemannian case. In practice, affine symmetric spaces often arise as follows:
Let $\theta:G\rightarrow G$ be an involutive automorphism of a connected Lie group, and let $G_+\subset G$
be a closed $\theta$-stable subgroup having the same identity component as the fixed point group $G^\theta$.
Then $\theta$ induces an involutive diffeomorphism $\theta_{x_0}$ of $M:=G/G_+$, and there is a unique $G$-invariant connection
$\nabla$ on $M$ making $\theta_{x_0}$ affine. Then $(M,\nabla)$ is an affine symmetric space. Note that $G$
might be different from the transvection group of $M$. In fact, the transvection group is a subgroup of a
quotient of $G$.

If $q:M_1\rightarrow M_2$ is an affine map between two affine symmetric spaces, then 
$q\circ\theta_x=\theta_{q(x)}\circ q$ for all $x\in M_1$. If $q$ is surjective, then it follows that the fibres
$q^{-1}(m)$, $m\in M_2$, are totally geodesic (possibly disconnected) affine symmetric submanifolds of $M_1$ and
that $q$ induces a surjective homomorphism between the transvection groups of $M_1$ and $M_2$. 
A submanifold $Y$ of a pseudo-Riemannian manifold $(M,g)$ is called coisotropic, if 
$T_yY^{\perp_{g_y}}\subset T_yY$ for all $y\in Y$.
The following
notion is a geometric counterpart of the notion of a quadratic extension (in the context of symmetric triples).

\begin{de}
Let $N$ be an affine symmetric space. A special affine 
fibration\index{fibration!special affine} over $N$ is a surjective affine
map $q:M\rightarrow N$, where $M$ is a pseudo-Riemannian symmetric space and the fibres of $q$
are flat, coisotropic, and connected.
\end{de}

Let $N$ be an affine symmetric space with corresponding symmetric pair $(\fl,\theta_{\fl})$. Its cotangent bundle $T^*N$ carries the structure of a pseudo-Riemannian
symmetric space such that its symmetric triple is $\fd_{0,0}(\fl,\theta_{\fl},0)$.
Then the bundle projection $p:T^*N\rightarrow N$ is the simplest example of a special 
affine fibration over $N$. 

Recall that $\fd_{0,0}(\fl,\theta_\fl,0)$ is a quadratic extension of $(\fl,\theta_\fl)$. We now want
to construct special affine fibrations that correspond to quite general quadratic extensions in the same
way as $p:T^*N\rightarrow N$ corresponds to $\fd_{0,0}(\fl,\theta_\fl,0)$. 
Any proper $\ZZ_2$-equivariant Lie algebra $(\fl,\theta_\fl)$
gives rise to a symmetric pair $(\fl_0,\theta_{\fl_0})$, where 
$\fl_0:=\fl/(\fz(\fl)\cap\fl_+)$ and $\theta_{\fl_0}$ is
induced by $\theta_\fl$.

\begin{pr}\label{crux}
Let $(\fl,\theta_\fl)$ be a proper $\ZZ_2$-equivariant Lie algebra. Let $(\fg,\theta,\ip)$ be a symmetric
triple equipped with the structure $(\fg,\fri,i,p)$ of a quadratic extension of $(\fl,\theta_\fl)$ by
some orthogonal $(\fl,\theta_\fl)$-module. Let $M$ be a pseudo-Riemannian symmetric space with symmetric
triple $(\fg,\theta,\ip)$. We assume in addition that at least one of the following two conditions is satisfied:
\begin{enumerate}
\item[(a)] $M$ is simply connected.
\item[(b)] $\fz(\fg)\subset \fri^\perp$.
\end{enumerate}
Then there is an affine symmetric space $N$, unique up to isomorphism, with associated symmetric pair $(\fl_0,\theta_{\fl_0})$ and a unique special affine fibration
$q:M\rightarrow N$ such that 
\begin{equation}\label{pickel}
dQ_e=p_0\ .
\end{equation} 
Here $Q$ is the homomorphism of transvection groups induced by $q$ and
$p_0$ is the composition of natural maps $\fg\rightarrow\fg/\fri\stackrel{p}\longrightarrow \fl\rightarrow \fl_0$.
The symmetric space $N$ can be written as a homogeneous space $N=L/L_+$, 
where $L,L_+$ are certain
Lie groups with 
Lie algebras $\fl_+\subset \fl$. 
\end{pr}
If $N=L/L_+$, then the transvection group of $N$ equals $L_0=L/(Z(L)\cap L_+)$. Note that $Z(L)\cap L_+$ acts trivially on $N$. Thus the group $L$ (and the Lie algebra $\fl$) that
corresponds to the data of the quadratic extension arises as a central extension of the geometrically visible
transvection group $L_0$ (of the Lie algebra $\fl_0$).

The proof of the proposition will be given in the appendix. The idea behind it is very simple. Let $G$ be the transvection group of $M$, and let $J\subset G$ be the analytic subgroup corresponding to the ideal $\fri^\perp\subset\fg$. Then $J$ acts on $M$, and we would like to set $q$ to be the projection
onto the orbit space $N=J\backslash M$. That the orbits are flat and coisotropic is a simple consequence of the properties of $\fri^\perp$. The main problem is to show that the orbit space is a manifold (the orbits have
to be closed, in particular). It is this point, where we need Condition (a) or (b). Without these conditions, it is not difficult to construct
examples with non-closed $J$-orbits. For them the resulting geometric structure will be a foliation only, not a fibration. 

We are mainly interested in admissible, hence balanced quadratic extensions. They always satisfy Condition (b).
Indeed, using Equation (\ref{blub}) we find $\fri=\fri(\fg)\subset R(\fg)\subset \fg^\prime$. Forming
orthogonal complements yields $\fz(\fg)\subset \fri^\perp$. Thus admissible quadratic extensions give rise
to special affine fibrations.

\begin{co}
Let $M$ be a pseudo-Riemannian symmetric space without semisimple 
local factors. Let $(\fg,\theta,\ip)$
be the corresponding symmetric triple. Then 
the canonical quadratic extension associated with $(\fg,\theta,\ip)$ (in the sense of Corollary \ref{C1})
defines 
a special affine fibration $q: M\rightarrow N$ over an affine symmetric space $N$.
\end{co}
We call this fibration the {\em canonical 
fibration}\index{fibration!canonical} of $M$. Its base 
$N$ is an important  invariant of the pseudo-Riemannian space $M$. Since 
there are a lot of non-admissible Lie algebras (see the end of Section~\ref{S24})
not all affine symmetric spaces $N$ can appear as the base of the canonical fibration of some pseudo-Riemannian
symmetric space $M$. 

Let $q:M\rightarrow N$ be a special affine fibration. Let us collect some of its
basic properties. 
\begin{enumerate}
\item $q:M\rightarrow N$ is locally trivial, i.e., it is a fibre bundle
with flat affine symmetric fibres. More precisely, there exist $k,l\in\NN_0$, an open covering $\{U_i\}$ of $N$,
and diffeomorphisms
$\Phi_i:q^{-1}(U_i)\rightarrow U_i\times\RR^k\times (S^1)^l$ 
such that for all $x\in U_i$ the restriction of $\Phi_i$ to $q^{-1}(x)$ is an affine
diffeomorphism from the fibre onto 
$\{x\}\times\RR^k\times (S^1)^l$. If $l=0$ we call the fibration $q$ very nice. Note that a very nice
fibration is not a vector bundle, in general. There is no distinguished zero section. 
\item $M$ is simply connected if and only if $N$ is simply connected and $q$ is very nice.
\item The fibres of $q$ are foliated by the null spaces of the restriction of the metric to the fibre. The leaves are
totally geodesic subspaces of dimension $n=\dim N$. If the leaves are closed we call $q$ nice. Of course, very nice
implies nice.
\item We look at the cotangent bundle $T^*N$ as bundle of abelian groups. There is a natural action $\sigma$ 
of the bundle $T^*N$ on the bundle 
$q: M\rightarrow N$
by translations on the fibres. Its orbits are precisely the null leaves described in 3.
\item If $q$ is nice, then the space of null leaves (= the orbit space of the $T^*N$-action) is
itself an affine symmetric space which is fibred over $N$. We obtain a factorisation $q=s\circ r$,
where $r: M\rightarrow P$ and $s: P\rightarrow N$ are affine maps which are fibre bundles
with flat affine symmetric fibres. Moreover, the fibres of $s$ come with a non-degenerate metric.
They are quotients of pseudo-Euclidean spaces by a (often trivial) discrete group of translations.
\end{enumerate}
One should observe the analogy to the properties of quadratic extensions.
We summarize the structure of a nice special affine fibration by the following diagram.
$$
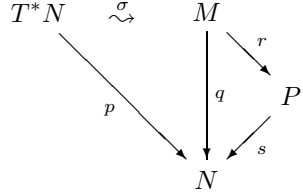
\begin{diagram}
\node{T^{*}N}\arrow[2]{se,b}{p} \node{\stackrel{\sigma}{\leadsto}} \node{M}\arrow[2]{s,r}{q} 
\arrow{se,t}{r}\\
\node{} \node{}  \node{} \node{P}\arrow{sw,b}{s}\\
\node{} \node{}  \node{N} 
\end{diagram}$$
Using Lemma \ref{hempel} it can be shown that the canonical
fibration of every indecomposable symmetric space $M$ is nice (in fact, the absence of flat
local factors that are not global factors is sufficient). Thus any such space $M$ carries a canonical structure of the kind indicated by the diagram. The data $(\fl,\theta_\fl,\fa,[\alpha,\gamma])$ appearing in the classification
scheme Thm.~\ref{symmcl} should be regarded as a complete set of (local) invariants describing this structure.
It would be an interesting project to work out the precise geometric meaning of each of these invariants.

We conclude this section with a certain converse of Proposition \ref{crux} saying that all special
affine fibrations come from quadratic extensions. Let us denote the special affine fibration constructed
in Proposition \ref{crux} by $q(M,\fri,i,p)$.

\begin{pr}\label{xurc}
Let $q:M\rightarrow N$ be a special affine fibration. Then there exists a structure of a quadratic
extension $(\fg,\fri,i,p)$ on the symmetric triple
$(\fg,\theta,\ip)$ of $M$ such that $q=q(M,\fri,i,p)$.
\end{pr}
We remark that the quadratic extension $(\fg,\fri,i,p)$ is not 
completely determined
by $q$. What is uniquely determined is $\fri_-\subset \fg_-$.
Indeed, if we identify $\fg_-=T_{x_0}M$, then 
$(\fri_-)^\perp\subset\fg_{-}$ has to be the tangent space to the fibre of $q$. We then have to choose
$\fri_+\subset\fg_+$ subject to the conditions
\begin{itemize}
\item[(a)] $\fri=\fri_+\oplus \fri_-\subset \fg$ is an isotropic ideal,
\item[(b)] $\fri^\perp/\fri$ is abelian.
\end{itemize}
The proposition says that such a choice is always possible.
Indeed, one can show that
$$ \fri_+:=\{X\in [\fg_-,(\fri_-)^\perp]\:|\: [X,(\fri_-)^\perp]=0\}\ .$$
always satisfies (a) and (b). Nevertheless, $\fri_+$ is not uniquely determined by $\fri_-$, (a), and (b), in general.

\subsection{Symmetric spaces of index one and two}\label{S33}
\index{symmetric space!Lorentzian}
In this 
section we will comment on some classification results for symmetric 
triples of small index. First we want to reformulate the classification 
of indecomposable non-semisimple Lorentzian symmetric 
triples by Cahen and Wallach \cite{CW} in terms of quadratic extensions. 
Indeed, this result follows easily from Theorem~\ref{symmcl}. To see 
this we just have to check which elements of (\ref{Esymmcl}) 
correspond to Lorentzian symmetric triples. Clearly, 
$\fd_{\alpha,\gamma}(\fl,\theta_{\fl},\fa)$ is Lorentzian if and only 
if $a+\dim \fl=1$, where $a$ is the index of $\ip_{\fa}$ restricted 
to $\fa_{-}$. Hence, either $\fl=0$ and 
$\fa=\fa_{-}=\RR^{1,0}$ or $(\fl,\theta_{\fl})=(\RR,-\Id)$ and $\ip_{\fa}$ 
restricted 
to $\fa_{-}$ is positive definite. The first case is trivial. Let us 
consider the second one, i.e.~suppose $(\fl,\theta_{\fl})=(\RR,-\Id)$. 
Take $\fa=\RR^{p,p}\oplus \RR^{0,2q}$, $p,q\ge0$. Let 
$e_1,\ldots,e_{2p}$ be the standard basis of $\RR^{p,p}$ and let 
$e'_1,\ldots,e'_{2q}$ be the standard basis of $\RR^{0,2q}$. We define 
an involution $\theta_{\fa}$ on $\fa$ by 
$\fa_{+}=\Span\{e_1,\ldots,e_{p},e'_1,\ldots,e'_{q}\}$ and 
$\fa_{-}=\Span\{e_{p+1},\ldots,e_{2p},e'_{q+1},\ldots,e'_{2q}\}$. 
For  $\lambda=(\lambda^{1},\ldots,\lambda^{p})\in (\fl^{*})^{p}\cong 
\RR^{p}$ and $\mu=(\mu^{1},\ldots,\mu^{q})\in (\fl^{*})^{q}\cong\RR^{q}$ we define a 
representation $\rho_{\lambda,\mu}$ of $\fl$ on $\fa$ by 
\begin{eqnarray} 
    \rho_{\lambda,\mu}(L)(e_{i})=\lambda^{i}(L)e_{i+p},&& 
\rho_{\lambda,\mu}(L)(e_{i+p})=\lambda^{i}(L)e_{i},\label{lm1}\\ 
   \rho_{\lambda,\mu}(L)(e'_{j})=\mu^{j}(L)e'_{j+q},&& 
\rho_{\lambda,\mu}(L)(e'_{j+q})=-\mu^{j}(L)e'_{j}\label{lm2} 
\end{eqnarray} 
for $L\in\fl$, $i=1,\ldots,p$ and $j=1,\ldots,q$. Then 
$\fa_{\lambda,\mu}:=(\rho_{\lambda,\mu},\fa)$ 
is an orthogonal $(\fl,\theta_{\fl})$-module and we can define 
$$\fd(p,q,\lambda,\mu):=\fd_{0,0}(\fl,\theta_{\fl},\fa_{\lambda,\mu}).$$

It is not hard to prove that every semisimple orthogonal $(\RR,-\Id)$-module for 
which the $-1$-eigenspace of the involution is positive definite is of 
the kind defined above. Since $\fl$ is one-dimensional we have 
$\cH(\fl,\theta_\fl,\fa)=\{(0,0)\}$ and indecomposability and 
admissibility conditions are easy to handle. Now Theorem~\ref{symmcl} 
gives: 
\begin{theo}[Cahen/Wallach \cite{CW}] \label{TCW}
    Every indecomposable non-semisimple Lo\-rentzian symmetric triple 
    is either one-dimensional or isomorphic to exactly one of the 
    symmetric triples $\fd(p,q,\lambda,\mu)$, $p,q\ge0$, $p+q>0$, $(\lambda,\mu)\in M_{p,q}$, where 
    $$M_{p,q}:=\left\{(\lambda,\mu)\in \RR^{p}\oplus \RR^{q} \left|
    \begin{array}{l}
0<\lambda^{1}\le \lambda^{2}\le \ldots \le 
    \lambda^{p},\ 0<\mu^{1}\le \mu^{2}\le \ldots \le \mu^{q},\\[0.5ex]
    \lambda^{1}=1\ \mbox{\rm if } p>0,\ \mu^{1}=1\ \mbox{\rm else }
    \end{array}\right. \hspace{-5pt}\right\}.$$ 
\end{theo}     
Cahen and Wallach \cite{CW} constructed explicit models for all simply connected Lorentzian symmetric spaces.  Let us describe the simply connected Lorentzian symmetric space $M$ associated with the symmetric triple $\fd:=\fd(p,q,\lambda,\mu)$. Since $\fd$ is isomorphic to the semidirect product of a Heisenberg algebra by $\RR$ it is not hard to see that the simply connected group $G$ with Lie algebra $\fd$ is isomorphic to $\fd=\fl^*\oplus\fa\oplus\fl$ with group multiplication
$$(Z,A,L)(\bar Z,\bar A,\bar L)=(Z+\bar Z+ \frac12 \left[e^{-\ad \bar L}(A),\bar A\right], e^{-\ad \bar L}(A)+\bar A,L+\bar L).$$ 
Here $\lb$ and $\ad$ are the operations in $\fd$.
The analytic subgroup $G_+$ of $G$ with Lie algebra $\fd_+$ then equals $\fa_+$. The projection $G\rightarrow G/G_+$ has the  global section
$$s: G/G_+\longrightarrow G,\quad (Z,A_++A_-,L)\cdot G_+\longmapsto (Z+\textstyle{\frac12} [A_+,A_-], A_-,L),$$
where $A_+\in\fa_+$ and $A_-\in\fa_-$. In particular, we can identify $G/G_+$ with $\fd_-$ (as a set). Let   $(z,a_1,\dots,a_p,a_1',\dots,a_q',l)$ denote the coordinates of a vector in $\fd_-\cong \fl^*\oplus\fa_-\oplus\fl\cong\RR\oplus\RR^{p+q}\oplus \RR$. In these coordinates, the metric on $G/G_+$ determined by the symmetric triple $\fd$ is given by
$$2dzdl+ \sum_{i=1}^p da_i^2 + \sum_{j=1}^q {da'_j}^2 + \left(  \sum_{i=1}^p\lambda_i^2a_i^2 -  \sum_{j=1}^q\mu_j^2a_j'^2 \right)dl^2.$$

Now let us turn to symmetric triples of index two. First classification 
results for this case were already obtained by Cahen and Parker in 
\cite{CP1} and \cite{CP}. In \cite{CP1} symmetric spaces of index two 
with solvable transvection group were studied. However, the 
classification presented there turned out to be incomplete. In his 
diploma thesis \cite{N} Th.~Neukirchner elaborated the claimed results, found the 
gaps 
and gave a revised classification of indecomposable symmetric spaces of 
index two that 
have a solvable transvection group. We tried to reproduce this result 
using our classification scheme and observed that also Neukirchner's 
classification is not quite correct. Besides minor errors a series of 
spaces is missing and some of the normal forms contain too much 
parameters. In \cite{KO2} we give the corrected classification result. 
We don't want to recall this result to its full extent here. However, let us give a rough classification of indecomposable non-semisimple symmetric 
triples of signature  $(2,n)$, which shows how our classification scheme works. 

Following Thm.~\ref{symmcl} we have first to determine all proper 
$\ZZ_2$-equivariant Lie algebras $(\fl,\theta_\fl)$ for which $\dim \fl_-\le 2$ 
holds. For $\fl$ besides $\RR, \RR^2, \su(2)$, $\fsl(2,\RR)$ and $\fh(1)$ the 
following Lie algebras appear: 
$$\fn(2)=\{[X,Y]=Z,\,[X,Z]=-Y\},\ \fr_{3,-1}=\{[X,Y]=Z,\,[X,Z]=Y\}.$$ 
The involutions $\theta_\fl$ are described in Table 1. The table also lists the associated simply connected 
affine symmetric spaces $N=L/L_+$\index{fibration!canonical} for all $(\fl,\theta_{\fl})$ (see Section \ref{S32} for the geometric meaning of $N$). 
The Lie algebra $\fn(2)$ is isomorphic to $\so(2)\ltimes\RR^2$, 
$(\fn(2),\theta_\fl)$  is a symmetric pair and the associated simply-connected 
affine symmetric space equals the universal covering $\widetilde{\cL(2)}$ 
of the space of (affine) lines in $\RR^2$. Analogously, 
$\fr_{3,-1}$ is isomorphic to $\so(1,1)\ltimes\RR^{1,1}$, 
$(\fr_{3,-1},\theta_\fl)$ is also a symmetric pair and the associated 
simply-connected affine symmetric space equals the space $\cL(1,1)$ of time-like 
(affine) lines in $\RR^{1,1}$. Note that we have two non-conjugate 
involutions on $\fsl(2,\RR)$. The two associated simply-connected 
symmetric spaces are the hyperbolic plane $H^{2}$ and the universal covering 
$\widetilde{S^{1,1}}$ of the unit sphere 
$S^{1,1}:=\{x\in\RR^{1,2}\mid \langle x,x \rangle_{1,2}=1\}$ in $\RR^{1,2}$. 

For all these $(\fl,\theta_\fl)$ we have to determine a set 
$\cA_{\fl,\theta_\fl}$ as described in Thm.~\ref{symmcl}. More exactly, 
since we want to classify only indecomposable symmetric spaces of signature $(2,n)$ we may restrict 
ourselves to the subset $\cA^{n}_ {\fl,\theta_\fl}\subset 
\cA_{\fl,\theta_\fl}$ of all $\fa\in\cA_{\fl,\theta_\fl}$ for 
which $\cH(\fl,\theta_\fl,\fa)_0$ is not empty and for which 
$(2,n)=\sgn\fa_-+(\dim\fl_-,\dim\fl_-)$. For all $(\fl,\theta_\fl)$ listed 
above $\cA^n_ {\fl,\theta_\fl}$ consists of finitely many families 
$(\rho_1,\fa_1),\dots, (\rho_k,\fa_k)$ of 
$(\fl,\theta_\fl)$-represen\-tations, where for each family 
$(\rho_j,\fa_j)$ the space $\fa_j$ is fixed and $\rho_j$ depends in a 
certain sense on 
$d_{\rho_j}$ continuous parameters for some $d_{\rho_j}\in\NN_0$. We list 
the spaces $\fa_1,\dots,\fa_k$ explicitly in the table, where we use the 
following notation: $\fa_+^{p,q}:=(\RR^{p,q},\Id)$,   
$\fa_-^{p,q}:=(\RR^{p,q},-\Id)$. For the families $\rho_1,\dots, 
\rho_{k}$ 
we omit a detailed description, we only give the number $d_{\rho_j}$ of 
parameters for each $\rho_j$. Next we have to compute 
$\cH(\fl,\theta_\fl,\fa_j)/G_{\fl,\theta_\fl,\fa_j}$ for each element of 
the family $(\rho_j,\fa_j)$, $j=1,\dots,k$. We do not give the result in 
the table, either. However, we give the number $d_{\cal H}$ of 
continuous parameters on which 
$\cH(\fl,\theta_\fl,\fa_j)/G_{\fl,\theta_\fl,\fa_j}$ depends for a 
generic element in the family $(\rho_j,\fa_{j})$. Here $d_\rho=0$ or $d_{\cal H}=0$ means that the corresponding set is discrete (in fact, it is finite). Note that the 
values of $d_\rho$ and $d_{\cal H}$ for $N=\RR^{2}$ in the table are 
not correct for small $n$, namely for those $n$ for which the value of 
$d_\rho$ given in the table would be negative. We don't consider these special 
cases in the table. 
\begin{table} 
\begin{tabular}{|c|ll|c|c|} 
\hline 
&&&&\\[-2ex] 
$N$ &\multicolumn{2}{|c|}{ $\fl,\ \cA^{n}_ {\fl,\theta_\fl}$} & 
$d_{\rho}$ & $d_{\cal H}$\\[1ex] 
\hline \hline 
&&&&\\[-2.0ex] 
$\RR^{1}$ &\multicolumn{2}{|l|}{$\fl=\fl_{-}=\RR$} && \\[0.5ex] 
&1.& $(\rho_{1,p},\fa_{+}^{p,n-p}\oplus \fa_{-}^{1,n-1}),\ 
p=0,\ldots,n-1$ & $n-1$ &0 
\\ [0.5ex] 
&2.& $(\rho_{2,p},\fa_{+}^{p,n-p}\oplus \fa_{-}^{1,n-1}),\ p=1,\ldots,n$ & $n-1$& 0 
\\ [0.5ex] 
&3.& $(\rho_{3,p},\fa_{+}^{p,n-p}\oplus \fa_{-}^{1,n-1}),\   
p=1,\ldots,n-1$ & $n-1$& 0 
\\ [0.5ex] 
\hline 
&&&&\\[-2.0ex] 
$\RR^{2}$ &\multicolumn{2}{|l|}{$\fl=\fl_{-}=\RR^{2}$} && \\[0.5ex] 
&1.& $(\rho_{1,p},\fa_{+}^{p,n-2-p}\oplus \fa_{-}^{0,n-2}),\ 
p=0,\ldots,n-2$ & $2n-8$& 0\\ 
&&(occurs only for $n\ge5$)&& 
\\ [0.5ex] 
&2.& $(\rho_{2,p},\fa_{+}^{p,n-1-p}\oplus \fa_{-}^{0,n-2}),\ 
p=0,\ldots,n-2$ & 
$2n-8$& 1 
\\ [0.5ex] 
&3.& $(\rho_{3,p},\fa_{+}^{p,n-1-p}\oplus \fa_{-}^{0,n-2}),\ 
p=1,\ldots,n-1$ & 
$2n-8$& 1 
\\ [0.5ex] 
\cline{2-5} 
&&&&\\[-2.0ex] 
 &\multicolumn{2}{|l|}{$\fl=\fh(1)=\{[X,Y]=Z\},\, 
\fl_{+}=\RR\hspace{-1pt}\cdot\hspace{-2pt}Z,\, \fl_{-}=\Span\{X,Y\}$} && \\[0.5ex] 
&1.& $(\rho_{1,p},\fa_{+}^{p,n-3-p}\oplus \fa_{-}^{0,n-2}),\ 
p=0,\ldots,n-3$ & 
$2n-10$ &2 
\\ [0.5ex] 
&2.& $(\rho_{2,p},\fa_{+}^{p,n-4-p}\oplus \fa_{-}^{0,n-2}),\ 
p=0,\ldots,n-4$ & 
$2n-12$ &3 
\\ [0.5ex] 
\hline 
&&&&\\[-2.0ex] 
$\widetilde{\cL(2)}$ &\multicolumn{2}{|l|}{$\fl=\fn(2),\ 
\fl_{+}=\RR\cdot\hspace{-1pt} Z,\ 
\fl_{-}=\Span\{X,Y\}$} && \\[0.5ex] 
&1.& $(\rho_{1,p},\fa_{+}^{p,n-2-p}\oplus \fa_{-}^{0,n-2}),\ 
p=0,\ldots,n-2$ & $n-2$& 0 
\\ [0.5ex] 
&2.& $(\rho_{2,p},\fa_{+}^{p,n-3-p}\oplus \fa_{-}^{0,n-2}),\ 
p=0,\ldots,n-3$ & 
$n-3$& 0 
\\ [0.5ex] 
&3.& $(\rho_{3,p},\fa_{+}^{p,n-2-p}\oplus \fa_{-}^{0,n-2}),\ 
p=0,\ldots,n-3$ & 
$n-3$& 0 
\\ [0.5ex] 
&4.& $(\rho_{4,p},\fa_{+}^{p,n-3-p}\oplus \fa_{-}^{0,n-2}),\ 
p=0,\ldots,n-4$ & 
$n-4$& 1 
\\ [0.5ex] 
\hline 
&&&&\\[-2.0ex] 
$\cL(1,1)$ &\multicolumn{2}{|l|}{$\fl=\fr_{3,-1},\ \fl_{+}=\RR\cdot Z,\ 
\fl_{-}=\Span\{X,Y\}$} && \\[0.5ex] 
&1.& $(\rho_{1,p},\fa_{+}^{p,n-2-p}\oplus \fa_{-}^{0,n-2}),\ 
p=0,\ldots,n-2$ & $n-2$& 0 
\\ [0.5ex] 
&2.& $(\rho_{2,p},\fa_{+}^{p,n-3-p}\oplus \fa_{-}^{0,n-2}),\ 
p=0,\ldots,n-3$ & 
$n-3$& 0 
\\ [0.5ex] 
&3.& $(\rho_{3,p},\fa_{+}^{p,n-2-p}\oplus \fa_{-}^{0,n-2}),\ 
p=0,\ldots,n-3$ & 
$n-3$& 0 
\\ [0.5ex] 
&4.& $(\rho_{4,p},\fa_{+}^{p,n-3-p}\oplus \fa_{-}^{0,n-2}),\ 
p=0,\ldots,n-4$ & 
$n-4$& 1 
\\ [0.5ex] 
\hline 
&&&&\\[-2.0ex] 
$\widetilde{S^{1,1}}$ 
&\multicolumn{2}{|l|}{$\fl=\fsl(2,\RR)=\{[H,X]=2Y, [H,Y]=2X,
$}&&\\ 
&\multicolumn{2}{|l|}{$[X,Y]=2H\},\ \fl_{+}=\RR\cdot H,\ 
\fl_{-}=\Span\{X,Y\}$} && \\[0.3ex] 
&\multicolumn{2}{|l|}{$(\rho,\fa_{+}^{n-2,n-2}\oplus \fa_{-}^{0,n-2})$}& $0$& 1 
\\ [0.5ex] 
\hline 
&&&&\\[-2.0ex] 
$H^2$ &\multicolumn{2}{|l|}{$\fl=\fsl(2,\RR),\ \fl_{+}=\RR\cdot X,\ 
\fl_{-}=\Span\{H,Y\}$} && \\[0.5ex] 
&\multicolumn{2}{|l|}{$(\rho_{klm}^{pqr},\fa_{+}^{n-2+p-q,0}\oplus 
\fa_{-}^{0,n-2}),\ k\in \NN^p,\, l\in \NN^q,\, m\in \NN^r,$}& $0$& 1\\[0.3ex] 
&\multicolumn{2}{|l|}{$k_1\le k_2\le\dots\le k_p,\ l_1\le\dots,\le l_q,\ m_1\le 
\dots \le m_r$} &&\\ 
&\multicolumn{2}{|l|}{$|k|+|l|+2|m|+q=n-2$} && 
\\ [0.5ex] 
\hline 
&&&&\\[-2.0ex] 
$S^2$ &\multicolumn{2}{|l|}{$\fl=\su(2)=\{[H,X]=2Y,[H,Y]=-2X, $}&&\\ 
&\multicolumn{2}{|l|}{$[X,Y]=2H\},\ \fl_{+}=\RR\cdot H,\ 
 \fl_{-}=\Span\{X,Y\}$} && \\[0.5ex] 
&\multicolumn{2}{|l|}{$(\rho_{klm}^{pqr},\fa_{+}^{0,n-2+p-q}\oplus 
\fa_{-}^{0,n-2}),$}& $0$& 1\\[0.3ex] 
&\multicolumn{2}{|l|}{$p,q,r,k,l,m$ as in the previous entry} &&\\ [0.5ex] 
\hline 
\end{tabular} 
\\[0.5ex] 
\begin{center} 
    {{\bf  Table 1.} Non-semisimple indecomposable symmetric triples of signature $(2,n)$} 
\end{center}     
\end{table}

In order to not present only vague data here let us study one case in more 
detail. We consider the case, where $\fl$ is the Heisenberg algebra 
$\fh(1)=\{[X,Y]=Z\}$. In Example~\ref{Eh1} we computed $\cH(\fl,\fa)_b$ 
for $\fl=\fh(1)$ and  any semisimple orthogonal $\fh(1)$-modul $\fa$. We identified   
$\cH(\fl,\fa)_b$ with a certain subset $Z_{\fl,b}$ of $C^2(\fl,\fa)$. Now 
let $\fa$ be a semisimple orthogonal $(\fl,\theta_\fl)$-module for 
$\theta_\fl$ given by $\fl_+=\RR\cdot Z,\,\fl_-=\Span\{X,Y\}$. Then 
$Z_{\fl,b}$ is invariant under the morphism of pairs 
$(\theta_\fl,\theta_\fa)$. Using Prop.~\ref{inva} we see that 
$\cH(\fl,\theta_\fl,\fa)_b$ corresponds bijectively to $\{\alpha\in Z_{\fl,b}\mid 
\alpha(Z,\fl)\subset\fa_-^\fl\}$. Furthermore, since $\fl$ is indecomposable, 
$\cH(\fl,\theta_\fl,\fa)_0$ corresponds bijectively to $\{\alpha\in Z_{\fl,b}\mid 
\alpha(Z,\fl)=\fa_-^\fl\}$. Summarising we see that 
\begin{eqnarray}\label{Eacht} 
Z_{\fl,0}:=\{\alpha\in C^2(\fl,\fa)\mid 
\alpha(X,Y)=0,\alpha(Z,\fl)=\fa_-^\fl\not=0\}& \longrightarrow& 
\cH(\fl,\theta_\fl,\fa)_0\\ 
\alpha &\longmapsto& [\alpha,0]\nonumber 
\end{eqnarray} 
is a bijection. Let us now determine a suitable set 
$\cA^n_{\fl,\theta_\fl}$. Let $(\rho,\fa)$ be an orthogonal 
$(\fl,\theta_\fl)$-module such that $\cH(\fl,\theta_\fl,\fa)_0\not=\emptyset$. 
In particular, $\rho$ is semisimple, which implies 
$\rho(Z)=0$ since $R(\fl)=\RR\cdot Z$. Hence, $\rho$ can be considered 
as a semisimple representation of the abelian Lie algebra $\fl_-$, which 
is determined by its weights. Moreover, we know from (\ref{Eacht}) that 
$\fa^\fl=\fa_-^{0,1}$ or  $\fa^\fl=\fa_-^{0,2}$.  For   
$\lambda=(\lambda^1,\dots,\lambda^p)\in  (\fl_-^*\setminus 0)^p$ and 
$\mu=(\mu^1,\dots,\mu^q) \in(\fl_-^*\setminus 0)^q$ we define a 
representation  $\rho_{\lambda,\mu}$ of $\fl_-$ on 
$\fa_+^{p,q}\oplus\fa_-^{0,p+q}$ by  
(\ref{lm1}) and (\ref{lm2}),  where we identify 
$\fa_+^{p,q}\oplus\fa_-^{0,p+q}$ with $\RR^{p,p}\oplus\RR^{0,2q}$. Now 
we define   
\begin{eqnarray*} 
    \fa_{1,\lambda \mu} 
&:=&(\rho_{\lambda,\mu}\oplus\rho_0,\fa_+^{p,q} 
    \oplus\fa_-^{0,p+q}\oplus \fa_-^{0,1}),\\ 
    \fa_{2,\lambda \mu}&:=&(\rho_{\lambda,\mu}\oplus\rho_0,\fa_+^{p,q}   
    \oplus\fa_-^{0,p+q}\oplus  \fa_-^{0,2}), 
\end{eqnarray*}     
where $\rho_0$ denotes the trivial 
representation   of dimension one and two, respectively. It is easy to 
see that $(\rho,\fa)$  is equivalent to one of these representations. Next 
we have to decide  which of these representations are isomorphic as objects 
of $M^{ss}_{\fl,\theta_{\fl}}$. The automorphism group of $(\fl,\theta_\fl)$ equals 
$$\Aut(\fl,\theta_\fl)=\Big\{ 
    \small{ 
\left( \begin{array}{cc} 
A&0\\ 
0&u 
\end{array}\right) 
} 
\ \Big|\ A\in GL(2,\RR),\ \det A=u\ 
\Big\}, $$ 
where the automorphisms are written with respect to the basis $X,Y,Z$. 
Let $\frak S_p$ denote the symmetric group of degree $p$. We define an action of $(\frak S_p\ltimes (\ZZ_2)^p)\times\GL(2,\RR)$ on $(\fl^*_-\setminus 0)^p$ by 
$(\lambda_{1},\ldots,\lambda_{p})\cdot \sigma\cdot\eps\cdot A=(\eps_{1}A^{*}\lambda_{\sigma(1)}, 
\ldots,\eps_{p}A^{*}\lambda_{\sigma(p)})$ for $(\lambda_{1},\ldots,\lambda_{p})\in (\fl^*_-\setminus 0)^p$, $\sigma\in {\mathfrak S}_{p}$, $(\eps_1,\ldots,\eps_{p})\in(\ZZ_2)^{p}$ and $A\in 
GL(2,\RR)$. We define the sets
$$\Lambda_p:=(\fl^*_-\setminus 0)^p/(\frak S_p\ltimes (\ZZ_2)^p),\quad \Lambda_{p,q}=\Lambda_p \times \Lambda_q.$$
Denoting the isomorphy relation in $M^{ss}_{\fl,\theta_{\fl}}$ by 
$\cong$ we obtain $\fa_{1,\lambda \mu}\not\cong \fa_{2,\lambda' 
    \mu'}$ and 
$$\fa_{i,\lambda \mu}\cong\fa_{i,\lambda' \mu'} 
    \ \Leftrightarrow \ p=p', q=q' \mbox{ and } [\lambda,\mu]= [\lambda',\mu']\in \Lambda_{p,q}/\GL(2,\RR) $$ 
for $i=1,2$. 
Hence, the elements of the following $2n-5$ families of $(\fl,\theta_\fl)$-modules 
constitute a suitable set $\cA^{n}_ {\fl,\theta_\fl}$: 
\begin{eqnarray*} 
    (\rho_{1,p},\fa_{+}^{p,n-3-p}\oplus 
\fa_{-}^{0,n-2})&:=&\{\fa_{1,\lambda\mu} \mid 
[\lambda,\mu]\in \Lambda_{p,n-3-p}/\GL(2,\RR) \},\\ 
(\rho_{2,\tilde p},\fa_{+}^{\tilde p,n-4-\tilde p}\oplus 
\fa_{-}^{0,n-2})&:=&\{\fa_{2,\lambda\mu}\mid 
[\lambda,\mu]\in \Lambda_{\tilde p,n-4-\tilde p}/\GL(2,\RR) \},\end{eqnarray*} 
where $p=0,\ldots,n-3$, $\tilde p=0,\ldots,n-4$.
In particular, for each family $\rho_{1,p}$, $p=0,\ldots,n-3$, we have 
$d_{\rho}=2n-10$ if $n\ge 5$ and $d_\rho=0$ for $n=3,4$. Similarly, 
for each $\rho_{2,\tilde p}$, $\tilde p=0,\ldots,n-4$, we have 
$d_{\rho}=2n-12$ if $n\ge6$ and $d_{\rho}=0$ for $n=4,5$. 

Finally, let us determine 
$\cH(\fl,\theta_{\fl},\fa)/G_{\fl,\theta_{\fl},\fa}$ 
for  $\fa\in\cA^n_{\fl,\theta_\fl}$. Take $\fa=\fa_{1,\lambda\mu}$. If $n\ge 5$, then 
in the 
generic case $(\lambda(X),\mu(X))\in\RR^{n-3}$ and 
$(\lambda(Y),\mu(Y))\in\RR^{n-3}$ are linearly independent and we get 
$\cH(\fl,\theta_{\fl},\fa)/G_{\fl,\theta_{\fl},\fa}\cong Z_{\fl,0}/\ZZ_2$. 
Hence, $d_{\cal H}=2$. In the non-generic case 
$G_{\fl,\theta_{\fl},\fa}$ becomes larger and we get $d_{\cal H}=0$. 
For $n=3,4$ we have $d_{\cal H}=0$. 
Now take $\fa=\fa_{2,\lambda\mu}$. If $n\ge 6$, then in the 
generic case $(\lambda(X),\mu(X))\in\RR^{n-4}$ and 
$(\lambda(Y),\mu(Y))\in\RR^{n-4}$ are linearly independent and we have 
$\cH(\fl,\theta_{\fl},\fa)/G_{\fl,\theta_{\fl},\fa}\cong 
Z_{\fl,0}/O(2)$. 
Hence, $d_{\cal H}=3$. In the non-generic case  we get $d_{\cal H}=1$. 
For $n=5$ we have $d_{\cal H}=1$ and for $n=4$ we get $d_{\cal H}=0$.
\index{symmetric space!pseudo-Riemannian|)}

\section{Special geometric structures on symmetric\\ spaces}

\subsection{Examples of geometric structures}\label{S41} 

We are now going to discuss pseudo-Riemannian symmetric spaces that are equipped with
certain geometric structures coming from complex and quaternionic geometry.

Let $(M,g)$ be a pseudo-Riemannian manifold. The Levi-Civita connection induces a connection on the bundle $\so(TM)\subset\End(TM)$ of endomorphisms
that are skew-symmetric w.r.t.~$g$. 
A K\"ahler structure on $(M,g)$ is a parallel section $I$ of $\so(TM)$ satisfying $I^2=-\Id_{TM}$.
In particular, $I$ is an integrable almost complex structure, and thus 
$(M,I)$ is a complex manifold. A {\em pseudo-Hermitian symmetric
space\index{symmetric space!pseudo-Hermitian}} is a tuple $(M,g,I)$, where $(M,g)$ is a pseudo-Riemannian symmetric space and $I$ is a K\"ahler structure on $(M,g)$. A hyper-K\"ahler
structure on $(M,g)$ is a pair of K\"ahler structures $(I,J)$ satisfying $IJ=-JI$.  A quaternionic K\"ahler structure
on $M$ arises if we weaken the parallelity conditions on $I,J$: it consists of a 3-dimensional parallel subbundle $E\subset\so(TM)$ that can be locally spanned by almost complex structures $I$, $J$, and $K:=IJ=-JI$. We have the
corresponding notions of a {\em hyper-K\"ahler symmetric space}\index{symmetric 
space!hyper-K\"ahler} $(M,g,I,J)$ and a {\em quaternionic K\"ahler symmetric 
space}\index{symmetric space!quaternionic K\"ahler} $(M,g,E)$. In particular, hyper-K\"ahler symmetric spaces form a subclass of all quaternionic K\"ahler symmetric spaces.

In the pseudo-Riemannian world all these structures have their ``para''-versions.
If we replace the condition $I^2=-\Id_{TM}$ for a K\"ahler structure by 
$I^2=\Id_{TM}$ we are lead to the notion of a para-K\"ahler
structure. A para-K\"ahler structure on $(M,g)$ is equivalent to a parallel splitting $TM=TM_+\oplus TM_-$ into totally isotropic
subbundles. Thus para-K\"ahler structures can exist for metrics of neutral signature $(m,m)$, only.
A pair $(I,J)$, where $I$ is a K\"ahler structure and $J$ is a para-K\"ahler structure such that $IJ=-JI$, is called
a hypersymplectic structure (sometimes also para-hyper-K\"ahler structure). Note that then $K:=IJ=-JI$
is a second para-K\"ahler structure on $M$. A parallel subbundle $E\subset\so(TM)$ locally spanned by
(not necessarily parallel) sections $I,J,K$ of this kind is called a para-quaternionic K\"ahler structure.
There are the corresponding notions of {\em 
para-Hermitian}\index{symmetric space!para-Hermitian}, {\em 
hypersymplectic}\index{symmetric space!hypersymplectic}, and {\em 
para-quaternionic K\"ahler} symmetric spaces\index{symmetric 
space!para-quaternionic K\"ahler}.

Recall the notion of an $(\fh,K)$-module $(V,\Phi_V)$ from Section \ref{S22}.
As usual, $(V,\Phi_V)$ is called irreducible if it has no proper submodules. Let $\widehat{(\fh,K)}$ be the set of equivalence classes of irreducible  $(\fh,K)$-modules.
Let $(V,\Phi_V)$ be arbitrary and fix $\pi\in \widehat{(\fh,K)}$. The $\pi$-isotypic component $V(\pi)\subset V$ is
the sum of all irreducible submodules of $V$ belonging to the equivalence class $\pi$. If $(V,\Phi_V)$ is
semisimple, then
$$ V=\bigoplus_{\pi\in \widehat{(\fh,K)}} V(\pi)\ .$$
In particular, we can consider $V$ as a $\widehat{(\fh,K)}$-graded vector space.

\begin{de}
Let $\Pi\subset\widehat{(\fh,K)}$.  An $(\fh,K)$-module $(V,\Phi_V)$ is called $\Pi$-graded,
if it is semisimple and $V(\pi)=0$ for all $\pi\not\in\Pi$.
\end{de}
Looking at the underlying $(\fh,K)$-module structure we can speak of
$\Pi$-graded $(\fh,K)$-equivariant
(metric) Lie algebras, $\Pi$-graded $(\fl,\Phi_\fl)$-modules, etc. In order to save words
we call a $\Pi$-graded $(\fh,K)$-equivariant
(metric) Lie algebra $(\fl,\Phi_\fl)$ simply a (metric) {\em Lie algebra with $\Pi$-grading}.
Thus, the term $\Pi$-grading stands for the whole equivariant structure $\Phi_\fl$, not only for the decomposition
into isotypic components. 

Let us explain the examples of $\Pi$-gradings relevant for the geometric structures discussed above.
In all these cases we can take $\fh=0$. Another interesting $\Pi$-grading for $(\fh,K)=(\RR,\ZZ_2)$ will
appear in Section \ref{S51}.
\begin{itemize}
\item $K=U(1)$, $\Pi=\{1,\sigma\}$. Here $1$ stands for the one-dimensional trivial representation, and $\sigma$ 
for the standard representation of $U(1)$ on $\CC\cong\RR^2$. A $\Pi$-grading of this kind is called {\em complex grading}.
\item $K=\RR^*$, $\Pi=\{1,\sigma,\sigma^*\}$. Here $\sigma$ stands for the standard representation of the multiplicative group $\RR^*$ on $\RR^1$, and $\sigma^*$ denotes the dual of $\sigma$. The corresponding $\Pi$-grading is called {\em para-complex grading}.
\item $K=Sp(1)$, $\Pi=\{1,\sigma\}$. Here $\sigma$ stands for the standard representation of $Sp(1)$ on $\HH\cong\RR^4$.
The corresponding $\Pi$-grading is called {\em quaternionic grading}.
\item $K=SL(2,\RR)$, $\Pi=\{1,\sigma\}$. Here $\sigma$ stands for the 2-dimensional standard representation of $SL(2,\RR)$.
The corresponding $\Pi$-grading is called {\em para-quaternionic grading}. 
\end{itemize} 

In all cases above we have a natural embedding of $\ZZ_2\cong O(1)\hookrightarrow K$ such that $\sigma(w)=-\Id$ for the non-trivial element $w\in\ZZ_2$. Therefore objects $V$ with such a $\Pi$-grading are special $\ZZ_2$-equivariant objects.
Thus they come with a splitting $V=V_+\oplus V_-$. We call such a $\Pi$-grading of a Lie algebra $\fg$
{\em proper} if $[\fg_-,\fg_-]=\fg_+$. In particular, metric Lie algebras with a proper $\Pi$-grading of this
kind are symmetric triples, which are equipped with an additional structure.

\begin{de}
A pseudo-Hermitian (para-Hermitian, hyper-K\"ahler, hypersymplectic) symmetric triple $(\fg,\Phi,\ip)$
is a metric Lie algebra $(\fg,\ip)$
with proper complex (para-complex, quaternionic, para-quaternionic) grading $\Phi$.
\end{de}

Then we have the following variant of Proposition~\ref{otto}.

\begin{pr}\label{buh}
There is a bijective map 
between isometry classes
of simply connected ``geostruc'' symmetric spaces and isomorphism classes of ``geostruc'' symmetric triples, where ``geostruc'' stands for
pseudo-Hermitian, para-Hermitian, hyper-K\"ahler, or hypersymplectic.
\end{pr}
Note that there is no analogous correspondence for (para-)quaternionic K\"ahler symmetric spaces.
Let us explain the correspondence for hyper-K\"ahler symmetric spaces. The other cases are similar.
Let $(M,g,I,J)$ be a 
hyper-K\"ahler symmetric space with base point $x_0$, and let $(\fg,\theta,\ip)$ be the associated
symmetric triple. Then $I,J,K:=IJ$ span a Lie algebra $\fk\cong \frak{sp}(1)$ 
which acts orthogonally on $\fg_-\cong T_{x_0}M$. This action commutes with the one 
of $\fg_+$. We extend the $\fk$-action to $\fg$ by the trivial action on $\fg_+$. 
For $X,Y\in\fg_-$, $Z\in\fg_+$, and $Q\in\fk$ we compute 
\begin{eqnarray*} 
\langle [QX,Y]+[X,QY], Z\rangle &=& 
\langle QX,[Y,Z]\rangle +\langle X,[QY,Z]\rangle \\ 
&=& 
\langle QX,[Y,Z]\rangle +\langle X,Q[Y,Z]\rangle=0\ . 
\end{eqnarray*} 
It follows that $\fk$ acts by derivations on $\fg$. Integrating the resulting 
homomorphism from $\frak{sp}(1)$ into antisymmetric derivations of $\fg$ we obtain the 
desired 
homomorphism $\Phi: Sp(1)\rightarrow \Aut(\fg)$. Vice versa,
 if $(\fg,\Phi,\ip)$ is a hyper-K\"ahler 
symmetric triple and $(M,g)$ is the simply connected symmetric space with symmetric 
triple 
$(\fg,\Phi(-1),\ip)$, then $I=\Phi(i)|_{\fg_-}$ and 
$J=\Phi(j)|_{\fg_-}$ 
are $\fg_+$-invariant anticommuting complex structures on $\fg_-\cong T_{x_0}M$ 
respecting the metric. They 
induce a $G$-invariant hyper-K\"ahler structure on $M$.

It is now easy to specify the classification scheme Thm.~\ref{class} for $(\fh,K)$-equivariant metric Lie algebras
to  pseudo-Hermitian (para-Hermitian, \dots) symmetric triples
(compare also Thm.~\ref{symmcl}). We don't want to write down the complete results. We only remark that any indecomposable non-semisimple ``geostruc'' symmetric triple  is isomorphic to
some $\dd$, where
\begin{itemize}
\item[(O1)] $(\fl,\Phi_\fl)$ is a Lie algebra with proper $\Pi$-grading, $\Pi$ choosen according to the geometric structure,
\item[(O2)] $\fa$ is a $\Pi$-graded semisimple orthogonal $(\fl,\Phi_\fl)$-module, and
\item[(O3)] $[\alpha,\gamma]\in \cH(\fl,\Phi_\fl,\fa)_b$ is indecomposable and satisfies Condition ($T_2$) in 
Proposition~\ref{miau}.
\end{itemize}

\subsection{Pseudo-Hermitian symmetric spaces}\label{S42}

We have seen in the previous section that in order to classify
indecomposable pseudo-Hermitian symmetric spaces\index{symmetric space!pseudo-Hermitian|(} one would have to classify the
objects (O1)--(O3) for $\Phi_\fl$ being a complex grading. This task is of a
similar complexity as the classification of all symmetric spaces discussed
in Section \ref{S3}. However, there are structural restrictions coming
from $U(1)$-equivariance and making some aspects of the theory simpler
then for general symmetric spaces. In addition, we will see that only very few symmetric spaces
of index 2 (as listed in Section \ref{S33}, Table 1) admit a K\"ahler structure.

We treat pseudo-Hermitian symmetric spaces together with 
para-Hermitian\index{symmetric space!para-Hermitian|(} ones because of the similarity
of their behaviour. Indeed, K\"ahler and para-K\"ahler structures could be viewed as different real
forms of only one complexified structure as the common complexification of $U(1)$ and $\RR^*$ is $\CC^*$.

\begin{pr}\label{gaga}
Let $(\fl,\Phi_\fl)$ be a Lie algebra with proper complex or para-complex grading.
Then the radical $\fr\subset\fl$ is nilpotent and acts trivially on every semisimple $(\fl,\Phi_\fl)$-module.
In particular, every solvable pseudo-Hermitian or para-Hermitian symmetric triple is nilpotent.
\end{pr}
\proof Differentiating $\Phi_\fl$ we can consider $(\fl,\Phi_\fl)$ as a $\fk$-equivariant Lie algebra, where $\fk$
is the Lie algebra of $U(1)$ or $\RR^*$, respectively. Properness of $(\fl,\Phi_\fl)$ implies that $\fl^\fk\subset \fl'$. Now we apply Lemma~\ref{666}.
\qed

Proposition \ref{gaga} has also the following consequence: For fixed $(\fl,\Phi_\fl)$, 
the set of isomorphism classes of (para-)complex graded semisimple orthogonal $(\fl,\Phi_\fl)$-modules is discrete. Indeed, such an orthogonal $(\fl,\Phi_\fl)$-module
is essentially determined by the action of the (semisimple) Levi factor of $\fl$ on it.

The smallest possible nonzero index of a pseudo-Hermitian symmetric triple is two.
For this case the objects (O1)--(O3) can be classified completely. Thanks to Proposition \ref{gaga}
this is considerably simpler than to classify all symmetric triples of index 2 (compare Section~\ref{S33}).

\begin{theo}[compare \cite{KO2}, Section 7.3]\label{kahl}
If $(\fg,\Phi,\ip)$ is an indecomposable pseudo-Hermitian symmetric triple of signature $(2,2q)$
that is neither semisimple
nor abelian,
then $(\fg,\Phi,\ip)$ is isomorphic to
$\dd$ for exactly one of the data in the following list:
\begin{enumerate}
\item $q=1:$ $\fl=\RR^2\cong \CC$, $\Phi_\fl=\sigma$ is the standard action of $U(1)$,
\begin{enumerate}
\item $\fa=\RR$ with the standard scalar product, $\ \rho,\Phi_\fa$ trivial,\\ 
$\alpha(z_1,z_2)=\Im \bar z_1z_2$, $z_1,z_2\in\CC$, 
$\ \gamma=0$;
\item all data as in (a) except for the opposite sign of the scalar product;
\end{enumerate}
\item $q=2:$ $\fl=\fh(1)$, $\Phi_\fl=\sigma\oplus 1$ via $\fl\cong\CC\oplus\fz(\fl)$ as a vector space,
\begin{enumerate}
\item[] $\fa=\CC$ with the real standard scalar product, $\ \rho$ trivial, $\Phi_\fa=\sigma$,\\ 
$\alpha(z,x)=z\cdot x$, 
$\alpha(z_1,z_2)=0$, $z,z_1,z_2\in\CC$, $x\in\fz(\fl)\cong\RR$, 
$\ \gamma=0$;
\end{enumerate}
\item $q=1+p,\,p\ge 0:$ $\fl=\su(2),\ \Phi_\fl=\Ad\circ i_1$, where $i_1: U(1)\rightarrow U(2)$ is the standard
embedding into the left upper corner,
\begin{enumerate}
\item[] $\fa=(\fa_1)^{p-r}\oplus(\fa_2)^r$, $\ 0\le r\le p$, where\\
$\fa_1=\CC^2$ with positive definite standard scalar product, $\rho_1$ is the standard representation of $\su(2)$, $\Phi_{\fa_1}=i_1$, and \\
$\fa_2=\su(2)$ with $-B$, $B$ the Killing form, $\rho_2=\ad$, $\Phi_{\fa_2}=\Ad\circ i_1$;\\
$\alpha=0$,
$\gamma(X_1,X_2,X_3)=cB([X_1,X_2],X_3),\ c\in\RR$.
\end{enumerate}
\item $q=1+p,\,p\ge 0:$ $\fl=\fsl(2,\RR), \ \Phi_\fl(z)=\Ad\left( i_2(\sqrt{z})\right)$, $z\in U(1)$, where 
$i_2: U(1)\cong SO(2)\rightarrow SL(2,\RR)$ is the natural embedding,
\begin{enumerate}
\item[] $\fa=(\fa_1)^{p-r}\oplus(\fa_2)^r$, $\ 0\le r\le p$, where\\
$\fa_1=\RR^2\otimes\CC$ with the scalar product given by $-\omega_{\RR^2}\otimes\omega_\CC$,  $\omega_{\RR^2},\omega_\CC$ being the standard symplectic forms of the factors, $\rho_1$ is the complexified standard representation of $\fsl(2,\RR)$, $\Phi_{\fa_1}(z)=(i_2\otimes\sigma)(\sqrt{z})$, and \\
$\fa_2=\fsl(2,\RR)$ with the Killing form $B$, $\rho_2=\ad$, $\Phi_{\fa_2}(z)=\Ad(i_2(\sqrt{z}))$;\\
$\alpha=0$,
$\gamma(X_1,X_2,X_3)=cB([X_1,X_2],X_3),\ c\in\RR$.
\end{enumerate}
\end{enumerate}
\end{theo}
The bases $N=L/L_+$ of the canonical fibrations\index{fibration!canonical} of the simply connected pseudo-Hermitian symmetric spaces corresponding to 1.-4.~are $\CC,\CC, S^2\cong\CC P^1$, and $H^2$, respectively. The complex gradings
of $\fl$ correspond to the natural complex structures of these spaces.
Note that the data defined in 4. do not depend on the choice of the square root of $z$.
The classification of complex graded semisimple orthogonal $(\fl,\Phi_\fl)$-modules in cases 3. and 4.
can be found already in \cite{CP}, Ch.~V, Prop.~3.3. 

We remark that the statement of Theorem \ref{kahl} differs slightly from the corresponding statement in \cite{KO2}.
In \cite{KO2} we determined all symmetric triples of signature $(2,2q)$ that admit the structure of a pseudo-Hermitian symmetric triple, whereas here we have determined isomorphism classes of pseudo-Hermitian symmetric triples.
Comparing both results we find that all of these symmetric triples admit exactly one complex grading (up to isomorphism). 

In a similar manner one can classify para-Hermitian symmetric triples of small index.
All para-Hermitian symmetric spaces of signature $(1,1)$ are locally isomorphic to either the flat space
$R^{1,1}$ or to the one-sheeted hyperboloid $S^{1,1}\subset\RR^{1,2}$, which is semisimple.
The complexifications of the pseudo-Hermitian symmetric triples of signature $(2,2q)$ listed in Theorem \ref{kahl} admit real forms that 
are para-Hermitian symmetric triples of signature $(q+1,q+1)$.
Concerning index 2, we obtain: 
\begin{pr}\label{sehrkahl}
There are exactly two isolated isomorphism classes and one 1-parameter family of isomorphism classes of indecomposable para-Hermitian symmetric triples of index at most $2$
that are neither semisimple
nor abelian,
namely
$\dd$ for:
\begin{enumerate}
\item $\fl=\RR^{2}\cong \RR^1\oplus(\RR^1)^*$, $\Phi_\fl=\sigma\oplus\sigma^*$,
$\ \fa=\RR$ with the standard scalar product, $\ \rho,\Phi_\fa$ trivial, 
$\ \alpha$ induced by the dual pairing, 
$\ \gamma=0$;
\item all data as in 1. except for the opposite sign of the scalar product of $\fa$;
\item $\fl=\fsl(2,\RR), \ \Phi_\fl=\Ad\circ i_3$, where 
$i_3: \RR^*\rightarrow GL(2,R)$ is the standard
embedding into the left upper corner,
$\ \fa=\{0\}$, $\ \alpha=0$,\\
$\gamma(X_1,X_2,X_3)=cB([X_1,X_2],X_3),\ c\in\RR$, $B$ being the Killing form.
\end{enumerate}
\end{pr}
Note that the family in 3.~corresponds to a 1-parameter family of para-Hermitian metrics on the symmetric space $T^*S^{1,1}$.

As in the cases of metric Lie algebras and general symmetric triples 
there seems to be no hope for a complete classification of pseudo-Hermitian and para-Hermitian triples without index restrictions.
The reason is that the Lie algebra structure of nilpotent pseudo-Hermitian and para-Hermitian symmetric triples 
can be arbitrarily complicated. 
E.g., a basic invariant
of a nilpotent Lie algebra $\fg$ is its nilindex\index{nilindex}, which is by definition the smallest non-negative integer $k$
such that $\fg^{k+1}=\{0\}$. Nilpotent Lie algebras of nilindex $k$ are sometimes also called $k$-step nilpotent
Lie algebras. The following series of examples shows that there are nilpotent pseudo-Hermitian symmetric triples with
an arbitrary large nilindex. This is in sharp contrast to the theory of
hyper-K\"ahler symmetric triples discussed in the next section, see Theorem~\ref{grob}.

\begin{ex}{\rm
For each $m\in\NN_0$ we define a pseudo-Hermitian symmetric triple $(\fg(m),\Phi,\ip)$ as follows:
As a vector space with complex grading we set 
$$\fg(m)=\CC^{m+1}\oplus\RR^m, \ \Phi=\sigma^{m+1}\oplus 1^m.$$
Let $E_i$, $i=1,\dots,m+1$, $Z_k$, $k=1,\dots,m$, be the standard basis vectors of $\CC^{m+1}$ and $\RR^m$,
respectively. Set $F_j:=iE_j$. Then $\{E_i,F_j,Z_k\}$ is a basis
of the real vector space $\fg(m)$. The nonzero Lie brackets between the basis vectors are defined as
$$ [E_i,F_j]=Z_{i+j-1},\ [Z_k,E_i]=F_{i+k},\ [Z_k,F_j]=-E_{k+j}.$$
Here $Z_l=E_q=F_q=0$ for $l>m$, $q>m+1$. Finally, the scalar product is given by $\CC^{m+1}\perp \RR^m$ and
$$ \langle E_i,F_j\rangle=0,\ \langle E_i,E_j\rangle=\langle F_i,F_j\rangle=\delta_{i+j,m+2}, \ \langle Z_k,Z_l\rangle=\delta_{k+l,m+1}\ .$$
}
\end{ex}
It is easy to check that $(\fg(m),\Phi,\ip)$ is indeed a pseudo-Hermitian symmetric triple. It is indecomposable and nilpotent of nilindex $2m+1$. Note that $\fg(1)$ and $\fg(2)$ appear in Theorem~\ref{kahl}
under 1.(a) and 2., respectively. Observe that the 
holonomy algebras 
$\fg(m)_+=\RR^m$ are abelian. However, there exist nilpotent pseudo-Hermitian 
symmetric triples having  even holonomy algebras of arbitrary large
nilindex.\index{symmetric space!pseudo-Hermitian|)}\index{symmetric 
space!para-Hermitian|)}\index{nilindex}

\subsection{Quaternionic K\"ahler  and hyper-K\"ahler symmetric 
spaces}\label{S43}
\index{symmetric space!quaternionic K\"ahler|(}\index{symmetric 
space!para-quaternionic 
K\"ahler|(}\index{symmetric space!hyper-K\"ahler|(} \index{symmetric 
space!hypersymplectic|(}
Let $(M,g,E)$ be a pseudo-Riemannian manifold of dimension $4n>4$ with quaternionic or para-quaternionic K\"ahler structure $E$. 
Then $(M,g)$ is Einstein\index{Einstein metric}. We have to distinguish between two cases: If the scalar curvature of $(M,g)$ is non-zero,
then $(M,g)$ is indecomposable and $E$ has no nontrivial parallel section. Otherwise $E$ can be spanned
by parallel sections, i.e., $(M,g)$ carries a hyper-K\"ahler
or hypersymplectic structure, respectively. For these facts we refer to \cite{AC2} and \cite{Besse}.
The latter reference deals only with the Riemannian case, but the arguments work in the indefinite case as well.
For symmetric spaces we have:
\begin{pr}[Alekseevsky/Cort\'es \cite{AC2}]\label{ac1}
Let $(M,g,E)$ be a (para-)quaternionic K\"ahler symmetric space of non-zero scalar curvature. Then the transvection
group $G$ of $(M,g)$ is simple.
\end{pr}
Indeed, the Ricci curvature of $(M,g)$ is essentially given by the Killing form of the Lie algebra $\fg$ of $G$. The Einstein property
now implies that the Killing form is non-degenerate. It follows that $G$ is semisimple. 
Simplicity of $G$ follows from the indecomposability of $(M,g)$. Complete lists of these spaces
can be found in \cite{AC2} for the quaternionic and in \cite{DJS} for the para-quaternionic case.

The above discussion and Proposition~\ref{buh} reduces the classification of (para-) quaternionic symmetric spaces to the classification
of hyper-K\"ahler (hypersymplectic) symmetric triples. First of all we have the following counterpart of
Prop.~\ref{gaga} and Prop.~\ref{ac1}.
\begin{pr}[\cite{KO3}, Prop.~2.1]\label{rhino}
Let $(\fl,\Phi_\fl)$ be a Lie algebra with proper quaternionic or para-quaternionic grading.
Then $\fl$ is nilpotent. 
In particular, every hyper-K\"ahler or hypersymplectic symmetric triple is nilpotent.
\end{pr}
\proof We can consider $(\fl,\Phi_\fl)$ as an $(\fk,\ZZ_2)$-equivariant Lie algebra, where $\fk$ is $\fsp(1)$  or $\fsl(2,\RR)$, respectively. Properness implies that this $(\fk,\ZZ_2)$-equivariant Lie algebra
satisfies the assumptions of Lemma \ref{667}, which says that $\fl$ has to be nilpotent.
\qed

That hyper-K\"ahler and hypersymplectic symmetric triples are solvable 
has been already observed in \cite{AC1}. Recently, we obtained an important sharpening of Proposition~\ref{rhino}. We shall discuss it at the end of
the present section, see Theorem~\ref{grob}. First we want to give an overview on the results on hyper-K\"ahler
symmetric triples obtained by Alekseevsky, Cort\'es, and the authors in \cite{AC1},\cite{Codd}, and \cite{KO3}.

Following our classification scheme we have to study the objects (O1)--(O3) for $\Phi_\fl$ being a quaternionic grading. Let us begin, however, with an alternative approach to hyper-K\"ahler symmetric triples due to Alekseevsky/Cort\'es \cite{AC1} that provides additional information. 
Let $(E,\omega)$ be a complex symplectic vector space. Any $S\in S^4E$
defines a complex linear subspace $\fh_S\subset\fsp(E,\omega)\cong S^2 E$ by
$$\fh_S=\Span\{S_{v,w}\in S^2E\:|\: v,w\in E\}\ ,$$
where $S_{v,w}$ is the contraction
of $S$ with $v$ and $w$ via the symplectic form $\omega$. 
If
\begin{equation}\label{cru}
S\in (S^4E)^{\fh_S}\ ,
\end{equation}
then $\fh_S\subset\fsp(E,\omega)$ is a Lie subalgebra and, moreover, there is a natural Lie bracket on 
$\fg_S:=\fh_S\oplus (\HH\otimes_\CC E)$ such that $\fh_S\subset\fg_S$ is a subalgebra and
$\fh_S$ acts on $\HH\otimes_\CC E$ by the natural action on the second factor. The remaining part of
the commutator maps $\HH\otimes_\CC E\times \HH\otimes_\CC E$ to $\fh_S$ as follows:
$$ [p\otimes v,q\otimes w]= \omega_\HH(p,q) S_{v,w}\ ,\qquad p,q\in\HH,\ v,w\in E,
$$
where $\omega_\HH$ is the alternating complex bilinear 2-form on $\HH$ such that $\omega_\HH(1,j)=1$.
Equation (\ref{cru}) is a system of quadratic equations for $S\in S^4E$. It admits families of particularly simple solutions,
namely all $S\in S^4E_+\subset S^4 E$, where $E_+\subset E$ is a Lagrangian subspace.
Let us call solutions of this kind tame. If $S$ is tame, then
the Lie algebra $\fh_S$ is abelian.

Let $J$ be a quaternionic structure on $E$ such that $J^*\omega=\bar\omega$.
Then $J$ induces a real structure on each of the spaces $S^4E$, $S^2E\cong\fsp(E,\omega)$, and 
$\HH\otimes_\CC E$. We denote all these structures by the the same symbol $\tau$.
If $S\in (S^4 E)^\tau$ satisfies (\ref{cru}), then
the real Lie algebra
$$ \fg_{J,S}:=(\fg_S)^\tau=(\fh_S)^\tau\oplus (\HH\otimes_\CC E)^\tau $$
carries a canonical structure of a hyper-K\"ahler symmetric triple. The $Sp(1)$-action on $\fg_{J,S}$
is given by
$\sigma\otimes 1$ on $(\HH\otimes_\CC E)^\tau$ and the trivial action on $(\fh_S)^\tau$.
This construction produces all 
hyper-K\"ahler symmetric triples.
\begin{pr}[Alekseevsky/Cort\'es \cite{AC1}]\label{ac2}
Let $(\fg,\Phi,\ip)$ be a hyper-K\"ahler symmetric triple, then there exist data $(E,\omega,J,S)$ as above, $S\in S^4 E$ being a $\tau$-invariant solution of (\ref{cru}), such that $(\fg,\Phi,\ip)\cong \fg_{J,S}$
as a hyper-K\"ahler symmetric triple. The tuple $(E,\omega,J,S)$ is uniquely determined by $(\fg,\Phi,\ip)$
up to complex linear isomorphisms.
\end{pr}
Thus the classification of hyper-K\"ahler symmetric triples is equivalent to the classification of 
all $\tau$-invariant
solutions of (\ref{cru}). Recall that there is a family of easy-to-find solutions called tame. Unfortunately, the claim in \cite{AC1} (repeated in \cite{AC2} and \cite{ABCV}) that all
solutions of (\ref{cru}) are tame is false. Recall that tame 
$\tau$-invariant solutions $S$ produce hyper-K\"ahler symmetric triples
with abelian holonomy\index{holonomy|(} algebra $\fg_+=(\fh_S)^\tau$.
Indeed, in the beginning of 2005 we found examples of hyper-K\"ahler symmetric triples with
non-abelian holonomy (see \cite{KO3} and Example~\ref{essig} below) and
noticed that the contradiction relies on a sign mistake in 
\cite{AC1}.  It appears that it is impossible to
find {\it all} solutions of (\ref{cru}) in a straightforward way.
In the recent paper \cite{Codd} Cort\'es reconsiders the situation and is
able to prove the following:
\begin{pr}[Cort\'es \cite{Codd}, Thm.~10]\label{ac3}
Let $S$ be a solution of (\ref{cru}). If $\fh_S$ is abelian, then $S$ is tame.
\end{pr}
Thus, Propositions \ref{ac1} and \ref{ac2} provide a classification of all hyper-K\"ahler symmetric triples
with abelian holonomy. In order to compare it with the other results of this paper we want to give the precise
formulation of this classification in terms of our standard models $\dd$.

In the following Examples \ref{oel} and \ref{essig} we consider quaternionic vector spaces $V$ in two different ways: as complex vector spaces equipped with a quaternionic structure $J$ and as $Sp(1)$-modules. Symmetric powers are symmetric powers of complex vector spaces. Then $J$ induces a real structure $\tau$ on the complex vector spaces $S^{2k}E$. The $\fl$-action on all appearing $(\fl,\Phi_\fl)$-modules $\fa$ is trivial.
\begin{ex}\label{oel}{\rm
We fix $n\in\NN$ and $S\in (S^4\HH^{n^*})^\tau$. Let $\fl=\fl_-=\HH^n$ be abelian.
The polynomial $S$ defines a symmetric bilinear form $b_S$ on the real vector space $(S^2\HH^{n})^\tau$.
We set $\fa_S=(\fa_S)_+:=(S^2\HH^{n})^\tau/\mathrm{rad}(b_S)$ and equip $\fa_S$ with 
the scalar product induced by $b_S$. We define $\alpha_S\in C^2(\HH^n,\fa_S)^{Sp(1)}$ by
$$ \alpha_S(v,w)=vJ(w)-wJ(v)\mod \mathrm{rad}(b_S)\ \in \fa_S\ ,\ \ v,w\in\HH^n\ .$$
Then $(\alpha_S,0)\in \cZ(\HH^n,\sigma^n,\fa_S)$. Moreover, $\alpha_S$ satisfies Condition ($T_2$) in 
Proposition~\ref{miau}, and $\fd_{\alpha_S,0}(\HH^n,\sigma^n,\fa_S)$ is a hyper-K\"ahler symmetric
triple. It has abelian holonomy $\fa_S$ and signature $(4n,4n)$. One can check that $\fd_{\alpha_S,0}(\HH^n,\sigma^n,\fa_S)$ is isomorphic to $\fg_{J,S}$ for $(E,\omega)=\HH^{n^*}\oplus\HH^n$, $\omega$ induced by the dual pairing.
Note that $\HH^{n^*}\subset E$ is a Lagrangian subspace, and thus $S$ is a tame solution of (\ref{cru}). 
}
\end{ex}

We call $S\in (S^4\HH^{n^*})^\tau$ indecomposable if $S\neq 0$ and $S\notin S^4 V^*\oplus S^4 W^*$ for all non-trivial decompositions
$\HH^n=V\oplus W$ into two quaternionic subspaces. We denote the set of all indecomposable $S$ by 
$(S^4\HH^{n^*})^\tau_0$. There is a natural right action of the group $GL(n,\HH)$ on $(S^4\HH^{n^*})^\tau_0$.
Moreover, if $S$ is indecomposable, then the cohomology class $[\alpha_S,0]\in \cH(\HH^n,\sigma^n,\fa_S)$
is balanced and indecomposable. Now we have the following consequence of Propositions \ref{ac1} and \ref{ac2}.

\begin{theo}[Alekseevsky, Cort\'es]\label{ac4}
The assignment
$$(S^4\HH^{n^*})^\tau_0\ni S\mapsto \fd_{\alpha_S,0}(\HH^n,\sigma^n,\fa_S)$$ 
yields
a bijection between the union of the orbit spaces $(S^4\HH^{n^*})^\tau_0/GL(n,\HH)$, $n\in\NN$,
and the set of isomorphism classes of indecomposable non-abelian hyper-K\"ahler symmetric spaces
with abelian holonomy.
\end{theo}

Note that this rather satisfactory classification result is not a classification in the sense of a list
since for large $n$ the orbit spaces $(S^4\HH^{n^*})^\tau_0/GL(n,\HH)$ are not explicitly known. 

\begin{ex}\label{essig}{\rm
Fix $n\in\NN$, $p\in\NN_0$. We define a Lie algebra $(\fl^n,\Phi_n)$ with proper quaternionic grading
as follows: $\fl^n_-:=\HH^n$, $\fl^n_+:=(S^2\HH^n)^\tau$, $\fl^n_+=\fz(\fl^n)$,  
$$[v,w]_\fl:=vJ(w)-wJ(v)\in\fl^n_+,\ \ v,w\in\HH^n.$$ 
We set $\fa_n=(\fa_n)_-:=S^3\HH^n$. The quaternionic structure of $\HH^n$ induces one
on $S^3\HH^n$. Thus $\fa_n$ is a quaternionic vector space. The standard $Sp(1)$-invariant complex
Hermitian form on $\HH^n$ of signature $(2p,2(n-p))$ induces a Hermitian form on $\fa_n$. 
We equip $\fa_n$ with the real
part of this Hermitian form and denote the resulting orthogonal $(\fl^n,\Phi_n)$-module by $\fa_{n,p}$.
We define $\alpha_n\in C^2(\fl^n,\fa_n)^{Sp(1)}$ by
$$\alpha_n(v,L)=vL \in S^3\HH^n, \ \ v\in\HH^n, L\in (S^2\HH^n)^\tau,\quad \alpha_n(\fl^n_-,\fl^n_-)=\alpha_n(\fl^n_+,\fl^n_+)=0. $$ 
The above mentioned Hermitian form on $\HH^n$ induces a natural identification of $\fl^n_+=(S^2\HH^n)^\tau$
with the Lie algebra $\fsp(p,n-p)$ and a scalar product $\ip_p$ on $\fl^n_+$. We denote the resulting Lie bracket on
$\fl^n_+$ by $\lb_p$. Eventually, we define a 3-form $\gamma_p\in C^3(\fl^n_+)\subset C^3(\fl^n)^{Sp(1)}$ by
$$ \gamma_p(L_1,L_2,L_3):=\langle [L_1,L_2]_p,L_3\rangle_p,\ \ \ L_i\in\fl_+^n\ . $$
Then $(\alpha_n,\gamma_p)\in\cZ(\fl^n,\Phi_n,\fa_{n,p})$. Moreover, the associated cohomology class is balanced,
indecomposable, and satisfies ($T_2$). Thus $\fd_{\alpha_n,\gamma_p}(\fl^n,\Phi_n,\fa_{n,p})$ is an indecomposable
hyper-K\"ahler symmetric triple. Because of the form of $\gamma_p$ its holonomy algebra $(\fl^n_+)^*\oplus \fl_+^n$
is non-abelian.
}
\end{ex}
These hyper-K\"ahler symmetric triples are natural generalisations of  Example~1 in \cite{KO3}, which is
isomorphic to $\fd_{\alpha_1,\gamma_0}(\fl^1,\Phi_1,\fa_{1,0})$ and has signature $(4,12)$. The paper \cite{KO3}
contains further examples with non-ablian holonomy. 

The smallest possible index of a non-abelian 
hyper-K\"ahler symmetric triple is~$4$. In the following classification result we use the notation of Examples~\ref{oel}, \ref{essig}.

\begin{theo}[\cite{KO3}, Thm.~7.4]
Let $(\fg,\Phi,\ip)$ be a non-abelian indecomposble hyper-K\"ahler symmetric triple of signature $(4,4q)$.
Then $q=1$ or $q=3$, and $(\fg,\Phi,\ip)$ is isomorphic to exactly only one of the following triples:
\begin{enumerate}
\item[q=1:]$\ \fd_{\alpha_{S_\lambda},0}(\HH,\sigma,\fa_{S_\lambda})$, 
$\lambda\in [-2, 2]$, $\ S_\lambda(z+wj):=z^4+\lambda z^2w^2+w^4$, $z,w\in\CC$;
\item[q=3:] $\ \fd_{\alpha_1,\gamma_0}(\fl^1,\Phi_1,\fa_{1,0})$.
\end{enumerate}
\end{theo}
The proof of the theorem in \cite{KO3} is based on Thm.~\ref{class}, i.e., it classifies directly the
relevant objects (O1)--(O3). It does not rely on Theorem~\ref{ac4}.

There is a completely parallel theory for hypersymplectic symmetric triples. E.g., in the Alekseevsky-Cort\'es
construction one has simply to replace the real structure $\tau$ induced by $J$ by a real structure coming from real structures on both factors $\HH$ and $E$. In other words, one can work from the beginning with $\RR^2\otimes_\RR E_0$
instead of $\HH\otimes_\CC E$,
where $E_0$ is a real symplectic vector space. Hypersymplectic triples with abelian holonomy can be classified
as in Theorem \ref{ac4}, see \cite{Codd,DJS,ABCV}. In addition, there is a variant of Example~\ref{essig} producing 
hypersymplectic symmetric triples with non-abelian 
holonomy (one has to replace $\HH^n$ by $\RR^2\otimes \RR^{2n}$,
the parameter $p$ disappears).

Note that all examples of hyper-K\"ahler symmetric triples presented so far have nilindex\index{nilindex} at most 5.
As the following theorem and its corollary show, this is not an accident.

\begin{theo}[\cite{O}]\label{grob}
Let $(\fl,\Phi_\fl)$ be a Lie algebra with proper quaternionic or para-quaternionic grading.
Then the nilindex of $\fl$ is at most 6. \index{nilindex}
\end{theo}

\begin{co}\label{borg}
Let $(\fg,\Phi,\ip)$ be a non-abelian hyper-K\"ahler or hypersymplectic symmetric triple.
As usual, let $\fl=\fg/\fri(\fg)^\perp$, where $\fri(\fg)\subset\fg$ is the canonical isotropic ideal.
Then there are only two possiblities for the nilindices\index{nilindex} of $\fg$, $\fg_+$, and $\fl$
as listed in the following table:
\renewcommand{\arraystretch}{1.5}
$$ \mbox{\begin{tabular}{|c|c|c|c|c|}
\hline
&{\rm holonomy}&$\fg$&$\fg_+$&$\fl$\\
\hline\hline
1.&{\rm abelian}&$3$&$1$&$1$\\
\hline
2.&{\rm non-abelian}&$5$&$2$&$2$\\
\hline
\end{tabular}}\ .
$$\index{holonomy|)}
\end{co} 
\proof Let us denote the nilindex of a nilpotent Lie algebra $\fh$ by $n(\fh)$.
It is easy to see that $n(\fg)$ is odd and that $2n(\fg_+)<n(\fg)$ for any nilpotent symmetric pair $\fg=\fg_+\oplus\fg_-$. Thus in our case $n(\fg)=3$ or $n(\fg)=5$ by Thm.~\ref{grob}.
The possibility that at the same time $n(\fg)=5$ and $n(\fg_+)=1$ is excluded by Thm.~\ref{ac4}.
The corollary now follows from the inequality
$ 2n(\fl)\le n(\fg)\le 2n(\fl)+1$,
which holds for any nilpotent metric Lie algebra and can be derived from the definition of 
$\fri(\fg)$.
\qed

\begin{co}
The base $N$ of the canonical fibration\index{fibration!canonical} (see Section~\ref{S32})
of a hyper-K\"ahler or hypersymplectic symmetric space $M$ is flat.
\end{co}
\proof By Corollary~\ref{borg} the Lie algebra $\fl$ is at most two-step nilpotent. It follows that
$N=L/L_+$ is flat.
\qed

Corollary \ref{borg} shows that the structure of hyper-K\"ahler and hypersymplectic symmetric triples
is strongly restricted. Therefore the goal of a full classification of these triples might be not unrealistic.
The (lengthy) proof of Theorem \ref{grob}, which we don't
want to explain here, gives further inside into the structure of these triples: it computes
certain universal Lie algebras $\frak L_n$ with proper (para-) quaternionic grading such that any Lie algebra with proper (para-)quaternionic grading in $n$ generators is a quotient of $\frak L_n$.

We conclude this section with the following open question. Is the dimension of every hyper-K\"ahler
symmetric space without flat local factors divisible by $8$? Of course, this is true for all known
examples. Note that indecomposable hypersymplectic symmetric spaces exist in all dimensions that are multiples of $4$.
\index{symmetric space!quaternionic K\"ahler|)}\index{symmetric 
space!para-quaternionic 
K\"ahler|)}\index{symmetric space!hyper-K\"ahler|)} \index{symmetric 
space!hypersymplectic|)}

\section{Further applications} 
\subsection{Extrinsic symmetric spaces}\label{S51}\index{symmetric 
space!extrinsic|(}
Let us consider a non-degenerate connected submanifold $M\subset {\RR}^{p,q}$. For $x\in M$ 
let $s_x$ be the reflection of $\RR^{p,q}$ at the normal space $T^\perp_x M$ of $M$ 
at $x$, i.e., $s_x$ is an affine isometry and $s_x|_{T_x M}=-\Id$,   
$s_x|_{T^\perp_x M}=\Id$. Here and throughout the section we consider $T_{x}M$ and $T^\perp_x M$ as affine subspaces of $\RR^{p,q}$. Then $M$ is called extrinsic symmetric if $s_x(M)=M$ 
holds for each point $x\in M$. Extrinsic symmetric 
spaces in ${\RR}^{p,q}$ are exactly those complete submanifolds whose second 
fundamental form is parallel. Extrinsic symmetric spaces in the Euclidean space are 
well understood. A classification in this case follows from Ferus' results 
discussed below and the classification of symmetric $R$-spaces due to Kobayashi and 
Nagano. The case of a pseudo-Euclidean 
ambient space seems to be more involved. 

In Section \ref{S31} we discussed the correspondence between 
pseudo-Riemannian symmetric spaces and symmetric triples. Here we will 
see that there is a similar correspondence for (a certain class of) 
extrinsic symmetric spaces. While symmetric triples are special 
(namely proper) $\ZZ_{2}$-equivariant metric Lie algebras the 
algebraic objects that we will use here are certain 
$(\RR,\ZZ_{2})$-equivariant metric Lie algebras. The group 
$\ZZ_{2}=\{1,-1\}$ acts on $\RR$ by multiplication. Hence for an 
$(\RR,\ZZ_{2})$-equivariant metric Lie algebra $(\fg,\Phi,\ip)$ we can 
regard $\Phi$ as a pair $(D,\theta)$ that consists of a derivation $D\in\Der(\fg)$ 
and an involution $\theta\in\Aut(\fg)$ such that $D\theta=-\theta D$. 
Assume that for such an $(\RR,\ZZ_{2})$-equivariant metric Lie algebra   
$D^{3}=-D$ holds.  Then the eigenvalues of $D$ are in $\{i,-i,0\}$. We 
put 
$$\fg^+:=\Ker D,\ \fg^-:=\Span\{X\in\fg\mid D^2(X)=-X\}$$ 
and define an involution $\tau_{D}$ on $\fg$ by 
$\tau_D:\fg\rightarrow\fg,\  \tau_D|_{\fg^+}=\Id,\ 
\tau_D|_{\fg^-}=-\Id.$ 
Obviously $\tau_D$ and 
$\theta$ commute, hence $\fg_+$ and $\fg_-$ are invariant under 
$\tau_D$. 
We introduce the notation 
\begin{equation}\label{not} 
\fg_+^+:=\fg_+\cap\fg^+,\ \fg_+^-:=\fg_+\cap\fg^-,\ 
\fg_-^+:=\fg_-\cap\fg^+,\ \fg_-^-:=\fg_-\cap\fg^-. 
\end{equation} 
\begin{de}\label{D51} 
    An extrinsic symmetric triple\index{symmetric triple!extrinsic} is an $(\RR,\ZZ_{2})$-equivariant metric 
    Lie algebra $(\fg,\Phi,\ip)$, $\Phi=(D,\theta)$, for which 
    \begin{itemize} 
        \item[(i)] the derivation $D$ is inner and satisfies $D^{3}=-D$, 
        \item[(ii)] the $\ZZ_{2}$-equivariant metric Lie algebras 
        $(\fg,\theta,\ip)$ and $(\fg_+,\tau_D|_{\fg_+},\ip|_{\fg_+})$ 
        are proper (i.e.~symmetric triples). 
    \end{itemize}         
\end{de} 
Two extrinsic symmetric triples are called isomorphic if they are 
isomorphic as $(\RR,\ZZ_{2})$-equivariant metric Lie algebras. 
\begin{re}\label{RO2}{\rm 
We consider  the subset 
$$\Pi=\{(\RR,\Phi_{1}),\ (\RR,\Phi_{2}),\ (\RR^{2},\Phi_{3})\}\subset \widehat{(\RR,\ZZ_{2})},$$ 
where $\Phi_1(r)=\Phi_2(r)=0, \Phi_{1}(z)=1$, $\Phi_{2}(z)=z$,  and 
$$ \Phi_{3}(r)=\small{ 
\left( \begin{array}{cc} 
0&-r\\ 
r&0 
\end{array}\right) 
},\ \Phi_{3}(z)=  \small{ 
\left( \begin{array}{cc} 
1&0\\ 
0&z 
\end{array}\right) 
}  $$
for $r\in\RR$ and $z\in\ZZ_{2}$.
In other words, $\Pi$ contains exactly those representations of 
$(\RR,\ZZ_2)$ that integrate to one of the following representation of $O(2)$:
the one-dimensional trivial representation, the one-dimensional representation via the 
determinant or the two-dimensional standard representation. 
Then an $(\RR,\ZZ_{2})$-equivariant (metric) Lie algebra $(\fg,\Phi,\ip)$, 
$\Phi=(D,\theta)$, satisfying $D^3=-D$ is the same as a $\Pi$-graded (metric) Lie 
algebra. 
}\end{re} 
Let $(\fg,\Phi,\ip)$, $\Phi=(D,\theta)$, be an extrinsic symmetric triple and choose 
$\xi\in\fg_{-}$ such that 
$D=\ad(\xi)$. We consider the subgroup 
$G_{+}:=\langle\, \exp(\ad X)|_{\fg_-}\mid X\in\fg_+\,\rangle$ of 
$O(\fg_{-})$ and define 
$$M_{\fg,\xi}:=G_+\cdot\xi\subset\fg_-.$$ 
\begin{pr}[\cite{Kext}]\label{POrbit} 
    \begin{enumerate} 
        \item        For any extrinsic symmetric triple 
        $(\fg,\Phi,\ip)$ and for each  $\xi\in\fg_{-}$ with $D=\ad(\xi)$ 
        the submanifold $M_{\fg,\xi}\subset\fg_-$ is an extrinsic symmetric 
        space. The abstract symmetric space $M_{\fg,\xi}$ with the induced metric 
is associated with the symmetric triple 
        $(\fg_+,\tau_D|_{\fg_+},\ip|_{\fg_+})$. 
        \item  Let $(\fg_{i},\Phi_{i},\ip_{i})$, $i=1,2$, be extrinsic 
        symmetric triples and let $M_{\fg_{i},\xi_{i}}$ be constructed 
        as above. Then there exists an affine isometry 
        $f:(\fg_{1})_{-}\rightarrow(\fg_{2})_{-}$ mapping $M_{\fg_{1},\xi_{1}}$ to 
$M_{\fg_{2},\xi_{2}}$ if and only if 
        $(\fg_{1},\Phi_{1},\ip_{1})$ and $(\fg_{2},\Phi_{2},\ip_{2})$ are 
        isomorphic. 
    \end{enumerate}         
\end{pr} 

Now let us turn to a construction that may be considered as a converse 
of Prop.~\ref{POrbit}, 1. Given an extrinsic symmetric space 
$M\subset \RR^{p,q}$ satisfying certain additional assumptions 
this construction yields an extrinsic symmetric triple $(\fg,\Phi,\ip)$ 
and an element $\xi\in\fg_{-}$ such that $M=M_{\fg,\xi}$. 
Let us first discuss these additional assumptions. 
\begin{de} 
A submanifold $M\subset \RR^{p,q}$ is called full if it is not contained 
in any proper affine subspace of $\RR^{p,q}$. It is called normal if the 
intersection of the 
normal spaces of all points 
of $M$ is not empty, i.e. $\bigcap_{x\in M} T^\perp_x M\not=\emptyset$. 
\end{de} 

To be full is not really a restriction for submanifolds $M$ of the Euclidean 
space $\RR^n$ since we may always consider the smallest affine 
subspace that contains $M$. Contrary to that, there are many 
submanifolds of the pseudo-Euclidean space $\RR^{p,q}$ that are contained in 
an affine subspace that is degenerate with respect to the inner product 
but not in a proper non-degenerate one. This is also the case if one restricts 
oneself to extrinsic symmetric spaces. As far as normality is 
concerned, Ferus proved in \cite{F1} that extrinsic symmetric spaces in the 
Euclidean 
space decompose into a product of an affine subspace and a normal extrinsic 
symmetric space. This seems to be not true for a pseudo-Euclidean 
ambient space. 

Let us now describe the construction. It was developed by Ferus 
\cite{F2,F3} who used it to prove that every full and normal 
extrinsic symmetric space in $\RR^n$ is a standard imbedded symmetric 
R-space\index{symmetric R-space}. 
In our language this means that any such extrinsic symmetric 
space arises from an extrinsic symmetric triple as described in Proposition 
\ref{POrbit},~1. Here we will present the more elementary 
description of this construction given by Eschenburg and Heintze 
\cite{EH}. We will also 
include the necessary modifications for the pseudo-Riemannian situation 
discussed in \cite{Kext, Kim}. 

Let $M\subset {\RR}^{p,q}$ be a full and normal extrinsic symmetric space. 
Fix a point $x_0\in M$. Since $M$ is normal we may assume 
that $0\in\RR^{p,q}$ is contained in the intersection of normal spaces 
$\bigcap_{x\in M} T^\perp_x M$.  This implies that  $s_x\in O(p,q)$ for all 
$x\in M$. 
We define the transvection group\index{transvection group} of $M\subset {\RR}^{p,q}$ by 
$$K:=\langle s_{x}\circ s_y \mid x,y\in 
M\rangle\subset O(p,q).$$ Obviously, $K$ is isomorphic to the transvection group 
of the abstract symmetric space $M$. Let $\fk\subset\so(p,q)$ be the Lie algebra 
of $K$.  Since it is 
isomorphic to the Lie algebra of the transvection group of the 
abstract  symmetric space $M$ it is the underlying Lie algebra of a symmetric 
triple $(\fk,\theta_{\fk},\ip_{\fk})$. 
Besides $\fk$ we can associate with $M\subset {\RR}^{p,q}$ the 
following 
metric Lie algebra $\fg$. As a vector space $\fg$ equals 
$$\fg=\fk\oplus V,\quad V=\RR^{p,q}.$$ 
Using the standard scalar product $\ip_{p,q}$ on $V$ we define 
a scalar product on $\fg$ by $\ip=\ip_\fk\oplus \ip_{p,q}$. 
Furthermore, we define a bracket 
$\lb:\fg\times\fg\rightarrow\fg$  by 
\begin{itemize} 
\item[(i)] $\lb$ restricted to ${\fk\times\fk}$ equals the Lie bracket of 
$\fk$, 
\item[(ii)] $[A,v]=A(v)$ for $A\in\fk$ und $v\in 
V$, 
\item [(iii)] $\lb$ restricted to $V\times V$ is given by the condition 
$\lb_{V\times 
V}: 
V\times V\rightarrow \fk$ and by 
$\langle A,[x,y] \rangle=\langle [A,x],y\rangle$ for all 
$A\in\fk$ and $x,y\in V$. 
\end{itemize} 
Using the fullness of $M\subset {\RR}^{p,q}$ one can prove that this 
bracket satisfies the Jacobi identity. Moreover, we have an involution $\theta$ on 
$\fg$ given by $\fg_+=\fk$, 
$\fg_-=V$ and an inner derivation $D:=\ad(x_0)$. 
We obtain: 
\begin{theo}[Ferus \cite{F2,F3}, see also \cite{EH,Kext,Kim}] 
  Let $M\subset \RR^{p,q}$ be a full and normal extrinsic 
symmetric space. Let $\fg,\ip$ and $\Phi=(D,\theta)$ be as constructed 
above. Then $(\fg,\Phi,\ip)$ is an extrinsic symmetric triple and 
$M=M_{\fg,x_{0}}$.   
\end{theo} 
In \cite{Kim} Kim discusses this construction also for the case of a non-normal 
extrinsic symmetric space. 

\begin{re}{\rm
    Let $(\fg,\Phi,\ip)$, $\Phi=(D,\theta)$, be an extrinsic symmetric 
    triple and suppose $D=\ad(\xi)$, $\xi\in\fg_{-}$. Then 
    $M_{\fg,\xi}$ is normal by construction. It is full if and only 
    if 
    \begin{equation} \label{irgend} [\fg_{+}^{-},\fg_{-}^{-}]=\fg_{-}^{+}
    \end{equation}
    (see (\ref{not}) for the notation used here). We call an 
    extrinsic symmetric triple satisfying (\ref{irgend}) full.
}\end{re}

Having explained the relation between extrinsic symmetric spaces and extrinsic 
symmetric triples we now want to apply our classification scheme to extrinsic 
symmetric triples. Recall that we can understand an extrinsic symmetric 
triple as a  $\Pi$-graded metric Lie algebra satisfying additional 
conditions, where $\Pi$ is given as in 
Remark~\ref{RO2}. We can proceed as in Section~\ref{S41} and we see that any 
extrinsic symmetric triple is isomorphic to some $\dd$, where now 
\begin{itemize} 
\item[(O1)] $(\fl,\Phi_\fl)$ is a Lie algebra with $\Pi$-grading such that the 
$\ZZ_2$-equivariant Lie algebras $(\fl,\theta_\fl)$ and 
$(\fl_+,\tau_{D_\fl}|_{\fl_+})$ are proper, 
\item[(O2)] $\fa$ is a $\Pi$-graded semisimple orthogonal $(\fl,\Phi_\fl)$-module, 
\item[(O3)] $[\alpha,\gamma]\in \HPhi_b$ is indecomposable and satisfies besides   
$(T_2)$ two further conditions $(T_2^+)$ and $(A_0^+)$ (see \cite{Kext}), which 
ensure that $\dd$ satisfies the properness condition (ii) in Def.~\ref{D51}, 
\item[(O4)] the derivation $D=-D_\fl^*\oplus D_\fa \oplus D_\fl$ is inner. 
\end{itemize} 
Condition (O4) is equivalent to 
$\exists l\in\fl_-\ \exists a\in \fa_-\subset 
C^0(\fl,\fa)\ \exists 
z\in\fl_-^{*}\subset C^1(\fl)\ :$
$$D_\fl=\ad_\fl(l),\ D_\fa=\rho(l),\ da=i(l)\alpha,\ 
dz=\langle a\wedge \alpha\rangle +i(l)\gamma.$$

Combining (O1) -- (O4) with Theorem \ref{class} we obtain a classification scheme 
for extrinsic symmetric triples. 

\begin{ex}[Full and normal extrinsic Cahen-Wallach spaces]{\rm  We want to 
    answer the question which indecomposable non-semisimple non-flat Lorentzian symmetric 
    spaces\index{symmetric space!Lorentzian} can be embedded into a pseudo-Euclidean space as full and normal 
    extrinsic symmetric spaces. Recall from Theorem~\ref{TCW} that 
    every indecomposable 
    non-semisimple non-flat Lorentzian symmetric space is associated 
    with a symmetric triple of the kind $\fd(p,q,\lambda,\mu)$, 
    $p,q\ge0$, $p+q>0$, $(\lambda,\mu)\in M_{p,q}$. 
   \begin{theo}[cf.~\cite{Kext}] \label{Text} The symmetric triple  $\fd(p,q,\lambda,\mu)$, 
    $p,q\ge0$, $p+q>0$, $(\lambda,\mu)\in M_{p,q}$ admits an associated 
    symmetric space $M^{1,n+1}$, $n=p+q$, that can be embedded as a 
    full and normal extrinsic symmetric space if and only if 
    \begin{enumerate}
\item	$q=0$ and $\lambda=(1,\ldots,1)$ or 
 \item $p=0$ and $\mu=   (1,\ldots,1)$.
 \end{enumerate}
     In both cases the ambient space is $\RR^{2,n+2}$ and there is a 
     one-parameter family of non-isomorphic embeddings.
 \end{theo}
    To verify the theorem we have to determine all  extrinsic symmetric 
    triples $(\fg,\Phi,\ip)$, $\Phi=(D,\theta)$, that satisfy
    (\ref{irgend}) and 
    $(\fg_{+},\tau_{D}|_{\fg_{+}},\ip_{\fg_{+}})\cong 
    \fd(p,q,\lambda,\mu)$ with $(\lambda,\mu)\in 
    M_{p,q}$ (see Section \ref{S33}). In particular, since
    $\fg_{+}$ is indecomposable and non-reductive the Lie algebra $\fg$ is 
    indecomposable and non-semisimple. This implies $\fg\cong\fd:=\dd$ for 
    a suitable $(\RR,\ZZ_{2})$-equivariant Lie algebra 
    $(\fl,\Phi_\fl)$, an orthogonal $(\fl,\Phi_\fl)$-module $\fa$ 
    and some $(\alpha,\gamma)\in \cZ(\fl,\Phi_{\fl},\fa)$. In 
    particular, we have $\fg_{+}\cong\fd_{+}=\fd_{\alpha|_{\fl_{+}}, 
    \gamma|_{\fl_{+}}}(\fl_{+},\Phi_{\fl}|_{\fl_{+}},\fa_{+})$. 
    Hence $\dim \fl^{-}_{+}=\dim \fl^{-}_{-}=1$ and $\fl^{+}_{+}=[\fl^{-}_{+}, 
    \fl^{-}_{+}]=0$.  Equation (\ref{irgend})  
    implies $[\fl^{-}_{+},\fl^-_{-}]=\fl_-^{+}$. Since 
    $\fl^{+}_{-}\not=0$ we also obtain  $\dim\fl^{+}_{-}=1$. 
    Now it is not hard to see that $\fl$ is isomorphic to one 
    of the Lie algebras 
    \begin{enumerate}
    \item $\fsl(2,\RR)=\{[H,X]=2Y,\,[H,Y]=2X,\,[X,Y]=2H\}$ 
    \item$\su(2)=\{[H,X]=2Y,\,[H,Y]=-2X,\,[X,Y]=2H\}$ 
    \end{enumerate}
    and that in both 
    cases  $\Phi_{\fl}=(D_{\fl},\theta_{\fl})$ is 
    given by $\fl_{+}=\RR\cdot H$, $\fl_{-}=\Span\{X,Y\}$ and 
    $D_{\fl}=(1/2)\cdot \ad X$.  Since 
    $H^{2}(\fl,\fa)=0$ for $\fl\in\{\fsl(2,\RR),\su(2)\}$ we may 
    assume $\alpha=0$. Since $\fd$ is indecomposable we get 
    $\fa^{\fl}=0$. Because of $D^{3}=-D$ the eigenvalues of $\rho(X)$ 
    are in $\{0,2i,-2i\}$. Consequently, $\fa$ is the direct sum of 
    submodules that are all equivalent to the adjoint representation 
    of $\fl$. Finally, we obtain the following result, which proves 
    Theorem~\ref{Text} and gives all embeddings explicitly. 
    
    \begin{pr}[cf.~\cite{Kext}] \label{PCW} If $M^{1,n+1}\hookrightarrow\RR^{r,s}$ is a full and normal extrinsic 
    symmetric space of Lorentz signature and if $M^{1,n+1}$ is solvable 
    and indecomposable, 
    then $(r,s)=(2,n+2)$ and we are in one of the following cases:
    \begin{enumerate}
    \item  $M^{1,n+1}$ is associated to the 
    symmetric triple $\fd(n,0,\lambda,0)$,  
    $\lambda=(1,\ldots,1)\in\RR^{n}$ and $M^{1,n+1}\subset\RR^{2,n+2}$ is extrinsic 
    isometric to $M_{\fd,\xi}$ for
$\fd= \dd$ with
    \begin{eqnarray*} &&\fl=\fsl(2,\RR),\ \Phi_{\fl} \mbox{ as 
    above }, 	\fa=\textstyle{\bigoplus_{i=1}^{n}}(\ad,\fl,B_\fl,\Phi_{\fa}),\  
     \Phi_{\fa}=(D_{\fl},-\theta_{\fl}),\\
     &&\alpha=0,\ \gamma=cB_{\fl}(\cdot,\lb),\ c\in\RR, 
     \end{eqnarray*}
     and $\xi=(1/2)\cdot X$. Here $B_{\fl}$ is the Killing form of $\fl$.
     \item $M^{1,n+1}$ is associated to the 
    symmetric triple $\fd(0,n,0,\mu)$,  
    $\mu=(1,\ldots,1)\in\RR^{n}$ and $M^{1,n+1}\subset\RR^{2,n+2}$ is extrinsic 
    isometric to $M_{\fd,\xi}$ for
$\fd= \dd$ with
    \begin{eqnarray*} &&\fl=\su(2),\ \Phi_{\fl} \mbox{ as 
    above }, 	\fa=\textstyle{\bigoplus_{i=1}^{n}}(\ad,\fl,-B_\fl,\Phi_{\fa}),\  
     \Phi_{\fa}=(D_{\fl},-\theta_{\fl}),\\
     &&\alpha=0,\ \gamma=cB_{\fl}(\cdot,\lb),\ c\in\RR, 
     \end{eqnarray*}
     and $\xi=(1/2)\cdot X$. Again $B_{\fl}$ denotes the Killing form of $\fl$.
     \end{enumerate}
     \end{pr}

     We remark that there are decomposable solvable non-flat Lorentzian symmetric 
    spaces that can be embedded as full and normal extrinsic symmetric 
    spaces such that the embedding is indecomposable, i.e., the 
    associated extrinsic symmetric triple is indecomposable.
     
    In \cite{Kext} we classify all extrinsic symmetric triples without simple 
    ideals that are associated with a Lorentzian extrinsic symmetric 
    space. Of course, that classification contains all 
    extrinsic symmetric triples that we obtained above. However, it 
    contains also extrinsic symmetric triples whose associated 
    extrinsic symmetric spaces are not full. 
}\end{ex} 
\begin{re}{\rm
 Let $(\fg,\Phi,\ip)$, $\Phi=(D,\theta)$, be an extrinsic symmetric 
 triple. Exponentiating $D$ defines a complex grading $\Phi_{0}$ 
 on $(\fg,\ip)$. Equation (\ref{irgend}) implies that this grading is 
 proper if $(\fg,\Phi,\ip)$ is full. Thus, forgetting about $\theta$, 
 a full extrinsic symmetric triple can be considered as a 
 pseudo-Hermitian symmetric triple such that the complex grading is 
 given by inner automorphisms. If 
 $(\fg_{+},\tau_{D}|_{\fg_{+}},\ip_{\fg_{+}})$ has signature $(p,q)$, 
 then the pseudo-Hermitian symmetric triple $(\fg,\Phi_{0},\ip)$ has 
 signature $(2p,2q)$. Moreover, Prop.~\ref{gaga} implies that $\fg$ has a 
 non-trivial Levi factor and that the radical of $\fg$ is nilpotent.
 In view of these facts, the reader should compare Prop.~\ref{PCW} 
 with Thm.~\ref{kahl}.
}\end{re}\index{symmetric 
space!extrinsic|)}

\subsection{Manin triples}\label{S52} 
Manin triples\index{Manin triple} are algebraic objects that are associated 
with Poisson-Lie groups\index{Poisson-Lie group}. 
A Poisson-Lie group is a Lie group equipped with a Poisson bracket 
that satisfies a compatibility condition with the group multiplication. 
The infinitesimal object associated with such a Poisson-Lie group $G$ is 
the Lie algebra $\fg$ of $G$ together with  
a 1-cocycle $\gamma:\fg\rightarrow 
\fg\wedge\fg$ that satisfies co-Jacobi identity, i.e., $\gamma^{*}: \fg^{*}\wedge\fg^{*}\rightarrow 
\fg^{*}$ is a Lie bracket on $\fg^{*}$. Such a pair $(\fg,\gamma)$ is 
called Lie bialgebra. Given such a Lie  bialgebra $(\fg,\gamma)$ it is 
easy to see that there exists a unique Lie algebra structure on the 
vector space $\fg\oplus\fg^{*}$ such that the inner product on  
$\fg\oplus\fg^{*}$ defined by the dual pairing is invariant. What we 
obtain is an example of a so-called Manin triple.
\begin{de}
    A Manin triple $(\fg,\fh_{1},\fh_{2})$ consists of a metric Lie 
    algebra $\fg$ and two complementary isotropic subalgebras 
    $\fh_{1}$ and $\fh_{2}$.
\end{de}    
Here we only consider finite-dimensional Lie bialgebras\index{Lie 
bialgebra} and Manin triples. In this case 
the above described assignment associating a Manin triple to a given 
Lie bialgebra is a one-to-one correspondence. For a detailed 
introduction to this subject we refer to \cite{KS}.

Manin triples $(\fg,\fh_{1},\fh_{2})$ for which $\fg$ is a complex 
reductive Lie algebra were classified by Delorme \cite{De}. Having 
learned about the difficulties in classifying non-reductive metric Lie 
algebras it is not surprising that there is little known about  
non-reductive Manin triples. However, there are some results in 
low dimensions. Figueroa-O'Farrill \cite{OF} studied  Manin triples in the context of 
conformal field theory and the $N=2$ Sugawara constructions. By 
rather heavy calculations he achieved a classification of 
complex six-dimensional Manin triples. A more conceptual poof of his 
result, which also yields a classification in the real case is due to 
Gomez \cite{G}.

Unfortunately, the method of quadratic extensions seems to be not 
adequate for the description of Manin triples. The difficulty is 
to handle both isotropic subalgebras at the same time. However, we 
can describe Manin pairs.
\begin{de}
    A Manin pair\index{Manin pair} $(\fg,\fh)$ consists of a metric Lie algebra of 
signature $(n,n)$ and  an $n$-dimensional isotropic subalgebra $\fh\subset \fg$.
\end{de}
Applying the method of quadratic extensions to Manin pairs we obtain the following result.
\begin{pr}\label{PManin}
Let  $\fl$ be a Lie algebra and  $\fl'\subset \fl$ a subalgebra. Let $\fa$ be an orthogonal $\fl$-module of signature $(m,m)$ and $\fa'\subset\fa$ an $m$-dimensional isotropic $\fl'$-invariant subspace. Let $(\alpha,\gamma)\in\cZ(\fl,\fa)$ satisfy
$\alpha(\fl',\fl')\subset\fa'$ and 
$\gamma(\fl',\fl',\fl')=0$. Then $\fd:=\fd_{\alpha,\gamma}(\fl,\fa)$ 
has  signature $(m+\dim \fl,m+\dim \fl)$ and 
$(\fd,\Ann(\fl')\oplus\fa'\oplus\fl')$ is a Manin  pair, where $\Ann(\fl')\subset \fl^*$ denotes the annihilator of $\fl'$.

Conversely, every Manin pair $(\fg,\fh)$ for which 
$\fg$ does not contain simple ideals is isomorphic to some 
pair $(\fd,\Ann(\fl')\oplus\fa'\oplus\fl')$ constructed in this way, 
where,  moreover, $\fd$ is balanced.
\end{pr} 
\proof
The first part of the proposition  is easy to check. The second part relies 
on the functorial   assignment (\ref{rueffel}) and the fact that the 
section  $s:\fl\rightarrow\fg$ defining $(\alpha,\gamma)\in \cZ(\fl,\fa)$ 
can be chosen in the following way. Decompose $\fh$ as a vector space 
into  $\fh=\fh_1\oplus\fh_2\oplus\fh_3$, where $\fh_1=\fh\cap\fri(\fg)$ 
and $\fh_2$ is a complement of $\fh_1$ in $\fh\cap\fri(\fg)^\perp$. 
Then choose $s:\fl\rightarrow\fg$ as described in Section \ref{S24} such 
that it satisfies  in addition $\fh_3\subset s(\fl)$ and $s(\fl)\perp \fh_2$.
\qed

Combining this with Theorem~\ref{class} we obtain a classification scheme for 
those Manin pairs $(\fg,\fh)$ for which $\fg$ does not have simple ideals.  
Let us remark that  now the absence of simple ideals is  a real restriction 
contrary to the case of metric Lie algebras without distinguished 
subalgebra.

In small dimensions Prop.~\ref{PManin}
is a helpful tool not only for classification of Manin pairs but also for classification of Manin triples, since a complementary isotropic subalgebra of 
$\fh\subset \fg$ (if it exists) can be determined by hand. E.g., using 
this method it is easy to recover the classification in dimension six 
without extensive calculations.

\section{Appendix: Some lemmas and proofs}

\subsection{Implications of $(\fh,K)$-equivariance}

We first consider the case of the trivial group $K=\{e\}$.

\begin{lm}\label{666}
Let $(\fl,\Phi_\fl)$ be an $\fh$-equivariant
Lie algebra such that $\fl^\fh\subset \fl'$. Let $\fr\subset\fl$ be the radical of the Lie algebra $\fl$.
Then $\fr$ is nilpotent and acts trivially on every semisimple $(\fl,\Phi_\fl)$-module.
\end{lm}
\proof 
Recall the notion of the nilpotent radical $R(\fg)$ of a Lie algebra $\fg$ from Section~\ref{S22}.
The ideal $R(\fg)$ is nilpotent and acts trivially on any semisimple $\fg$-module.
We consider the Lie algebra $\fh\ltimes\fl$, where the action of $\fh$ on $\fl$ is given by $\Phi_\fl$.
A semisimple $(\fl,\Phi_\fl)$-module can be considered as a semisimple $\fh\ltimes\fl$-module.
It is therefore sufficient to show that $\fr\subset R(\fh\ltimes\fl)$.

Formula (\ref{blub}) implies
that $\Phi_\fl(\fh)\fr+\fr\cap\fl'\subset R(\fh\ltimes\fl)$. 
Since $\fh$ acts semisimply on $\fl$ and $\fr\subset\fl$ is $\fh$-invariant we have $\fr=\fr^\fh\oplus \Phi_\fl(\fh)\fr$.
By assumption $\fr^\fh\subset\fr\cap\fl'$.
We conclude that $\fr\subset R(\fh\ltimes\fl)$. The lemma follows.
\qed  

We now look at the case $K=\ZZ_2$. Any $(\fh,\ZZ_2)$-equivariant Lie algebra has a decomposition
$\fl=\fl_+\oplus\fl_-$ w.r.t.~the $\ZZ_2$-action. 

\begin{lm}\label{667} Let $\fh$ be a semisimple Lie algebra. Let $(\fl,\Phi_\fl)$ be an $(\fh,\ZZ_2)$-equivariant
Lie algebra such that 
$\ 
(a)\,$ $\fl_+=\fl^\fh$, and 
$\ (b)\,$ $\fl_+=[\fl_-,\fl_-]$. 
Then $\fl$ is nilpotent.
\end{lm}
\proof Let $\fr$ be the radical of $\fl$. Then the semisimple Lie algebra $\fs:=\fl/\fr$ inherits an $(\fh,\ZZ_2)$-equivariant structure $\Phi_\fs$. Since any derivation of $\fs$ is inner we may identify
$\Phi_\fs(\fh)$ with a semisimple subalgebra $\fk\subset\fs$ acting on $\fs$ by the adjoint representation. $\Phi_\fs(\fh)$ respects the decomposition
$\fs=\fs_+\oplus\fs_-$. Thus $\fk\subset\fs_+$. By Condition (a) we have $\fs_+=\fs^\fk$.
Thus $\fk$ is abelian, hence zero. Applying (a) again we obtain $\fs=\fs_+$, thus $\fs_-=\{0\}$.
Now Condition (b) yields $\fs=0$. Thus $\fl$ is solvable. Now we consider $\fl$ as an $\fh$-equivariant solvable
Lie algebra and apply Lemma \ref{666}.
\qed

\subsection{Proof of Proposition \ref{inva}}
Any semisimple $(\fh,K)$-module $V$ has a canonical decomposition 
$V=V^{(\fh,K)}\oplus V^1$, where $V^1$ is the sum of of all 
irreducible submodules of $V$ carrying a non-trivial $(\fh,K)$-action.
We denote the corresponding components of any $v\in V$ by $v^0,v^1$.
Tensor products of semisimple $(\fh,K)$-modules are again semisimple. In particular, the natural
$(\fh,K)$-actions on $C^p(\fl,\fa)$, $C^q(\fl)$ are semisimple.

Let us first prove injectivity of the natural map $\HPhi\rightarrow \cH(\fl,\fa)^{(\fh,K)}$.
Let $(\alpha_i,\gamma_i)\in \cZ(\fl,\Phi_{\fl},\fa)=\cZ(\fl,\fa)^{(\fh,K)}$, $i=1,2$, be two cocycles
that represent the same element in $\cH(\fl,\fa)$, i.e., there exists 
$(\tau,\sigma)\in C^1(\fl,\fa)\oplus C^2(\fl)$ such that
$$ (\alpha_2,\gamma_2)=(\alpha_1+d\tau,\,\gamma_1+d\sigma+\langle(\alpha_1+\textstyle{\frac12} d\tau)\wedge\tau\rangle)
\ .$$
The invariance of $\alpha_i$ implies that $d\tau=(d\tau)^0=d\tau^0$. It follows that
$$(d\tau\wedge\tau)^0=(d\tau^0\wedge\tau)^0= d\tau^0\wedge\tau^0\ .$$
This and the invariance of $\gamma_i$ and $\alpha_1$ implies 
$$ \gamma_2=\gamma_{1}+d\sigma^0+\langle(\alpha_1+\textstyle{\frac12} d\tau^0)\wedge\tau^0\rangle\ .$$
Thus $(\alpha_2,\gamma_2)=(\alpha_1,\gamma_1)(\tau^0,\sigma^0)$. Since $(\tau^0,\sigma^0)\in 
{\cal C}^{1}_Q(\fl,\Phi_{\fl},\fa)$ the cocycles $(\alpha_i,\gamma_i)$, $i=1,2$, 
represent the same class in $\HPhi$.
This shows injectivity.

Surjectivity is a little bit more involved. Let $(\alpha,\gamma)\in \cZ(\fl,\fa)$ such that
$[\alpha,\gamma]\in \cH(\fl,\fa)^{(\fh,K)}$. This means that for all $X\in\fh$, $k\in K$ there exist elements
$(\tau_X,\sigma_X),(\tau_k,\sigma_k)\in C^1(\fl,\fa)\oplus C^2(\fl)$
such that
\begin{eqnarray}
k\alpha&=&\alpha+d\tau_k\, ,\qquad k\gamma\,=\,\gamma+d\sigma_k+\langle(\alpha+\textstyle{\frac12} d\tau_k)\wedge\tau_k\rangle\, ,\label{sand}\\
X\alpha&=&d\tau_X\, ,\qquad\quad\  
X\gamma\,=\,d\sigma_X+\langle\alpha\wedge\tau_X\rangle\, .\label{korn}
\end{eqnarray}
We have to show that $[\alpha,\gamma]$ can be represented by an $(\fh,K)$-invariant cocycle.

We first consider the case that $K$ is connected. Then $(\fh,K)$-invariance is the same as $\fk\ltimes\fh$-invariance,
where $\fk$ is the Lie algebra of $K$ which acts by derivations on $\fh$. Since the action of $\fk\ltimes\fh$
on $\fl$ and $\fa$ is semisimple it factors over a reductive quotient of $\fk\ltimes\fh$.
Thus without loss of generality we may assume that $K$ is trivial and $\fh=\fs\oplus\frak t$, where $\fs$
is semisimple and $\frak t$ is abelian. Let $C$ be the Casimir operator of $\fs$.
For any semisimple $\fh$-module $V$ there exists an element $X_V\in \frak t$
such that the operator $D_V:=X_V+C$ acts invertibly on any $\fh$-submodule
of $V^1$. In particular, $\im D_V=V^1$, $\ker D_V=V^\fh$.
Now we fix $D:=D_V$ for the semisimple $\fh$-module  $V=C^2(\fl,\fa)\oplus C^1(\fl,\fa) \oplus C^3(\fl) \oplus C^2(\fl)$.

The operator $D$ is of the form $D=X+\sum Y_iZ_i$ for some $X,Y_i,Z_i\in\fh$. 
Now (\ref{korn}) implies that there exists an element $\tau_D\in C^1(\fl,\fa)$ such that
$D\alpha=d\tau_D$. Thus $D\alpha\in B^2(\fl,\fa)\cap V^1$, where $B^2(\fl,\fa)\subset C^2(\fl,\fa)$ denotes
the $\fh$-submodule of coboundaries. Therefore there exists a coboundary $d\tau$ such that
$D\alpha=Dd\tau$. We now consider the quadratic cocycle $(\alpha',\gamma'):=(\alpha,\gamma)(-\tau,0)$.
Then $D\alpha'=0$, thus $\alpha'$ is $\fh$-invariant. The cocycle $(\alpha',\gamma')$ satisfies (\ref{korn}) again. We
claim that
\begin{equation}\label{trouble}
D\gamma'=d\sigma'_D+\langle\alpha'\wedge\tau'_D\rangle
\end{equation}
for some quadratic cochain $(\tau'_D,\sigma'_D)$ satisfying $d\tau'_D=0$.
It is obvious that (\ref{trouble}) is valid for $D$ replaced by $X\in \fh$. We have to show
that it also holds for $D$ replaced by a second order monomial $YZ$ 
for $Y,Z\in \fh$. Thus we can assume (\ref{trouble}) for $Z$
and apply $Y$. Note that $Y\alpha'=0$ by $\fh$-invariance. Setting
$\sigma'_{YZ}:=Y\sigma'_Z$ and $\tau'_{YZ}:=Y\tau'_Z$ we obtain
$$ YZ\gamma'=dY\sigma'_Z+Y\langle \alpha'\wedge \tau'_Z\rangle
=d\sigma'_{YZ}+\langle \alpha'\wedge \tau'_{YZ}\rangle\ .
$$
This justifies (\ref{trouble}).
We denote by $B^3(\fl)\subset C^3(\fl)$ and $Z^1(\fl,\fa)\subset C^1(\fl,\fa)$ the corresponding
submodules of coboundaries and cocycles. Then (\ref{trouble}) says
$$D\gamma'\in (B^3(\fl)+\langle\alpha'\wedge Z^1(\fl,\fa)\rangle)\cap V^1=:W^1\ .$$
Thus we can solve
$$ D\gamma'= D(d\sigma'+\langle\alpha'\wedge\tau'\rangle)$$
in $W^1$, in particular with $d\tau'=0$. Now we set $(\alpha'',\gamma''):=(\alpha',\gamma')(-\tau',-\sigma')$.
Note that $\alpha''=\alpha'$, thus $D\alpha''=0$. We have
$ D\gamma''= D\gamma'-Dd\sigma'-D\langle\alpha'\wedge\tau'\rangle=0$.
Thus $(\alpha'',\gamma'')$ is $\fh$-invariant and $[\alpha,\gamma]=[\alpha'',\gamma'']$. This implies
surjectivity of the natural map $\HPhi\rightarrow \cH(\fl,\fa)^{(\fh,K)}$ in the case of connected $K$.

If $K$ is disconnected and $[\alpha,\gamma]\in \cH(\fl,\fa)^{(\fh,K)}$ we can assume by the above
that $(\alpha,\gamma)\in \cZ(\fl,\fa)^{(\fh,K_0)}$, where $K_0\subset K$ is the identity component.
We now work in the vector space $V=[C^2(\fl,\fa)\oplus C^1(\fl,\fa) \oplus C^3(\fl) \oplus C^2(\fl)]^{(\fh,K_0)}$,
which carries an action of the finite group $G=K/K_0$. 
By injectivity of the canonical map $\cH(\fl,\Phi_\fl|_{(\fh,K_0)},\fa)\rightarrow \cH(\fl,\fa)^{(\fh,K_0)}$,
the
Equations (\ref{sand}) are valid with $k$ replaced
by $g\in G$ and with $\tau_g,\sigma_g$ invariant under $(\fh,K_0)$. We 
now proceed similar as in the connected case.
We will use the projection operator $P=\frac{1}{|G|}\sum_{g\in G} g$.
We set 
$\tau:={\frac{1}{|G|}}\sum_{g\in G}\tau_g$ and $(\alpha',\gamma'):=(\alpha,\gamma)(\tau,0)$.
Then we have $\alpha'=P\alpha$, i.e., $\alpha$ is $G$-invariant. 
Equation (\ref{sand}) for $(\alpha',\gamma')$
provides for $g\in G$ certain $(\fh,K_0)$-invariant elements $(\tau_g',\sigma_g')$ with $d\tau_g'=0$. 
We now define
$$\tau':= {\frac{1}{|G|}}\sum_{g\in G}\tau'_g\:,\ \qquad \sigma':=\frac{1}{|G|}\sum_{g\in G}\sigma'_g$$
and set $(\alpha'',\gamma''):=(\alpha',\gamma')(\tau',\sigma')$. Then $\alpha''=\alpha'$ is $G$-invariant. 
Using $d\tau'_g=0$ and (\ref{sand}) we eventually obtain
$$ \gamma''=\gamma' +d\sigma'+\langle \alpha'\wedge\tau'\rangle
=P\gamma'\ .$$
Thus $(\alpha'',\gamma'')$ is $G$-invariant. In other words, $(\alpha'',\gamma'')\in \cZ(\fl,\fa)^{(\fh,K)}$.
Since $[\alpha,\gamma]=[\alpha'',\gamma'']\in\cH(\fl,\fa)$ this finishes the proof of surjectivity
of the 
canonical map $\HPhi\rightarrow \cH(\fl,\fa)^{(\fh,K)}$.
\qed

\subsection{Proof of Proposition \ref{crux}}

We first assume Condition (b). Let $G$ be the transvection group of $M$, and let $J\subset G$ be the analytic subgroup corresponding to the ideal $\fri^\perp\subset\fg$. By the definition of a 
quadratic extension
(Definition \ref{qex}) the Lie algebra $\fri^\perp/\fri\cong\fa$ is abelian, i.e., 
$[\fri^\perp,\fri^\perp]\subset \fri$. Moreover, 
$$\langle [\fri^\perp,\fri],\fg \rangle =\langle \fri^\perp,[\fri,\fg] \rangle
\subset \langle \fri^\perp,\fri \rangle=\{0\}\ .$$
It follows that $\fri^\perp$ is 2-step nilpotent:
\begin{equation}\label{hampel} 
[\fri^\perp,[\fri^\perp,\fri^\perp]]=0\ .
\end{equation}
The following lemma together with Condition (b) implies that $J\subset G$ is closed.

\begin{lm}\label{hempel}
Let $G$ be a connected Lie group with Lie algebra $\fg$.
Let $\fj\subset \fg$ be a nilpotent ideal containing the center
$\fz(\fg)$. Then the analytic subgroup $J\subset G$ corresponding to $\fj$ is closed.
\end{lm}
\proof The assertion of the lemma is well-known, 
if $\fj\subset \fg$ is the maximal nilpotent ideal 
(see e.g. \cite{Bou3}, Ch.~III \S 9).
In this case $J$ is called the nilradical of $G$.
We look at the adjoint group $G_1:=\Ad_\fg(G)\cong G/Z(G)$ acting on
the Lie algebra $\fg$. Let $N_1$ be the nilradical of $G_1$. 
Then $N_1\subset G_1$ is closed.
By Engel's Theorem, there is a basis of $\fg$ such that all elements of $N_1$ are 
represented by upper
triangular matrices with $1$'s on the diagonal with respect to this basis. 
It follows that the exponential is a diffeomorphism from the Lie algebra $\fn_1$ of $N_1$ to $N_1$.
Now $\Ad_\fg(J)\subset N_1$. Via the exponential we
see that
$\Ad_\fg(J)\cong \ad_\fg(\fj)$ is closed in $N_1\cong\fn_1$. Hence $\Ad_\fg(J)$ is closed in $G_1$.  
It follows that $J\cdot Z(G)$ is closed in $G$.
Since $\fz(\fg)\subset\fj$ the group $J$ is the identity component of $J\cdot Z(G)$. Therefore, $J$ is closed as well.
\qed

We set $L:=G/J$. We identify the Lie algebra of $L$ via $p$ with $\fl$. Let $Q_L: G\rightarrow L$ be the natural
projection. Let $G_+\subset G$ be the stabilizer of the base point $x_0$, and let $\theta$ be the corresponding 
involution of $G$. Then $G_+\subset G^\theta$. 
Since $J$ is $\theta$-stable $\theta$ induces an involution $\theta_L: L\rightarrow L$. Then $L_+:=Q_L(G_+)$ is contained
in $L^{\theta_L}$ and has Lie algebra $\fl_+$. Since $L^{\theta_L}$ has 
the same Lie algebra we see that $L_+\subset L$ is closed. Then $N:=L/L_+$ is an affine symmetric space and the map
$q:M\cong G/G_+\rightarrow N$ induced by $Q_L$ is affine and surjective. The transvection group of $N$ is equal
to $L_0=L/(Z(L)\cap L_+)$. The natural map $Q$ between the transvection groups of $M$ and $N$ is therefore
the composition of $Q_L$ with the projection of $L$ to $L_0$. It is now evident that Condition (\ref{pickel}) is
satisfied.

The fibres of $q$ are connected since they are precisely the orbits of the
connected group $J$. Let us show that they are coisotropic and flat. By homogoneity it is sufficient to
look at tangent space $T_{x_0}M$ at the base point which can be identified with $\fg_-$. Then 
$(\fri^\perp)_-$
corresponds to the tangent space of the orbit. Since $\fri_-\subset (\fri^\perp)_-$ the orbits are coisotropic.
For $X,Y,Z\in \fg_-\cong T_{x_0}M$ the curvature tensor is given by
$R(X,Y)Z=-[[X,Y],Z]$. Flatness of the orbits now follows from (\ref{hampel}).

The uniqueness of $N$ and $q$ is a simple consequence of Condition (\ref{pickel}) and the 
required connectedness of the fibres.

If Condition (b) is not satisfied, but $M$ is simply connected, we make the following modifications of the above proof. First of all it is more convenient to work with the universal covering group of the transvection
group of $M$. We now denote this universal cover by $G$. The involution $\theta$ of $\fg$ induces one of $G$ which we 
again denote by $\theta$. $G$ acts on $M$. Let $G_+\subset G$ be the stabilizer of the base point $x_0\in M$.
The Lie algebra of $G_+$ is $\fg_+$.
Since $M$ and $G$ are simply connected $G_+$ is connected. It follows that $G_+\subset G^\theta$. 
(This is the crucial point where the simple connectedness of $M$ enters the proof. 
If $M$ were not simply connected we would
have $G^\theta
\subsetneq G_+$.)
Let $\tilde p: \fg\rightarrow \fl$ be the homomorphism
induced by $p$. Let $L$ be the connected and simply connected group with Lie algebra $\fl$. The surjection $\tilde p$
integrates to a surjection $Q_L: G\rightarrow L$. We set $J:=\ker Q_L$ which is closed in $G$. The group $J$
has Lie algebra $\fri^\perp$. Moreover, $J$ is connected since $L\cong G/J$ and $G$ are simply connected.
With these changes understood the proof runs as in the first case. Lemma \ref{hempel} is no longer needed.

\printindex
\end{document}